\documentclass[11pt,reqno]{amsart}
\usepackage[T1]{fontenc}
\usepackage{lmodern,amssymb,mathtools,microtype}
\usepackage[a4paper,margin=1.1in]{geometry}
\usepackage{xcolor}
\usepackage[pdfdisplaydoctitle]{hyperref}
\hypersetup{colorlinks,linkcolor={red!50!black},citecolor={blue!50!black},
  urlcolor={blue!75!black},
  pdftitle={Common crystal lattices and upper global bases for quantum cluster algebras},
  pdfauthor={Jiarui Fei and Tian Liu}}
  
\newtheorem{theorem}{Theorem}[section]
\newtheorem{proposition}[theorem]{Proposition}
\newtheorem{lemma}[theorem]{Lemma}
\newtheorem{corollary}[theorem]{Corollary}
\theoremstyle{definition}

\newtheorem{remark}[theorem]{Remark}
\numberwithin{equation}{section}
\allowdisplaybreaks

\newcommand{\Z}{\mathbb Z}
\newcommand{\Q}{\mathbb Q}
\newcommand{\N}{\mathbb N}

\title[Common crystal lattices and upper global bases]
{Common Crystal Lattices and Upper Global Bases\\ for Quantum Cluster Algebras}
\author{Jiarui Fei}
\address[Jiarui Fei]{School of Mathematical Sciences, Shanghai Jiao Tong University, China}
\email{jiarui@sjtu.edu.cn}
\thanks{The authors were supported in part by the National Natural Science Foundation of China (No.~12571038).}
\author{Tian Liu}
\address[Tian Liu]{School of Mathematical Sciences, Shanghai Jiao Tong University, China}
\email{tian.liu@sjtu.edu.cn}
\date{}

\begin{document}
\raggedbottom
\begin{abstract}
For a skew-symmetrizable quantum seed whose mutable part admits a reddening sequence, we construct a lattice common to all seeds. After the frozen variables are inverted, the integral bar-invariant lifts of its distinguished residues form the common triangular basis. If some frozen variables are kept polynomial, the regular lifts span precisely when every lift and the classical theta function of the same label have equal orders along those divisors.

We quantize the boundary derivations at reachable sinks and sources. In each nonzero divided power of a regular lift, the coefficient of the unique term of greatest string height is a quantum binomial of strictly greater Laurent degree than all other coefficients. For a nonempty set of sinks, let $F$ be the product of their frozen variables. The positive module embeds in copies of the restricted dual of the negative half, with regular lifts mapping to upper global vectors, precisely when all sufficiently long raising words kill $F$ and its cyclic submodule meets the torus on the unchosen frozen variables nontrivially. If one copy suffices, the algebra is a quantum Schubert cell up to a twist of multiplication.

For adapted integral coroots, an integrable quantum group extension exists precisely when each classical derivation extends to an integrable $\mathfrak{sl}_2$-action with its prescribed coroot, or equivalently when the coroot reflections preserve the regular theta cone. Its upper Kashiwara operators preserve the common lattice. The regular lifts give its upper global basis precisely when the two frozen orders agree for every label at the unchosen polynomial divisors of weight zero for all coroots. The degrees of the source monomials characterize existence and determine the lowering operators. A source can have multiplicity greater than one only if it is unchosen. When all frozen vertices are polynomial and chosen, spanning holds and every integrable extension comes from exact pairs. For the quantized coordinate ring of a simply connected complex semisimple group, a rescaling by the two weights recovers Kashiwara's upper global basis.
\end{abstract}
\maketitle
\setcounter{tocdepth}{1}
\tableofcontents

\section*{Introduction}
Let $G$ be a simply connected complex semisimple group. Its
quantized coordinate ring $\mathcal O_q(G)$ carries Kashiwara's
\emph{upper global basis} and has a quantum cluster structure
\cite{KashiwaraGlobal,OQY}. Can one recover the upper global basis
from products of cluster monomials? We approach this question
through the raising operators attached to frozen vertices of a
quantum upper cluster algebra. For $\mathcal O_q(G)$, the
construction recovers Kashiwara's basis after a rescaling determined
by the left and right weights.

Lusztig first constructed the \emph{canonical basis} of the positive
half of a quantized enveloping algebra in simply laced finite type.
He compared the Poincar\'e--Birkhoff--Witt bases associated with
different reduced expressions for the longest Weyl group element
\cite{LusztigCanonical}. Kashiwara constructed a \emph{crystal basis}
for the negative half and for each irreducible integrable highest
weight module in symmetrizable Kac--Moody type \cite{Kashiwara}.
The crystal operators preserve a lattice at $q=0$, where $q$ is
the quantum parameter, and act by partial inverses on a basis of its
reduction modulo $q$. The \emph{bar involution} sends $q$ to
$q^{-1}$. Each residue basis vector has a unique bar-invariant
lift belonging to both the lattice and the integral form. These lifts
form the \emph{lower global basis}. Under the usual
identification of the two halves, the lower global basis of the
negative half agrees with Lusztig's canonical basis
\cite{GrojnowskiLusztig,Lusztig}. Dualizing degree by degree gives
the upper global basis of the restricted dual.

In quantum cluster algebras, Berenstein and Zelevinsky used bar
invariance to construct a triangular basis from the standard monomials
of an acyclic seed. They showed that the basis is independent of
the acyclic seed when its mutable vertices are ordered compatibly
with the exchange matrix \cite[Theorems~1.4 and~1.6]{BZTri}.
Qin formulated triangularity in terms of multiplication by the
variables of a seed. A \emph{common triangular basis} is triangular
in every seed, and its degrees transform by tropical mutation. For
injective-reachable quantum cluster algebras, Qin gave sufficient
conditions for a triangular basis in one seed to be common to all
seeds. He also showed how existence is preserved under changes of
coefficients and quantization between similar seeds
\cite[Theorems~6.1.5--6.1.7]{QinTri}. For several families of
quantum unipotent subgroups he identified the common triangular
basis with the dual canonical basis after normalization and
localization at the frozen variables \cite[Theorem~1.2.1]{QinTri}.

For each Weyl group element in symmetric Kac--Moody type, Kang,
Kashiwara, Kim and Oh constructed a monoidal categorification of the
corresponding quantum unipotent coordinate ring using graded modules
over quiver Hecke algebras. They deduced that every quantum cluster
monomial is an upper global basis vector up to a power of $q^{1/2}$
\cite{KKKO}.
Qin later proved the cluster monomial theorem in symmetrizable type.
After inverting the frozen quantum minors, he identified the common
triangular basis of the resulting quantum unipotent cell with the
localized dual canonical basis rescaled by weight
\cite[Theorems~9.5.1 and~9.5.2]{QinDCB}. Qin subsequently used
freezing and changes of coefficients to construct common triangular
bases for other families in Lie theory, including decorated double
Bott--Samelson cells and quantum double Bruhat cells.
For symmetric Cartan matrices he also proved positivity
of their structure constants \cite[Theorems~1.1, 7.1 and~7.4]{QinSimilarity}.

The first author introduced crystal structures whose colors are reachable frozen
vertices. For skew-symmetric seeds he described the crystal operators
using representations of quivers with potential. The vertices of the
crystal are the presentation weights for which the general cokernel
is supported on the mutable vertices \cite{FeiCrystal}. He later
constructed boundary derivations of the upper cluster algebra directly
from a skew-symmetrizable seed. To describe one
such derivation, let $x_i$ be a frozen variable kept polynomial.
We call $i$ a \emph{reachable sink} if some seed has all arrows
between $i$ and the mutable vertices pointing towards $i$. In
that seed the derivation is $\partial/\partial x_i$ multiplied by
a monomial in the adjacent mutable variables. It is independent of
the choice of seed with $i$ a sink. At a \emph{reachable source}
the arrows point away from $i$, giving a second family of
derivations. Each family satisfies the Serre relations for its
boundary Cartan matrix \cite[Theorem~B]{Fei}. We call a seed in which
the chosen vertex is a sink or a source \emph{optimizing} for the
corresponding derivation.

Suppose the mutable exchange matrix admits a \emph{reddening sequence},
a sequence of mutations after which all mutable $c$-vectors for
principal coefficients are nonpositive. The regular theta functions
span a subalgebra called the \emph{middle algebra}. Berenstein and
Kazhdan introduced \emph{perfect bases} to read crystal arrows from
the leading terms of the raising operators
\cite[Definition~5.30]{BK}. The theta basis of the middle algebra is
\emph{normalized perfect} for both families of boundary derivations
\cite[Theorem~B]{Fei}. For a boundary derivation $\partial$, the $\partial$-degree of a
basis vector $b$ is the greatest $r$ for which $\partial^r b\ne0$. If $r>0$, the expansion of
$\partial b$ has a unique basis vector of degree $r-1$, with
coefficient $r$. All other terms have smaller degree. This vector
is the next vertex of the string, and it determines $b$ in turn.

In the quantum algebra we construct the lattice before choosing any
boundary operators. Limits of products of cluster monomials give it a distinguished
residue basis, and each residue in that basis has a unique integral
bar-invariant lift. Keeping some frozen variables polynomial singles
out the lifts that are regular along their divisors. Suppose these
lifts span the resulting integral form. The quantum sink and source
operators then give crystal arrows on their residues. If the positive
action extends to an integrable quantum group action with the prescribed
integral coroot weights, its upper global basis for these residues and
the given bar consists of the regular lifts.

\subsection{The common lattice and its distinguished lifts}
Let $U$ be the mutable vertex set and $V$ the full vertex set.
We use the row convention of \cite[\S1]{Fei} and write monomial
exponents as row vectors. The exchange at a mutable vertex $i$
is determined by the row $B_{i,V}$.
Write $B_\mu=B_{U,U}$ and fix a positive integral diagonal matrix
$D$ for which $B_\mu D$ is skew-symmetric. Suppose $B_\mu$ admits a
reddening sequence.
Let $v$ be the quantum parameter, $u=v^{-1}$, and
$R=\mathbb Z[v^{\pm1}]$.
To compare different choices of frozen variables, we begin in the
\emph{polarized seed}, whose matrix of mutable rows is
$(B_\mu\mid -I\mid I)$. Thus each mutable vertex has two frozen
vertices with opposite unit columns. For each seed $t$, reddening
gives a successor $t[1]$ whose mutable variables have degrees in
$t$ with mutable parts the negatives of the coordinate vectors,
up to permutation. We form normalized products of cluster monomials
in $t$ and $t[1]$, keeping their leading degree $m$ fixed.
The least choice of successor exponents gives Qin's injective
pointed element \cite[Definition~5.4.1]{QinTri}. We let the exponents
of all mutable successor variables tend to infinity. At each Laurent
exponent the coefficients have a
$u$-adic limit (Lemma~\ref{I:limit}). The resulting series
$S_m^t$ is \emph{pointed} at $m$, with leading monomial of
exponent $m$ and coefficient one. Every Laurent polynomial in
$t$ has a unique formal expansion in the $S_m^t$.

To take residues in these expansions, let $A_0$ be the ring obtained from
$\mathbb Z[u]$ by inverting the polynomials with constant term
$\pm1$. An element of the localized upper algebra over
$\mathbb Q(v)$ lies in the \emph{common lattice} exactly when all
coefficients in its expansion in one seed belong to $A_0$.
A \emph{tropical point} specifies a degree in every seed, with
adjacent degrees related by tropical mutation. For degrees
corresponding to the same tropical point, the coefficients at $u=0$
agree. The residue with coefficient one at a single tropical point
$g$ and zero elsewhere has a unique integral bar-invariant lift
$L_g$ (Theorem~\ref{I:localized}). The bar inverts $v$, fixes the
based monomials of a seed, and reverses products. Correction for bar
invariance commutes with mutation (Proposition~\ref{I:covariant}).
Reddening for the opposite matrix bounds the support of the corrected
series (Lemma~\ref{I:finite-support}). In every frozen pattern the exponents of a lift
have the form $m-aB_{U,V}$ with $a\in\mathbb N^U$. For each $a$ the
coefficient depends only on $B_\mu$, $D$ and the mutable part of $m$. For fixed mutable
exchange rows $B_{U,V}$, changing the compatible quantization
preserves the Laurent supports and the specializations $L_g(1)$
(Corollary~\ref{I:symmetrizer}).

The lattice and its residues thus determine the $L_g$ before we
consider triangularity. Comparing expansions in a seed and its
predecessor shows that the lifts form an $R$-basis of the quantum
upper algebra with all frozen variables inverted. Multiplication by each seed variable is
triangular on this basis. Thus under reddening $\{L_g\}$ is the
common triangular basis for every compatible quantization
(Theorem~\ref{I:qin-comparison}).

Let $P$ be the set of frozen vertices whose variables are kept
polynomial, and call the corresponding upper algebra the target.
We obtain its lattice by intersecting the
localized lattice with the algebra over $\mathbb Q(v)$. The order
of $L_g$ at a frozen divisor $f$ is its least Laurent exponent at
$f$. Thus $L_g$ is regular if and only if these orders are
nonnegative for every $f\in P$. Specializing $L_g$ at $v=1$
can only increase its order, and the resulting order is at most that
of the classical theta function $\vartheta_g$. The regular lifts
span over $R$ precisely when the quantum lift and the classical
theta function have the same order for every $g$ and $f\in P$
(Propositions~\ref{I:frozen-comparison} and~\ref{I:spanning-reductions}).

A term of displacement $a$ has $f$-exponent $m_f-aB_{U,f}$.
A sufficient condition for equality of the two orders at $f$ is that
the row $-B_{U,f}^{\mathsf T}$ lie in the closure of the support
of the cluster $g$-vector fan for the row matrix $B_\mu^{\mathsf T}$.
If this holds for all $f\in P$, the regular
lifts form a basis, and the frozen-order inequalities describe their labels as
the integral points of a rational polyhedral cone
\cite[Corollary~3.2]{Fei}. In particular, this holds when every
vertex in $P$ is a reachable sink or source
(Corollary~\ref{I:all-colors}). The regular lifts span the target
precisely when their residues span its common lattice modulo $u$,
and this last condition may be checked over $\mathbb Q$
(Corollary~\ref{I:balance}). In that case the integral form, the
lattice and its bar conjugate form a \emph{balanced triple}. In other words, reduction modulo $u$ maps their intersection
isomorphically onto the residue space.

\subsection{The quantum boundary action}
To define the boundary operators we also need the symmetrizing
integers at the frozen vertices. Fix positive integers $d_a$ such
that $D_V=\operatorname{diag}(d_a)_{a\in V}$ extends $D$ and
$BD_V$ is skew-symmetric. Choose reachable sinks with polynomial
variables as the \emph{colors} of the action and write $C=(c_{ij})$
for their boundary Cartan matrix. Their operators $E_i$ make the
target integral form a module over the divided-power integral form
of $U_v^+(C)$. Reachable sources give a second action with its own
Cartan matrix (Theorem~\ref{B:serre}).

In a seed optimizing a chosen sink $i$, write $x^{m_i}$ for the
image of $x_i$ under the classical derivation and set $s_i=m_i-e_i$,
where $e_i$ is the $i$th coordinate vector. On a based monomial $X^a$,
\[ E_i X^a=[a_i]_{v^{d_i}}X^{a+s_i}, \]
where $[r]_z=(z^r-z^{-r})/(z-z^{-1})$ is the symmetric quantum integer.
A centered commutator defines $E_i$ on
the whole algebra and specializes at $v=1$ to the classical
derivation. For a source we use the opposite algebra.

Let $L_g$ be regular at $i$. Its \emph{height} for $E_i$ is
the greatest $r$ for which $E_i^{(r)}L_g\ne0$. For $0\le j\le r$,
the expansion of $E_i^{(j)}L_g$ contains a unique lift of height
$r-j$, with coefficient
$\left[\begin{smallmatrix}r\\j\end{smallmatrix}\right]_{v^{d_i}}$.
Every other term has smaller height, and its coefficient has
$v$-degree strictly less than $d_i j(r-j)$, the degree of this
quantum binomial (Theorem~\ref{B:biperfect}). The bound follows by
comparing the terms at the degree and the codegree in two products,
one with the boundary monomial on the left and one with it on the
right (Lemma~\ref{B:product-bounds}). In particular, the last
nonzero divided power sends a regular lift to a single regular lift.

At $v=1$ the quantum upper algebra specializes to the classical upper
algebra over $\mathbb Z$ for every choice of polynomial frozen
variables (Proposition~\ref{B:fiber}). If the regular lifts span,
their specializations form a normalized perfect basis for both
families of derivations in characteristic zero
(Corollary~\ref{B:specialization}). On the coordinate ring of a maximal
unipotent subgroup, a \emph{biperfect basis} carries two crystal
structures coming from the left and right simple-root actions
\cite[Definition~2.2]{BKK}.

Our result settles the
triangular-basis assertion of \cite[Conjecture~8.14]{FeiCrystal} for the
skew-symmetric ice quivers with potential of
Corollary~\ref{B:presentation-biperfect}. More precisely, work over
$\mathbb C$, assume mutable reddening, a compatible quantization
and a finite-dimensional full Jacobian algebra, and let $\delta$
range over the \emph{$\mu$-supported} presentation weights, those
whose general cokernel is supported on the mutable vertices.
If the $L_{-\delta}(1)$ form a basis of the unlocalized upper
algebra, as assumed in the conjecture, this basis and the regular
theta basis have the same labels and are \emph{BK-biperfect} for
the reachable ordinary colors and the dually reachable starred
colors, with their compatible gradings, in the sense of
\cite[Definition~8.9]{FeiCrystal}. The theta-basis
assertion was proved in \cite[Corollary~4.9]{Fei}. The basis
hypothesis holds whenever every frozen vertex admits an optimizing seed.

Assume that the regular lifts span. For each color the strict degree
bound gives a basis of strings in the common lattice. Modulo $u$,
a step along a string sends a distinguished residue to the next one
or to zero. The reverse step induces the partial inverse (Theorem~\ref{B:strings}).
For a fixed color the strings are either all finite or all infinite.
In the optimizing seed, write $L_{-s_i}$ for the lift whose degree is
the negative of the shift. The strings are infinite precisely when this lift is
regular along every polynomial frozen divisor
(Corollary~\ref{B:string-lengths}).

\subsection{The positive half}
For the positive action we seek a $U_v^+$-module embedding into a
direct sum of copies of the restricted graded dual of the negative
half at $u=v^{-1}$, taking each regular lift to an upper global
vector. On this dual the positive generators act by transposes of
left multiplication by the negative generators \cite{Kimura}.

For a nonempty set of sink colors, let $F$ be the based product of
their variables with positive multiplicities and $\mathcal T$ the
quantum torus on the unchosen frozen variables in the localization.
Write $N(F)$ for the supremum of the lengths of words in the raising
operators that act nontrivially on $F$. By
Theorem~\ref{S:single-product} the embedding exists precisely when
\[ N(F)<\infty \qquad\text{and}\qquad (U_v^+F)\cap\mathcal T\ne0. \]
We can then choose it to induce the target integral form, bar and
common lattice. The regular lifts map bijectively onto the ambient
upper global vectors in the image. These vectors form a basis of
the image precisely when the regular lifts span the target integral form.

The first condition is equivalent to \emph{joint local nilpotence},
meaning that for each vector every sufficiently long word in the
raising operators kills it.
With $d$ the number of mutable vertices plus the number of chosen
sinks, $d+1$ repetitions of any ordering of the colors suffice to
test it (Theorem~\ref{S:frozen-product}). Let $A$ be the result of
successive maximal divided powers on $F$ over $d$ repetitions. The
second condition asks whether its degree is supported on the unchosen
frozen vertices. Both tests use tropical mutation and the valued
exchange rows in optimizing seeds
(Corollary~\ref{S:tropical-test}). If all frozen vertices are chosen,
the second condition is $1\in U_v^+F$.

One copy suffices precisely when the action is jointly locally
nilpotent and its joint kernel consists of the constants. Evaluation
then identifies the target, up to a twist of multiplication, with a
quantum Schubert cell whose flag minors are polynomial. The action
determines its Weyl group element, and the regular lifts give its
upper global basis (Theorem~\ref{C:cell}). For each reduced
expression, the variables of one extended cluster map to the chamber
minors (Proposition~\ref{C:chamber}). The comparison uses Kimura's
multiplication formula \cite[Theorem~6.19]{Kimura} and Qin's
triangularity theorem in the initial seed
\cite[Proposition~9.3.1 and Theorem~9.3.2]{QinDCB}.

More generally, under joint local nilpotence, invert the unchosen
frozen variables. Among the mutable variables in the joint kernel,
the smallest set to freeze and invert to obtain copies of one Schubert
cell consists of the mutable factors of $A$ (Theorem~\ref{S:minimal}).

In the embedding given by the frozen-product criterion, both ambient
upper Kashiwara operators of color $i$ preserve the original target
precisely when $L_{-s_i}$ is regular there. This lift is a monomial
in the successor of an optimizing seed, so its frozen exponents decide
regularity. The operators then preserve the integral form and common
lattice. If the regular lifts span the target integral form, the
restrictions agree with the distinguished string shifts. If $L_{-s_i}$ is regular for every $i$,
the Weyl group element is the longest element, and target spanning
is equivalent to the regular joint highest lifts spanning the joint
kernel over $\mathbb Q(v)$ (Theorem~\ref{S:shifts}).

\subsection{Integrable extensions}
When does the positive action determine lowering operators on the
target? The \emph{exchange grading} assigns to each monomial its
exponent modulo the integral span of the mutable rows of $B$. Let
$\Gamma$ be this quotient and $\gamma_j$ the degree of $E_j$.
Prescribe \emph{adapted integral weights}, homomorphisms
$\omega_i:\Gamma\to\mathbb Z$ with $\omega_i(\gamma_j)=c_{ij}$.
The positive action extends to an integrable type-one action of
the quantum group of the derived Kac--Moody algebra of $C$ over
$\mathbb Q(v)$ precisely when $E_i^{(r)}$ maps the space of
$i$-weight $-r$ isomorphically onto that of $i$-weight $r$ for
every $i$ and $r\ge1$. These isomorphisms determine the lowering
operators uniquely (Theorem~\ref{X:criterion}).

There are two equivalent tests for existence, both independent of
target spanning. The first asks whether, for each $i$, the reflection
$g\mapsto g-\omega_i(\deg X^g)s_i$ preserves the rational polyhedral
cone of regular theta labels in a seed optimizing $i$. It suffices to
test finitely many generators (Theorem~\ref{X:polyhedral}).
The cone also determines the only possible coroot grading. The
$i$-weight of each seed variable other than $X_i$ is its theta
string length, and the condition that $E_i$ have weight two fixes
the weight of $X_i$ (Corollary~\ref{X:reconstruct-weights}).

On each coroot weight space, injectivity and surjectivity of
$E_i^{(r)}$ over $\mathbb Q(v)$ can each be tested at $v=1$ over
$\mathbb Q$. Thus the quantum extension exists precisely when each
classical boundary derivation extends to an integrable
$\mathfrak{sl}_2$-action with its prescribed coroot
(Theorem~\ref{X:classical}). All simple divided powers then preserve
the target integral form, and the lowering operators specialize to
the classical ones. The boundary Cartan matrix has only finite and
affine components. Each affine component acts at level zero and
its raising operators fail joint local nilpotence
(Corollary~\ref{X:indefinite}).

In every integrable extension both upper Kashiwara operators
preserve the common lattice. An unchosen polynomial variable can
be inverted while preserving the action precisely when it has
weight zero for every chosen coroot (Corollary~\ref{X:localization}).
Inverting all such variables leaves only the chosen sinks and their
source partners polynomial. The regular lifts span in this localization,
whose stable lattice intersects the original target in its common lattice
(Theorem~\ref{X:partial-scale} and the proof following it).
In particular, spanning is automatic if each unchosen polynomial
variable has nonzero weight for some chosen coroot.

When the regular lifts span the original target integral form,
integrability is equivalent to the requirement that each string have
length equal to the coroot weight of its highest lift. In that case
the regular lifts form the upper global basis for the distinguished
residues, and the negative divided powers satisfy the strict
Laurent-degree bounds along the reversed strings
(Theorem~\ref{X:criterion}).

Given an integrable extension and at least one color, the
positive-half embedding exists precisely when successive maximal
raises on $F$ over $d$ periods give a Laurent monomial in the
unchosen frozen variables. Each component then has finite type and
some unchosen polynomial variable has positive weight for a coroot
in that component (Corollary~\ref{X:simultaneous}). If $h_i$ is the
grading in the product rule for $E_i$, the full action makes the
target a module algebra for the symmetric coproduct precisely when
$2h_i=-d_i\omega_i$ for every color (Proposition~\ref{X:product}).

A reachable sink $i$ and a reachable source $k$ form an
\emph{exact boundary pair} when the
source monomial at $k$ has pointed degree $e_k-s_i$ in a seed
optimizing $i$, where $e_k$ is the coordinate vector of $k$.
In the quiver-with-potential setting of \cite[Definition~4.28]{FeiCrystal},
this condition is equivalent to \emph{$\tau$-exactness} of the
corresponding boundary representations \cite[Proposition~4.6]{Fei}.
For such a pair the boundary derivations are
the two classical Chevalley operators \cite[Theorems~6.5 and~6.8]{Fei}.
Theorem~\ref{B:paired} constructs their quantum counterparts on every
upper algebra in which the paired variables are polynomial.
Conversely, if the centered operators of a specified sink and source
are the Chevalley operators of a type-one rank-one action, the pair is
exact and its commutator determines the coroot.
If the regular lifts span, the sink and source arrows on their
residues are the upper Kashiwara arrows (Corollary~\ref{B:upper}).
When all frozen variables are polynomial and covered by exact pairs,
an opposite involution of the valued pattern with frozen set $I$
induces a symmetry of the upper global basis via the relative
Donaldson--Thomas square. For a finite-type set of colors, the
frozen permutation of the given involution and finitely many
classical highest labels determine whether this symmetry acts by
Lusztig's involution (Corollary~\ref{B:relative-dt}).

An integrable extension recovers a unique reachable source $k$ for
each chosen sink $i$. After deleting the invertible frozen vertices,
its source monomial has $r_i=d_k/d_i$ times the degree required for
exactness. With adapted coroots, $r_i=1$ whenever $k$ is chosen.
If every polynomial vertex is chosen, the partners form a permutation
of exact pairs in the reduced seed, and the boundary strings recover
its coroots (Theorem~\ref{X:partial-scale} and
Corollary~\ref{X:partners}). Rescaling an unchosen source row gives
examples with $r_i>1$ (Section~\ref{X:recovered-lowering}).

These degree and multiplicity conditions also suffice for an extension
with some adapted integral coroots. The mutable tropical
Donaldson--Thomas transformation compares the normalized rows of
$i$ and $k$. The lift $S_i=L_{-s_i}$ must have order $-1$ at $k$ and
zero at every other polynomial vertex. Its graded
commutator gives $F_i$. On a seed variable the coroot weight is
the source height at $k$ minus the sink height at $i$
(Theorem~\ref{X:root-lowering}). For fixed partners the polynomial
patterns admitting the action form an interval, and the frozen orders
determine the subinterval on which the regular lifts give the upper
global basis (Corollary~\ref{X:restored-coefficients}).

For a connected, simply connected complex semisimple group $G$,
Qin conjectured that the common triangular basis of $\mathcal O_q(G)$
agrees with Kashiwara's upper global basis up to powers of $q^{1/2}$
\cite[Conjecture~1.6]{QinSimilarity}. For input weight $\delta$ and
output weight $\gamma$, multiplication by
$q^{((\delta,\delta)-(\gamma,\gamma))/4}$ takes the regular
distinguished lifts to Kashiwara's vectors (Theorem~\ref{O:comparison} and
Corollary~\ref{O:qin}). The principal
minor of each dominant weight is both a distinguished lift and the
normalized highest vector for the two regular actions. These common
vectors determine the comparison by Proposition~\ref{I:recognition}.

\section{Construction in the polarized seed}\label{I:construction}

Different boundary operators may be calculated in different seeds.
To compare their actions on residues, we need a lattice that can be
recognized in any seed. We first invert all frozen variables, so
that multiplying a cluster monomial by a frozen Laurent monomial
allows us to prescribe the frozen part of its degree. In
Section~\ref{I:polynomial} we recover the lattice for polynomial
frozen variables by imposing regularity at the corresponding divisors.

Let $V$ be the vertex set and list its $n$ mutable vertices $U$
first. For $B=(b_{ij})$ we use the row convention of
\cite[\S1]{Fei}. Write $E=B_{U,V}$
for the extended exchange matrix and $B_\mu=B_{U,U}$ for its mutable
part. Thus the row $B_{j,V}$ gives the exchange at $j\in U$.
Exponents and displacements are row vectors, and $e_j$ denotes the
$j$th coordinate row. Choose a positive diagonal symmetrizer
$D=\operatorname{diag}(d_j:j\in U)$ with integer entries, so that
$B_\mu D=-DB_\mu^{\mathsf T}$. An integral skew form $\Lambda$ on exponents
is \emph{compatible} if $E\Lambda=-(D\mid 0)$, where
$\Lambda(a,b)=a\Lambda b^{\mathsf T}$. In particular the rows of
$E$ are linearly independent. In the convention of
\cite[\S3]{BZ} the extended matrix is $-E^{\mathsf T}$.
The quantum torus of a seed $t$ has an
$R=\mathbb Z[v^{\pm1}]$-basis of monomials $X_t^a$ indexed by
$a\in\mathbb Z^V$. The variable at $j$ is $X_j=X_t^{e_j}$. Suppressing $t$, the product is
$X^aX^b=v^{\Lambda(a,b)}X^{a+b}$. The \emph{bar involution}
inverts $v$, fixes the based monomials $X^a$, and reverses products
\cite[\S\S3--4]{BZ}.
A \emph{localized cluster monomial} is a based monomial in some seed
whose mutable exponents are nonnegative. Its frozen exponents may be
arbitrary integers.

For the lattice and its reduction we use $u=v^{-1}$ and
$k=\mathbb Q(v)$. Let $A_0$ be the localization of $\mathbb Z[u]$
obtained by inverting the polynomials with constant coefficient
$\pm1$. Evaluation at zero identifies $A_0/uA_0$ with $\mathbb Z$.

The mutable exchange rows define the \emph{dominance order} on exponents
by $m-rE\preceq_t m$ for $r\in\mathbb N^n$, where $0\in\mathbb N$.
We omit the subscript when the seed is understood, and write $\prec_t$
for the strict order.
A Laurent polynomial is \emph{pointed} at $m$ if it has the form
$X^m+\sum_{0\ne r\in\mathbb N^n}c(r)X^{m-rE}$. The largest
exponent $m$ is its \emph{degree}. Degrees in different seeds
represent the same \emph{tropical point} if tropical mutation takes one
to the other. In row coordinates, mutation at $j$ sends $m_j$ to $-m_j$ and
$m_a$ to $m_a+[b_{ja}]_+m_j-b_{ja}[m_j]_+$ for $a\ne j$,
where $[z]_+=\max(z,0)$ \cite[\S3.3]{Fei}. The fully
localized quantum upper cluster algebra $\mathcal U_R^{\mathrm{loc}}$
is the intersection of the seed tori in their common division ring.
An element of this algebra is \emph{compatibly pointed} if in every
seed it is pointed at the representative of a single tropical point.

The exponents of a pointed element all have the same class in
$\Gamma=\mathbb Z^V/\mathbb Z^U E$. This class is its
degree for the \emph{exchange grading}. To see that the grading
extends to every seed, let $\lambda=(\lambda_a)_{a\in V}$ be the
column of degrees in a seed $t$. If $E_t\lambda=0$, the two
exchange monomials at $j$ have the same degree. Giving the new
variable degree $-\lambda_j+[B_{j,V}]_+\lambda$ yields
$E_{\mu_jt}\lambda'=0$ by matrix mutation.

For a compatibly pointed element, Lemma~\ref{I:predecessor} recovers
the least exponent from its degree after a suitable change of seed.
An element is \emph{copointed}
at its \emph{codegree} $\eta$ if
$X^\eta$ has coefficient one and all its other exponents lie in
$\eta+(\mathbb N^n\setminus\{0\})E$.
When multiplying two elements $z,z'$ pointed at $m,m'$, we divide by the
coefficient of $X^{m+m'}$ so the product stays pointed. We denote
this \emph{normalized product} by $[zz']_t=v^{-\Lambda(m,m')}zz'$.

To compare frozen patterns with the same mutable matrix, we first
adjoin a pair of frozen vertices for each mutable vertex. In this
\emph{polarization} the mutable exchange rows form the matrix
$E_{\mathrm{pol}}=(B_\mu\mid -I_n\mid I_n)$.
In the indicated decomposition a compatible form is given by
$\Lambda_{12}=-D$, $\Lambda_{21}=D$, and
$\Lambda_{22}=-B_\mu D$, with all other components zero.

To construct an element of each degree we use cluster monomials
pointed in the negative mutable directions. The principal extension
has mutable rows $(B_\mu\mid-I_n)$. After mutation, write its mutable
rows as $(B_\mu(t)\mid-C_t^{\mathsf T})$. In the usual notation the
columns of $C_t$ are the $c$-vectors, and $C_t$ is initially
the identity. We assume that
$B_\mu$ admits a \emph{reddening sequence}, a sequence of mutations
after which all these vectors are nonpositive.
Equivalently, the seed is \emph{injective-reachable}. Every seed $t$
then has a \emph{successor} $t[1]$ whose mutable cluster variables
have mutable degrees $-e_1,\ldots,-e_n$ in $t$, up to permutation
\cite[\S2.3, especially Corollary~2.3.4]{QinTropical}.

Fix a seed $t$ in any compatible frozen pattern. Given
$N\in\mathbb N^n$, take the cluster monomial in the mutable variables
of $t[1]$ of mutable degree $-N$ in $t$. Multiplication by a
frozen monomial and normalization give an \emph{injective cluster
monomial} $I_N$ of full degree $(-N,0)$. To keep the full degree
fixed at $m=(g,c)$ as $N$ varies, take the normalized products
$P_{m,N}=[X_t^{(g+N,c)}I_N]_t$. For $N\ge[-g]_+$ componentwise
both factors are localized cluster monomials. Write $C_N(r)$ for
the coefficient of $X_t^{(-N,0)-rE}$ in $I_N$, and set
$\langle g,r\rangle=gDr^{\mathsf T}$. Compatibility gives the
coefficient $u^{\langle N+g,r\rangle}C_N(r)$ at $m-rE$ in
$P_{m,N}$. We take the limit at each displacement and define
\[ S_m^t=\sum_{r\in\mathbb N^n}\left(\lim_{\min_i N_i\to\infty}u^{\langle N+g,r\rangle}C_N(r)\right)X_t^{m-rE}. \]
At each displacement $r$, Lemma~\ref{I:limit} shows that the
$u$-adic limit exists and is a rational function of $u$. If $n=0$, we take
$S_m^t=X_t^m$.

The coefficient at $m$ in $S_m^t$ is one. Every Laurent polynomial
$x$ in this chart over $k$ therefore has a unique formal
expansion $x=\sum_m c_m^t(x)S_m^t$ with indices in
$\operatorname{supp}_t(x)-\mathbb N^n E$, where
$\operatorname{supp}_t(x)$ is the set of exponents occurring in $x$.
Indeed, the Laurent support is finite and the rows of $E$ are independent,
so at each exponent only finitely many larger degrees can contribute.
Subtracting their contributions determines the next coefficient.
We require all these coefficients to belong to $A_0$ and define
\[ \mathcal L_t=\{x\in k\otimes_R\mathcal U_R^{\mathrm{loc}}:c_m^t(x)\in A_0\text{ for all }m\}. \]
For $x\in\mathcal L_t$, evaluation at zero gives its residue
coordinates $c_m^t(x)\bmod uA_0$. These coordinates are indexed by
degrees in $t$. The theorem identifies the same residues in every
seed by indexing them by tropical points. Here a lattice is the
$A_0$-span of a $k$-basis.
An element of the localized upper
algebra over $k$ is \emph{integral} if it belongs to
$\mathcal U_R^{\mathrm{loc}}$.

\begin{theorem}\label{I:localized}
Suppose that the skew-symmetrizable mutable matrix $B_\mu$ admits a
reddening sequence and has symmetrizer $D$. For every compatible
frozen pattern with mutable data $(B_\mu,D)$, the module $\mathcal L_t$
is a lattice independent of $t$. Its residue coordinates identify its reduction modulo $u$
with the free abelian group on tropical points.
The unit vector at each tropical point has a unique integral
bar-invariant lift in this lattice. In a chart where its degree is
$m$, denote the lift by $L_m$. It has the expansion
\[ L_m=\sum_{r\in\mathbb N^n} B_{m_{\rm mut}}(r)X^{m-rE} \]
with finite support. We also write $L_\gamma$ for the lift indexed
by the tropical point $\gamma$ represented by $m$. The coefficients
are obtained by imposing bar invariance on the formal expansion
in the polarized seed and are the same in all compatible patterns
with the given mutable matrix and symmetrizer. The $L_m$ are
compatibly pointed and form an $A_0$-basis of the lattice and
an $R$-basis of $\mathcal U_R^{\mathrm{loc}}$.
\end{theorem}

The bar correction first determines the $B_h(r)$ as coefficients
of a formal series in a single seed. We show that the corrected series agree under
mutation. We then use a reddening sequence to bound their supports,
so that they define elements of the upper algebra. When $n=0$
the based Laurent monomials themselves give both bases. Henceforth
in the construction assume $n>0$.

\subsection{Degrees in a seed and its predecessor}

We use the quantum Laurent phenomenon \cite[\S5]{BZ}. Quantum
cluster variables are compatibly pointed by the quantum separation
formula \cite[Theorems~5.3 and~6.1]{Tran}, sign coherence
\cite[Corollary~5.5]{GHKK}, and tropical mutation of $g$-vectors
\cite[Proposition~4.2]{NZ}.
The same holds for localized cluster monomials by based multiplication.
To recognize a cluster monomial from its tropical point, we recover
its least exponent in $t$ from its degree in a
\emph{predecessor}, a seed $t[-1]$ whose successor is $t$ up
to relabeling. Write $\mathbb Z^V(t)$ for the exponent lattice with
coordinates in seed $t$.

\begin{lemma}\label{I:predecessor}
Let $s=t[-1]$ be the predecessor of $t$, and let $\psi:\mathbb Z^V(t)\to\mathbb Z^V(s)$
send the coordinate vector of each variable of $t$ to its degree
in $s$. This linear isomorphism reverses dominance.
If $z\ne0$ is Laurent in both seeds, $\psi$ takes the minimal
exponents of its expansion in $t$ to the maximal exponents of its
expansion in $s$, preserving their coefficients. In particular,
a compatibly pointed element with the same tropical point as a
localized cluster monomial equals that monomial.
\end{lemma}
\begin{proof}
The degree--codegree correspondence and its consequence for cluster
monomials are due to Qin
\cite[Propositions~3.3.11--3.3.12 and Lemma~3.4.12]{QinTropical}.
Write $G$ for the matrix of $\psi$, so $\psi(a)=aG$. Its mutable
part is $-P$ for a permutation matrix $P$, and frozen variables
are fixed. Thus $G$ is unimodular. The rows of $E_tG$ vanish in
the exchange grading of $s$, so $E_tG=H E_s$ for some integral
matrix $H$. Comparing leading monomials in the commutation relations
gives $G\Lambda_sG^{\mathsf T}=\Lambda_t$. Multiplying
$E_tG=HE_s$ on the right by $\Lambda_sG^{\mathsf T}$ and using
compatibility yields $HDP^{\mathsf T}=-D$. Thus $-H$ has one
positive integer in each row and column. Their product is one, so
$H=-P$. Hence $\psi$ reverses dominance.

Each variable of $t$ is pointed in $s$. Its inverse has a formal
expansion obtained, after factoring out the leading monomial, by a
geometric series in positive displacements. Products of these series are
defined coefficientwise since $E_s$ has independent rows. Since
$\psi$ preserves the skew forms, every based
monomial $X_t^a$ has a formal expansion pointed at $\psi(a)$.
Substituting these expansions in $z$ gives its Laurent expansion
in $s$, as we see by clearing the $t$-denominators with a single
monomial. The maximal exponents among the images of the support are
distinct and cannot cancel. Since $\psi$ reverses dominance, these
are exactly the images of the minimal exponents in $t$, with their
original coefficients.

For the last assertion, choose a seed $t$ in which the cluster
monomial is $X_t^a$. Its degree in the predecessor is $\psi(a)$.
A compatibly pointed element with the same tropical point has these
same two degrees. Its largest and smallest exponents in $t$ are
therefore both $a$, and the coefficient there is one. Hence the
element is $X_t^a$.
\end{proof}

\subsection{The limit of injective cluster monomials}

Fix a seed $t$ of any compatible frozen pattern. The normalized
product of $I_{N'}$ and $I_N$ is $I_{N+N'}$, since all three are
monomials in the successor seed and their pointed degrees add in $t$.
This identity gives the limit in the definition of $S_m^t$. The
polynomial $C_N=\sum_r C_N(r)y^r$ records the coefficients of
$I_N$. Here $y^r$
stands for displacement $r$ from the pointed degree. Translating that
degree introduces powers of $u$, which we record by
$\rho_g(y^r)=u^{\langle g,r\rangle}y^r$ for $g\in\mathbb Z^n$.
The product of displacement terms is
$y^r*y^s=u^{-\langle rB_\mu,s\rangle}y^{r+s}$. It extends to formal series
since every nonnegative displacement has only finitely many decompositions
as a sum of two nonnegative vectors.

Write $G_N=\rho_N C_N$ for the rescaled coefficient series of $I_N$.
Suppose two pointed factors have mutable degrees
$g,h$. By compatibility, the terms of displacements $r,s$ in
their normalized product acquire the power of $u$ with exponent
$\langle g,s\rangle-\langle h,r\rangle-\langle rB_\mu,s\rangle$.
Applying this to $I_{N'} I_N$ gives
\begin{equation}\label{I:product} G_{N+N'}=\rho_{2N}(G_{N'})*G_N . \end{equation}
\begin{lemma}\label{I:limit}
For each displacement $r$, the coefficient $G_N(r)$ has a limit
in the ring $\mathbb Z((u))$ of formal Laurent series as
$\min_i N_i\to\infty$. The limiting series
$G$ is the unique series with constant term one satisfying
\[ G=\rho_{2\mathbf1}(G)*G_{\mathbf1},\qquad \mathbf1=(1,\ldots,1). \]
Its coefficients belong to the subring of $\mathbb Q(u)$ generated
by $R$ and $(1-u^{2a})^{-1}$ for $a>0$.
\end{lemma}
\begin{proof}
For $r\ne0$, the fixed point equation expresses
$(1-u^{2\langle\mathbf1,r\rangle})G(r)$ in terms of coefficients of
smaller total displacement. Since $\langle\mathbf1,r\rangle>0$,
recursion gives a unique solution with constant term one over the
stated ring. It remains to prove convergence to this solution.

On formal series let $F_N(A)=\rho_{2N}(A)*G_N$. By
\eqref{I:product}, $F_NF_{N'}=F_{N+N'}$. Thus $F_N(G)$ has
constant term one and is fixed by $F_{\mathbf1}$, giving $F_N(G)=G$.
Write $G^{*-1}$ for the inverse in the formal series algebra. We have
\begin{equation}\label{I:finite-approx} G_N=\rho_{2N}(G^{*-1})*G. \end{equation}
At a fixed displacement $r$ only finitely many displacements $a\le r$
of the first factor can contribute. For each $a\ne0$, the factor
$u^{2\langle N,a\rangle}$ tends to zero $u$-adically as every
$N_i$ tends to infinity. The contribution from $a=0$ is $G(r)$.
\end{proof}

\subsection{Formal expansions and the bar correction}

We now work in the polarized pattern. To lift the residue indexed by
$m$, we seek bar-invariant pointed series in which the expansion
of $S_m^t$ has coefficient one at
$m$ and all lower coefficients vanish at $u=0$.
The displacement series of $S_m^t$ is $S_g=\rho_gG$, where
$g=m_{\rm mut}$.

For a term based at $m-aE$, the translation
$\tau_a(y^r)=y^{a+r}$ expresses its displacements relative to $m$.
The coefficients of $S_g$ belong to $A_0[u^{-1}]$ by
Lemma~\ref{I:limit}. This ring has the additive decomposition
$A_0[u^{-1}]=R^{\mathrm{bar}}\oplus uA_0$, where
$R^{\mathrm{bar}}=\mathbb Z[v+v^{-1}]$. Indeed, the constant and
negative powers in the Laurent expansion at zero determine a unique
bar-invariant Laurent polynomial. Subtracting it leaves a rational
function of positive order whose denominator, after cancellation of
powers of $u$, has constant term $\pm1$.

There is a unique family of formal series $B_g$ with constant term
one and coefficients in $R^{\mathrm{bar}}$ such that
\begin{equation}\label{I:correction} S_g=B_g+\sum_{a\ne0}z_{g,a}\tau_a B_{g-aB_\mu},\qquad z_{g,a}\in uA_0. \end{equation}
For $r\ne0$, subtract from
$S_g(r)$ the terms $z_{g,a}B_{g-aB_\mu}(r-a)$ with $0\le a\le r$,
$a\notin\{0,r\}$. Split the result into its two components. They are
$B_g(r)$ and $z_{g,r}$, respectively. Induction on $\sum_i r_i$
gives existence and uniqueness simultaneously for every $g$.
The same induction applies when the series are indexed by full
degrees, with the translate at $a$ based at $m-aE$.
The separation of the bar-invariant part at each degree is the
recursion underlying Lusztig's lemma \cite[Theorem~1.1]{BZTri}.
Here the initial series is the limit of injective cluster monomials.
We prove mutation invariance and finiteness of the corrected series
in Section~\ref{I:mutation}.

The correction is unchanged by a triangular change of series whose
off-diagonal coefficients vanish at zero. Indeed, suppose that
$S_g=T_g+\sum_{a\ne0}t_{g,a}\tau_aT_{g-aB_\mu}$ with
$t_{g,a}\in uA_0$. Each displacement has only finitely
many ordered decompositions into nonzero displacements. The inverse
change thus has the same triangular form, and uniqueness in
\eqref{I:correction} gives the same $B_g$ starting from $T_g$.

\section{Mutation invariance of the common lattice}\label{I:mutation}

\subsection{Mutation in one variable}

To compare pointed expansions in adjacent seeds, we group terms that
differ by a multiple of the exchange row at the mutated vertex. Mutation acts
on each group separately. Fix $t'=\mu_j t$ and a degree $m$ in
$t$ with mutable part $g$. Each group of displacements is a line
$\nu+\mathbb N e_j$, where $\nu\in\mathbb N^n$ and $\nu_j=0$.
We call $\nu$ the \emph{transverse displacement}. On this line
we collect the coefficients into a polynomial or formal series in
an indeterminate $x$, assigning $x^a$ to $\nu+ae_j$. We use
the indeterminate $y$ for the same purpose in $t'$.

The exchange relation and its powers determine the change of variable
on each line. Write $w=u^{d_j}$. For
$l\ge0$ let $Q_l(x)=\prod_{b=0}^{l-1}(1+w^{l-1-2b}x)$, with
$Q_{-l}=Q_l^{-1}$ for $l>0$. The common exponent
of the old cluster variable is $h$, and the shift from the new pointed
degree is $A$. In terms of $g$ and $\nu$, they are given by
\[ s=g_j,\quad c=-\sum_{i\ne j}\nu_ib_{ij},\quad \lambda=\sum_{i\ne j}\nu_i[b_{ij}]_+,\quad h=s+c,\quad A=[-s]_++\lambda. \]
Here $s$ depends only on the pointed degree, while $c$ and
$\lambda$ depend only on the transverse displacement.

\begin{lemma}\label{I:fibre}
Measured from the tropical degrees in the two seeds, quantum mutation
takes a polynomial $f(x)$ at transverse displacement $\nu$ to
$y^A Q_h(y) f(y^{-1})$. Over $R$, or over $A_0[u^{-1}]$, the polynomials
with nonnegative displacements on both sides have basis
\[ x^a Q_{[-h]_+}(x),\qquad 0\le a\le A-[-h]_+ . \]
Mutation sends the member indexed by $a$ to
$y^{A-[-h]_+-a}Q_{[h]_+}(y)$.
\end{lemma}
\begin{proof}
Write $E'$ for the exchange matrix of $t'$, and let
$p_i=-e_iE$, $p'_i=-e_iE'$ be the displacements of the $i$th
exchange monomials in the two seeds. First take the linear branch of
mutation that sends $e_j$ to $-e'_j+[p_j]_+$ and fixes the
coordinates at the other vertices. Here positive parts are taken
coordinatewise. This map sends $p_j$ to $-p'_j$ and $p_i$ to
$p'_i+[b_{ij}]_+p'_j$ for $i\ne j$. The tropical image of $m$
is its linear image plus $\min(s,0)p'_j$. Thus reversal of the line and translation of its origin give
$y^A f(y^{-1})$. The remaining factor comes from the exchange
relation, which contributes $1+y$ for exponent one. Compatibility makes the
coordinates other than $j$ pair trivially with $p'_j$
under $\Lambda'$. Thus the factor depends only on the old
$j$th exponent, and the based product for
exponents $a,b$ gives $Q_{a+b}(y)=Q_a(w^by)Q_b(w^{-a}y)$.
This also gives the formula for negative exponents by division
\cite[\S\S4--5]{BZ}. Since $Q_l$ has both extreme coefficients one
and $x^lQ_l(x^{-1})=Q_l(x)$ for $l\ge0$, polynomiality on both
sides is equivalent to $\deg f\le A$ and $Q_{[-h]_+}\mid f$.
These conditions give the stated basis and its image. Its index range
contains zero because $A\ge0$ and
$A+h=[s]_++\sum_{i\ne j}\nu_i[-b_{ij}]_+\ge0$.
\end{proof}

For each full degree $\ell$ in $t$, write $\ell'$ for its
tropical image in $t'$. Take $M_\ell=X_t^\ell$ if $\ell_j\ge0$
and $M_\ell=X_{t'}^{\ell'}$ if $\ell_j<0$. In $t$ its only
displacements from $\ell$ lie on the $j$-axis, where its
polynomial is $Q_{[-\ell_j]_+}(x)$. In $t'$ its polynomial on
the corresponding axis is $Q_{[\ell_j]_+}(y)$, by
Lemma~\ref{I:fibre}. In particular $M_\ell$ is bar invariant.

These $M_\ell$ are the normalized standard monomials for mutation at
$j$ with all the other variables inverted
\cite[Definition~7.2]{BZ}. They form an $R$-basis of the
intersection of the two seed tori. Indeed, group the exponents by
cosets of the subgroup generated by the exchange row at $j$.
The linear map in the proof of Lemma~\ref{I:fibre} identifies the
cosets in the two seeds. Clearing a common denominator in the line
variable preserves each coset, so Laurentness can be checked on one
coset at a time. There the old $j$th exponent $h$ is constant.
The mutation formula takes $f(x)\in R[x^{\pm1}]$ to a Laurent
polynomial precisely when $Q_{[-h]_+}(x)$ divides $f(x)$ in
this ring. The basis $x^a Q_{[-h]_+}(x)$, $a\in\mathbb Z$,
consists of the $M_\ell$ on this coset.

We write $\ell\preceq_{t,t'}m$ when $\ell\preceq_t m$ and
$\ell'\preceq_{t'}m'$, and use $\prec_{t,t'}$ for the strict
order. Every interval is finite, since an interval for
$\preceq_t$ is an integral coordinate box in the displacement
variables. For $\ell=m-(\nu+ae_j)E$, the polynomial of
$M_\ell$ on the line at $\nu$, measured from $m$, is
$x^aQ_{[-h]_+}(x)$. Its least power on the primed line is
$A-[-h]_+-a$. Thus $\ell\preceq_{t,t'}m$ precisely for the
indices $0\le a\le A-[-h]_+$ in Lemma~\ref{I:fibre}.
In particular a pair of corresponding polynomials with nonnegative
displacements in the two charts has the same coefficients in the
$M_\ell$ on each line.

\subsection{Compatibility of the bar correction with mutation}

To use the polynomial bases above, we replace $S_m^t$ by a series
that is polynomial on each line in both charts. The change affects
one factor, which comes from displacements on the $j$-axis.
We will see that it leaves the bar correction unchanged.

Let $\Psi$ be the restriction of $G$ to the $j$-axis.
Restriction respects multiplication because all displacements are
nonnegative. In \eqref{I:finite-approx}, take $N_j=0$ and let
the other $N_i$ tend to infinity. Only displacements on the
$j$-axis survive in the first factor, so
\begin{equation}\label{I:partial} H=\lim_{\substack{N_j=0\\N_i\to\infty\ (i\ne j)}}G_N=\Psi^{*-1}*G. \end{equation}
For $n=1$ we take $N=0$, giving $H=1$. Each coefficient of
$H$ belongs to $A_0[u^{-1}]$, since computing a coefficient of
the formal inverse requires only finitely many coefficients of $\Psi$.
Regard $H$ as the series $H^t=\sum_r H(r)X_t^{-rE}$ of degree
zero, and set $T_m^t=M_m H^t$.

\begin{lemma}\label{I:relative-mutation}
On each line of transverse displacement $\nu$, the series
$T_m^t$ is a polynomial, and mutation takes it to the corresponding
polynomial of $T_{m'}^{t'}$.
\end{lemma}
\begin{proof}
When $N_j=0$, mutation fixes $(-N,0)$. The injective monomials
of that degree in the two charts have the same tropical point and
coincide by Lemma~\ref{I:predecessor}. The monomial $X_t^{(N,0)}$
uses only the shared variables and has degree opposite to that of
$I_N$, so the normalizing scalar is one. Thus
$P_{0,N}=X_t^{(N,0)}I_N$ is the same element in the two charts,
of degree zero in each. Its displacement coefficients in $t$ are $G_N$.

The product $M_mP_{0,N}$ is pointed at $m$ in $t$ and at
$m'$ in $t'$. Its coefficients converge to those of $T_m^t$
and $T_{m'}^{t'}$ in the respective charts, because $M_m$ is
Laurent in both. Lemma~\ref{I:fibre} bounds the degree on each line
in either chart independently of $N$. We can therefore pass to the
limit in its polynomial mutation formula, after clearing the fixed
denominator $Q_{[-h]_+}(y)$ when $h<0$.
\end{proof}

To recover the contribution of the axis, write
$\Psi(x)=\sum_{a\ge0}\psi_a x^a$. The axis of $I_{e_j}$ has
$s=h=-1$ and $A=1$, so its polynomial is $Q_1(x)=1+x$
by Lemma~\ref{I:fibre} and pointing. The corresponding polynomial
of $G_{e_j}$ is $1+wx$. On the axis the product $*$ is ordinary
multiplication, and $G=\rho_{2e_j}(G)*G_{e_j}$ gives
$\Psi(x)=(1+wx)\Psi(w^2x)$. Thus $\psi_0=1$ and
$(1-w^{2a})\psi_a=w^{2a-1}\psi_{a-1}$, whence
\[ \psi_a=\frac{w^{a^2}}{\prod_{b=1}^a(1-w^{2b})}\qquad(a>0). \]

Let $W_m^t=\sum_{a\ge0}w^{sa}\psi_a X_t^{m-ae_jE}$, where $s=m_j$.
Then $S_m^t=W_m^t H^t$. Indeed, multiplying terms of displacements
$ae_j$ and $b$ contributes the factor
$u^{\langle g,b\rangle-\langle ae_jB_\mu,b\rangle}$. Together with
the factor $w^{sa}$ in $W_m^t$, this gives the coefficient of
$\rho_g(\Psi*H)=\rho_g G$ by \eqref{I:partial}.
Iterating the functional equation for $\Psi$ gives
$\Psi(w^sx)=Q_{[-s]_+}(x)\Psi(w^{|s|}x)$.
The $j$th coordinate of $m-ae_jE$ is still $s$, so the
polynomial of $M_{m-ae_jE}$, measured from $m$, is
$x^aQ_{[-s]_+}(x)$. Expanding $W_m^t$ in these monomials and
multiplying on the right by $H^t$ gives
\begin{equation}\label{I:relative} S_m^t=T_m^t+\sum_{a>0}w^{a|s|}\psi_a\,T_{m-ae_jE}^t . \end{equation}
All lower coefficients belong to $uA_0$. The two families
therefore have the same bar correction by \eqref{I:correction}.

\begin{proposition}\label{I:covariant}
At every transverse displacement $\nu$, the series
$\sum_{a\ge0}B_g(\nu+ae_j)x^a$ belongs to the polynomial module of
Lemma~\ref{I:fibre}. Mutation takes it to the corresponding polynomial
of the corrected series at $m'$ in $t'$.
\end{proposition}
\begin{proof}
Expand $T_m^t=\sum_{\ell\preceq_{t,t'}m}d_{m\ell}M_\ell$,
one line at a time. By Lemmas~\ref{I:fibre} and~\ref{I:relative-mutation},
the primed expansion has the same coefficients. They belong to
$A_0[u^{-1}]$. On the $j$-axis only $M_m$ can occur, and its
coefficient is one by pointing.

For this order write $d=\xi\beta$, with diagonal entries one,
off-diagonal entries of $\xi$ in $uA_0$ and those of $\beta$
in $R^{\mathrm{bar}}$. The factorization exists and is unique by
the same coefficient splitting as in \eqref{I:correction}. Indeed,
to determine $\xi_{m\ell}$ and $\beta_{m\ell}$, subtract
$\sum_{\ell\prec_{t,t'}k\prec_{t,t'}m}\xi_{mk}\beta_{k\ell}$
from $d_{m\ell}$ and split the remainder. Each sum is finite, and
induction on the size of the interval determines both entries.

Form $Z_m=\sum_{\ell\preceq_{t,t'}m}\beta_{m\ell}M_\ell$ in
each chart. Its coefficients are bar invariant, and on each line
only finitely many $M_\ell$ occur by the bound on $a$.
The identity $T_m^t=Z_m+\sum_{\ell\prec_{t,t'}m}\xi_{m\ell}Z_\ell$
is valid coefficientwise in either chart, and every term in the sum
has smaller degree in both. By \eqref{I:relative} and uniqueness of
the bar correction, the expansions of $Z_m$ are the corrected
series at $m$ and $m'$, respectively. On each line mutation
identifies the two sums term by term in the $M_\ell$, proving
the assertion.
\end{proof}

\subsection{Finite Laurent expansions and the residue basis}

To pass from finiteness along each line to finiteness of the whole
support, we use a reddening sequence for the opposite matrix. A related
support bound for quantum theta functions is given in
\cite[Theorem~9.3]{TangTheta}.

\begin{lemma}\label{I:finite-support}
Each series $B_g$ has finite support. The maximum of each
displacement coordinate on its support is determined by tropical
mutation along any reddening sequence for $-B_\mu$.
\end{lemma}
\begin{proof}
We bound the support by tracking the maximum of a linear functional
through successive mutations. For an integral column covector $\theta$, let
$\sigma_g(\theta)=\sup\{r\theta:B_g(r)\ne0\}$,
allowing $+\infty$. In a nonzero polynomial at fixed $\nu$, write
$\ell_-,\ell_+$ for the smallest and largest $j$-displacements.
Lemma~\ref{I:fibre} shows that the smallest and largest displacements in the new seed are
$A-\ell_+$ and $A+h-\ell_-$, respectively. This remains true
when the exchange formula involves division by a polynomial. Let
$\theta'_j=-\theta_j$ and, for $i\ne j$, let
$\theta'_i=\theta_i+\theta_j[\operatorname{sgn}(\theta_j)b_{ij}]_+$.
Choosing the appropriate extreme according to the sign of $\theta'_j$
gives
\begin{equation}\label{I:support-transport} \sigma'_{g'}(\theta')=\sigma_g(\theta)+\theta'_j[\operatorname{sgn}(\theta'_j)s]_+. \end{equation}
Indeed the constant and transverse parts of the new maximum come
from $A+h=[s]_++\sum_i[-b_{ij}]_+\nu_i$ when $\theta'_j>0$,
and from $A=[-s]_++\sum_i[b_{ij}]_+\nu_i$ when $\theta'_j<0$.
For $\theta'_j=0$ we use that mutation preserves the set of transverse
displacements with nonzero polynomial. The formula therefore applies
also when the support is infinite.

The change of $\theta$ is coefficient mutation for the opposite
principal matrix $-B_\mu$. A reddening sequence for $-B_\mu$ is obtained
by reversing a reddening sequence for $B_\mu$,
negating the final framed matrix and relabeling by its terminal
permutation, since matrix mutation commutes with negation.
Along this sequence the
covector initially equal to $e_i^{\mathsf T}$ evolves into the
transpose of the $i$th row of the $c$-matrix. At the end this matrix
is a negative permutation matrix, so the covector is nonpositive.
The final series has nonnegative displacements and constant term one.
The supremum in this direction is therefore zero.

Reading \eqref{I:support-transport} backwards from this zero value
now determines $\sigma_g(e_i^{\mathsf T})$ as a finite integer.
It is the largest $i$th displacement in the initial seed. These
bounds for $1\le i\le n$ confine the support to a finite box.
\end{proof}

It remains to show that the Laurent polynomials obtained in this way
span the upper algebra. For this it suffices to know their pointed
degrees along a path from a seed to its predecessor.

\begin{lemma}\label{I:pointed-spanning}
Let $t$ be injective-reachable. Suppose that a family of integral
upper elements has exactly one member at each full degree in $t$ and is
compatibly pointed along a mutation path from $t[-1]$ to $t$.
Then the family is an $R$-basis of $\mathcal U_R^{\mathrm{loc}}$.
\end{lemma}
\begin{proof}
This is Qin's pointed-basis criterion
\cite[Theorems~4.3.1 and~6.2.1(1)]{QinTropical}. We prove the
quantum statement using the two adjacent charts of Lemma~\ref{I:fibre}.
In each seed $s$ on the path, successive subtraction gives a formal
expansion of an upper element $z$ in the pointed family. Its
coefficients belong to $R$ and its degrees lie in
$\operatorname{supp}_s(z)-\mathbb N^n E_s$.

Consider two adjacent seeds $s,s'$ on the path. Write $z$ as a
finite sum of the monomials $M_\ell$ for these two charts. A family
member pointed at $m$ and $m'$ has an expansion in the same
$M_\ell$, with diagonal coefficient one and every other degree
strictly smaller for $\preceq_{s,s'}$, by Lemma~\ref{I:fibre}.
The inverse of this triangular change of coefficients is defined
over $R$, since each entry uses only the finitely many degrees
in one interval of $\preceq_{s,s'}$. Applying it to the finite
expansion of $z$ gives the same coefficients in both charts.
Indeed, each resulting degree lies below one of the finitely many
degrees in the $M_\ell$-expansion of $z$ in both charts. Only
degrees between a fixed Laurent exponent and one of these upper
bounds can contribute at that exponent. These intervals are finite
in either chart, so the two sums are the unique formal pointed
expansions of $z$.

The coefficients of the formal expansion thus agree along the entire
path. Any member with nonzero coefficient has degrees
$g$ in $t$ and $g^-$ in $t[-1]$ satisfying
$g\preceq_t b$ and $g^-\preceq_{t[-1]} b^-$ for some
$b\in\operatorname{supp}_t(z)$ and
$b^-\in\operatorname{supp}_{t[-1]}(z)$.
By Lemma~\ref{I:predecessor}, $\psi^{-1}(g^-)$ is the codegree in $t$
and $\psi$ reverses dominance. Hence
\[ \psi^{-1}(b^-)\preceq_t\psi^{-1}(g^-)\preceq_t g\preceq_t b. \]
There are finitely many pairs $b,b^-$. Each dominance interval
is finite, since independence of the rows of $E$ identifies a nonempty
interval with an integral coordinate box. The expansion is therefore
finite. Its uniqueness gives linear independence.
\end{proof}

\begin{proof}[Proof of Theorem~\ref{I:localized}]
Lemma~\ref{I:finite-support} gives a finite Laurent polynomial for each
degree $m$ in the initial chart. It meets only finitely many lines transverse to the next mutation.
Proposition~\ref{I:covariant} identifies its polynomials on these lines
with those in the adjacent seed.
Iterating along a mutation sequence therefore gives the same integral
upper element, pointed at the corresponding tropical point in every
seed. Since the construction gives one such element at every full
degree, Lemma~\ref{I:pointed-spanning} gives an $R$-basis.

Consider another compatible extended matrix $E'$ with
the same $B_\mu,D$. For $m,m'$ with the same mutable coordinates,
the change of coefficients is
\[ \sum_r c(r)X^{m-rE}\ \longmapsto\ \sum_r c(r){X'}^{m'-rE'}. \]
Along the same mutation path,
the mutable parts of the two tropical degrees agree at each step.
The one-variable formula in Lemma~\ref{I:fibre} involves only that
mutable part, the mutable matrix and $d_j$. The translated elements
are therefore Laurent and compatibly pointed in the new frozen
pattern. They form an $R$-basis by Lemma~\ref{I:pointed-spanning}.
The same change of coefficients takes $I_N$ to the injective monomial
of full degree $(-N,0)$ in the new pattern. Indeed, in the successor
seed both have a single term of the prescribed degree. Thus the two
patterns have the same displacement series $S_g$ and the same
bar correction \eqref{I:correction}.

In any of these charts, \eqref{I:correction} changes the family
$S_m^t$ to the $L_m$ with diagonal one and off-diagonal
coefficients in $uA_0$. The inverse change has the same property
because each coefficient is a finite sum over decompositions of a displacement.
For an upper element $x$ the expansion in the $L_m$ is finite.
The coefficients of this expansion all lie in $A_0$ if and only if
all the $c_m^t(x)$ do. In that case the two expansions agree modulo
$u$. In particular only finitely many residue coordinates are
nonzero. It follows that $\mathcal L_t=\bigoplus_m A_0L_m$.
Proposition~\ref{I:covariant} identifies the lifts at corresponding
tropical degrees across each mutation. Their span and their residues
are therefore independent of $t$.
We denote the localized lattice by $\mathcal L^{\mathrm{loc}}$.
If an integral bar-invariant element belongs to this lattice, its
coordinates in the $L_m$ are bar-invariant Laurent polynomials
regular at $u=0$. They are integers. Its residue therefore determines
it uniquely, giving the asserted uniqueness of the lifts.
\end{proof}

The construction also respects a decomposition of the mutable matrix.
Suppose $U=U_1\sqcup U_2$ and $B_\mu=B_{U_1,U_1}\oplus B_{U_2,U_2}$. For each $U_i$ use the induced seed
with all the frozen vertices, and write $E^{(i)}$ for its exchange
matrix, extended by zero on the other mutable vertices. Write
$m=m_1+m_2$, where the mutable part of $m_i$ is supported on
$U_i$, assigning the frozen exponents to either factor. Then
$L_m=[L^{(1)}_{m_1}L^{(2)}_{m_2}]_t$. Indeed, mutations in one
component leave the variables and mutable exchange rows of the other
unchanged. A reddening sequence restricts to each component, and
the normalized injective monomials factor. In the
product of terms of degrees $m_1-rE^{(1)}$ and $m_2-sE^{(2)}$,
compatibility gives the same pairing $\Lambda(m_1,m_2)$ for every
$r,s$. Taking the limits therefore factors $S_m^t$ as the
normalized product of the two component series. Multiplying the expansions of the component lifts in these series
gives the unit residue at $m$. Their normalized product is also
integral and bar invariant, so uniqueness in Theorem~\ref{I:localized}
identifies it with $L_m$.

\begin{corollary}\label{I:symmetrizer}
For each connected component of $B_\mu$ with vertex set $U_\alpha$, let
$D_\alpha^\circ$ be its positive integral symmetrizer with relatively
prime entries. Write $D|_{U_\alpha}=s_\alpha D_\alpha^\circ$, and let
$B^\alpha_{g_\alpha}(r_\alpha;v)$ be the displacement coefficients
for $B_{U_\alpha,U_\alpha},D_\alpha^\circ$. In any compatible frozen pattern,
the coefficient of $X^{m-rE}$ in $L_m$ is
\[ \prod_\alpha B^\alpha_{g_\alpha}(r_\alpha;v^{s_\alpha}), \]
where $g_\alpha,r_\alpha$ are the restrictions of $m_{\rm mut},r$
to that component. For fixed $E$, the Laurent support in every seed
and the specialization $L_m(1)$ are independent of the compatible
quantization. In particular the least exponent of each frozen variable
is independent of this choice.
\end{corollary}
\begin{proof}
On a connected component the ratios of adjacent symmetrizer entries
are prescribed by $B_{U_\alpha,U_\alpha}$. In its polarization, multiplying
$D_\alpha^\circ$ by $s_\alpha$ multiplies the compatible form
specified in Section~\ref{I:construction} by the same factor.
The substitution $v\mapsto v^{s_\alpha}$ therefore respects based
products and exchange relations. It carries each normalized injective
monomial to the corresponding monomial for the scaled form. It also
commutes with the coefficientwise $u$-adic limit and with bar
correction, since it preserves the bar involution and sends $uA_0$
into $uA_0$. The component product above and coefficient transport
in Theorem~\ref{I:localized} now give the formula. A product of these
Laurent polynomials is nonzero precisely when each factor is nonzero,
and at $v=1$ each substitution has the same value.
\end{proof}

A single injective monomial $I_{\mathbf1}$ suffices to determine
$L_m$. Within the bounds of Lemma~\ref{I:finite-support},
Lemma~\ref{I:limit} determines the coefficients of $G$ by rational
arithmetic. At each displacement, the principal part and constant term
at $u=0$ determine the bar-invariant coefficient in \eqref{I:correction}.

\subsection{Multiplication by cluster monomials and comparison with Qin}

Normalized left multiplication by a cluster monomial of $t$
translates the degree of $S_m^t$. On the distinguished lifts the
same translation gives the leading term, and every other coefficient
vanishes modulo $u$. We also obtain a formula for right multiplication
by a cluster monomial of the successor seed. The left multiplication
formula identifies the lifts with Qin's common triangular basis.

\begin{proposition}\label{I:multiplication}
Let $X_t^b$ be a localized cluster monomial of $t$, and let $J$
be a localized cluster monomial of $t[1]$ of full degree $h$ in
$t$. In the basis of Theorem~\ref{I:localized} their normalized
products have the expansions
\begin{align*}
 [X_t^b L_m]_t&=L_{b+m}+\sum_{\ell\prec_t b+m} a_\ell L_\ell,\\
 [L_m J]_t&=L_{m+h}+\sum_{\ell\prec_t m+h} b_\ell L_\ell,
\end{align*}
where all coefficients in the sums belong to $u\mathbb Z[u]$.
Normalized multiplication on either side by a frozen Laurent monomial
simply translates the degree.
\end{proposition}
\begin{proof}
Write $g=m_{\rm mut}$. Normalized left multiplication by $X_t^b$
multiplies the coefficient at displacement $r$ by
$u^{\langle b_{\rm mut},r\rangle}$, so it carries $S_m^t$ to
$S_{m+b}^t$. For right multiplication, write $h_{\rm mut}=-N$
with $N\ge0$. Then $J$ is a normalized frozen translate of $I_N$
and has the same displacement coefficients $C_N$. The normalized
product formula preceding \eqref{I:product} gives the displacement
series of $[S_m^tJ]_t$ as
\[ \rho_{g+N}(G)*\rho_{g-N}(G_N)=\rho_{g-N}\bigl(\rho_{2N}(G)*G_N\bigr)=S_{g-N}. \]
The last equality uses the fixed point identity $F_N(G)=G$ in
Lemma~\ref{I:limit}. Thus $[S_m^tJ]_t=S_{m+h}^t$.

In the expansion of $L_m$ in the $S^t$, the leading coefficient is
one and every coefficient of $S_{m-aE}^t$ with $a\ne0$
belongs to $uA_0$, by \eqref{I:correction}. Normalizing the entire
product at $m$ introduces on this summand the additional factor
$u^{\langle b_{\rm mut},a\rangle}$ on the left and
$u^{\langle N,a\rangle}$ on the right. Both exponents are
nonnegative. Passing back to the $L$-basis preserves the leading
coefficient and the vanishing of every other coefficient modulo $u$.
The products themselves belong to the integral upper algebra, so their
remaining coefficients are in $R\cap uA_0=u\mathbb Z[u]$.
Finally, a frozen Laurent monomial pairs trivially under $\Lambda$
with every displacement $-aE$. Its normalized product therefore preserves
each displacement coefficient, and the coefficients of $L_m$
depend only on $m_{\rm mut}$.
\end{proof}

The products $P_{m,N}$ used to define $S_m^t$ are Qin's products
of cluster and injective monomials, normalized here to degree $m$.
For $g=m_{\rm mut}$ the least choice $N=[-g]_+$ gives his injective
pointed element \cite[Definition~5.4.1 and (19)]{QinTri}.
Qin applies bar correction to pointed elements in the completed
quantum torus \cite[Theorem~6.1.3]{QinDCB}. Here the coefficientwise
limit of the $P_{m,N}$ supplies the series to be corrected.
In fact every $P_{m,N}$ with $N\ge[-g]_+$ has the same residue
as $L_m$. Indeed, Lemma~\ref{I:predecessor} and
Theorem~\ref{I:localized} identify $I_N$ with the distinguished
lift of its degree. The congruence
$P_{m,N}\equiv L_m\bmod u\mathcal L_t$ now follows from
Proposition~\ref{I:multiplication}.

This gives a second description of the common lattice. For each
$m$ choose any $N(m)\ge[-m_{\rm mut}]_+$, and expand an element
of the localized upper algebra over $k$ in the pointed products
$P_{m,N(m)}$, recursively in the dominance order of $t$.
The element lies in $\mathcal L^{\mathrm{loc}}$ exactly when all
coefficients of this formal expansion belong to $A_0$. Their
values at $u=0$ are its distinguished residue coordinates.
Indeed, Proposition~\ref{I:multiplication} gives a unitriangular
transition from the products to the $L_m$, with lower coefficients
in $u\mathbb Z[u]$. Its formal inverse has the same property
because only finitely many chains contribute at a fixed displacement.
Thus either transition preserves the condition on coefficients
and induces the identity on residues.

For fixed $B_\mu,D$, the change of frozen coefficients in
Theorem~\ref{I:localized} identifies lifts that are \emph{similar}
in Qin's sense. His \emph{common triangular basis} is pointed at
all tropical points and is triangular in every seed $t$
\cite[Definitions~2.20--2.21 and~4.4]{QinSimilarity}. Its elements are
bar invariant, it contains the localized cluster monomials in $t$
and $t[1]$, and normalized left multiplication by each seed variable
is unitriangular with lower coefficients in $u\mathbb Z[u]$.

\begin{theorem}\label{I:qin-comparison}
Let $B_\mu$ be skew-symmetrizable and admit a reddening sequence.
For every compatible quantum seed with mutable matrix $B_\mu$, the
distinguished lifts $L_m$ of Theorem~\ref{I:localized} form the
unique common triangular basis of the fully localized quantum upper
cluster algebra over $\mathbb Z[v^{\pm1}]$.
\end{theorem}
\begin{proof}
Theorem~\ref{I:localized} gives a bar-invariant basis with one
compatibly pointed element at each tropical point. For a localized
cluster monomial, choose the basis element with the same degree in
its own seed. Lemma~\ref{I:predecessor} identifies the two elements.
Proposition~\ref{I:multiplication} gives the required unitriangularity
in every seed.

For uniqueness, use Qin's unitriangularity argument
\cite[Lemma~6.3.2(i)--(ii)]{QinTri}. In any triangular basis the
expansion of $P_{m,[-g]_+}$, with $g=m_{\rm mut}$, has diagonal
coefficient one and lower coefficients in $u\mathbb Z[u]$.
Indeed, successive normalized multiplications of the injective
monomial by mutable variables contribute factors $u^{d_i a_i}$
at displacements $a$. Normalized frozen products are bar invariant,
so their lower coefficients vanish and frozen inverses translate
degrees as well. The description of the lattice by the $P_{m,[-g]_+}$
therefore gives each triangular basis element the distinguished
residue of its degree. Theorem~\ref{I:localized} identifies its
unique integral bar-invariant lift with $L_m$.
\end{proof}

\section{The lattice with polynomial frozen variables}\label{I:polynomial}

\subsection{Frozen orders and sufficient conditions for spanning}

Let $P$ be a set of frozen vertices, and require the frozen variables
$X_f$, $f\in P$, to occur with nonnegative exponents.
The lifts $L_g$ that are regular at these divisors belong to the
resulting upper algebra. We ask when they span it. The \emph{frozen
order} $\nu_f(x)$ of a nonzero Laurent polynomial $x$ is the least
exponent of $X_f$ in its expansion in based monomials. For an element
of the localized upper algebra this order is independent of the seed,
since the mutation formula of Lemma~\ref{I:valuation-mutation} fixes the
coordinate covector at $f$. We set $\nu_f(0)=+\infty$.
The integral form $\mathcal U_{R,P}$ therefore consists of the
localized upper elements of nonnegative order at every $f\in P$.
Its \emph{common lattice} is the intersection
\[ \mathcal L_P=\mathcal L^{\mathrm{loc}}\cap\bigl(k\otimes_R\mathcal U_{R,P}\bigr). \]
In a reference seed let $\mathcal B_P$ be the set of full degrees
$g$ of the lifts regular at every $f\in P$.

Write $g=(h,a)$ for the mutable and frozen parts of a degree, and
let $e_f^*=e_f^{\mathsf T}$ be the coordinate covector at $f$.
The $f$-coordinate of the displacement $-rE$ is $rw_f^{\mathsf T}$,
where $w_f=-(E_{U,f})^{\mathsf T}$. Thus
$\nu_f(L_g)=a_f+\min_{B_h(r)\ne0}rw_f^{\mathsf T}$.
In a sum of distinct lifts, cancellation among terms of least frozen
order could make the sum regular even when some summands have poles.
If $w_f$ is nonnegative in some seed, the pointed monomials in
that seed give the order of every linear combination. For such
divisors Qin proved that the regular elements of a compatibly pointed
basis span the polynomial upper algebra
\cite[Proposition~2.15]{QinSimilarity}.

Let $\mathcal F(B_\mu^{\mathsf T})$ be the \emph{cluster $g$-vector
fan} for the mutable matrix $B_\mu^{\mathsf T}$ in the row convention.
Its cones are spanned by the mutable $g$-vectors of the clusters,
expressed in the reference seed. This is the $g$-vector fan for
$-B_\mu$ in the column convention of \cite{NZ}. Its support describes
which vectors $w_f$ become nonnegative in some seed. By
continuity, the same spanning statement holds throughout its closure.

We shall use the same valuation calculation to compare frozen orders
and to check integrality of Laurent coefficients. Give the scalar
field $k$ either the trivial valuation or the order
$\operatorname{ord}_\pi$ at an irreducible $\pi\in R$, and denote
this valuation by $\eta$. In both cases $\eta(v)=0$. For a real
column $\ell$ and a nonzero Laurent polynomial in the seed $t$, set
\[ \nu_{t,\ell}\Bigl(\sum_a c_aX_t^a\Bigr)=\min_{c_a\ne0}\{\eta(c_a)+a\ell\}. \]

\begin{lemma}\label{I:valuation-mutation}
For $t'=\mu_jt$, the valuation $\nu_{t,\ell}$ on the common
division ring equals $\nu_{t',\ell'}$, where $\ell'_a=\ell_a$
for $a\ne j$ and
\[ \ell'_j=-\ell_j+\min\{[B_{j,V}]_+\ell,[-B_{j,V}]_+\ell\}. \]
In particular, an element Laurent over $k$ in every seed has
coefficients in $R$ in every seed as soon as it does in one.
\end{lemma}
\begin{proof}
In each of two Laurent polynomials take the lexicographically least
exponent among the terms of minimum value. Their sum is attained by
just one pair of minimum-value terms. Its coefficient in the product
has the predicted value because the quantum commutation factors are
units of value zero. This proves multiplicativity and gives a valuation
on the division ring of fractions of the torus. The two monomials in
the exchange numerator have distinct exponents, since $B_{j,V}\ne0$.
The value of the new cluster variable is therefore $\ell'_j$, also
when the two values are equal. Expanding a Laurent polynomial of $t'$
gives $\nu_{t,\ell}\ge\nu_{t',\ell'}$ on that torus. The inequality
reverses on the old torus, since the reverse mutation interchanges
the numerator exponents and recovers $\ell$. Thus the valuations agree
on the intersection of the two tori. Write a fraction as $p q^{-1}$
in the old torus and multiply $p,q$ on the right by the same
sufficiently large power of $X_j$. Both polynomials now have
nonnegative $j$th exponents and hence belong to the intersection.
This proves equality on the division ring. Finally take $\ell=0$ and
$\eta=\operatorname{ord}_\pi$ for each irreducible $\pi$. The
least coefficient order is independent of the chart, and $R$ is a
unique factorization domain.
\end{proof}

\begin{proposition}\label{I:regular}
Suppose that $w_f$ lies in the closure of the support of
$\mathcal F(B_\mu^{\mathsf T})$ for every $f\in P$.
Then $\{L_g:g\in\mathcal B_P\}$ is an $R$-basis of $\mathcal U_{R,P}$.
In particular, if $\mathcal F(B_\mu^{\mathsf T})$ has dense support, the
conclusion holds for every compatible frozen pattern with mutable part
$B_\mu$ and every choice of $P$.
\end{proposition}
\begin{proof}
Fix a nonzero finite sum $z=\sum_g c_gL_g$ with distinct labels and
$c_g\in R\setminus\{0\}$. We show that
$\nu_f(z)=\min_g\nu_f(L_g)$ at each $f\in P$.
Use the trivial scalar valuation in Lemma~\ref{I:valuation-mutation}
and abbreviate $\nu_{t,\ell}$ to $\nu_\ell$.

Write $w(\ell)=-(E\ell)^{\mathsf T}$. Suppose first that $w(\ell)\ge0$ in this seed. The pointed term of
each $L_g$ has least weight. Among the summands of least weight,
choose one whose pointed degree is maximal in the dominance order.
No other summand can contain its pointed monomial. Such a summand
would have the same least weight and strictly larger pointed degree.
The coefficient of that monomial therefore survives, proving
\begin{equation}\label{I:minimum} \nu_\ell(z)=\min_g\nu_\ell(L_g). \end{equation}

Under mutation at $k$, take $\ell'$ from
Lemma~\ref{I:valuation-mutation}. For $w=w(\ell)$ and
$w'=-(E'\ell')^{\mathsf T}$, matrix mutation
gives $w'_k=-w_k$ and
$w'_j=w_j+[-b_{jk}]_+w_k+b_{jk}\min(w_k,0)$ for $j\ne k$.
This is tropical $g$-vector mutation for the row matrix $B_\mu^{\mathsf T}$
\cite[Proposition~4.2]{NZ}, using sign
coherence \cite[Corollary~5.5]{GHKK}. Along a mutation sequence to a
cluster, the rule carries its cone onto the nonnegative orthant in that
cluster's coordinates. Thus \eqref{I:minimum} holds
whenever $w(\ell)$ belongs to the support of the fan. In the same chart the value of a compatibly pointed
element is the value of its pointed monomial. Thus two such elements
of the same tropical degree have equal values. This also compares a
quantum element with a classical one, since the covector mutation
rule is the same in both tori.

The real linear map $\ell\mapsto w(\ell)$ is surjective because
$E$ has full row rank. If $J$ is a linear right inverse and $w_a\to w_f$ in the support
of the fan, then $\ell_a=e_f^*+J(w_a-w_f)$ satisfies
$w(\ell_a)=w_a$ and $\ell_a\to e_f^*$.
For fixed $x$, the function $\ell\mapsto\nu_\ell(x)$ is the minimum
of finitely many linear functions in the reference seed. Passing to
the limit gives \eqref{I:minimum} at $f$ and preserves the equality of
values for any two compatibly pointed elements of the same label.

Finally, expand an element of $\mathcal U_{R,P}$ in the localized
basis. The minimum formula at each $f\in P$ shows that every summand
is regular at all these divisors. Linear independence follows from
linear independence in the localized algebra.
\end{proof}

The frozen orders of the classical theta functions give inequalities
on their tropical labels. We will show when the same inequalities
describe the regular lifts $L_g$. To compute the orders of the theta
functions, we complete the mutable exchange rows $E$ to a
square valued matrix $B$ \cite[\S1]{Fei}. For each frozen vertex $f$
choose $d_f$ divisible by every $d_j$ with $j\in U$, and let
$b_{fj}=-(d_f/d_j)b_{jf}$. These entries are integral and satisfy
$B_{f,U}D=d_fw_f$. Setting the entries between frozen vertices
to zero makes $B\operatorname{diag}(d_a:a\in V)$ skew-symmetric.

For a frozen vertex $f$, a seed is \emph{optimizing} if the mutable
entries of its frozen row are all nonpositive or all nonnegative.
In the first case $f$ is a sink relative to the mutable vertices,
and in the second it is a source. We call the polynomial frozen vertices chosen for these
sink or source actions the \emph{colors} of the boundary actions in
Section~\ref{B:boundary-action}.

The vector $w_f$ belongs to the support of the fan precisely when
it is nonnegative in some seed. Indeed, covector mutation fixes
$e_f^*$, and $w(e_f^*)$ gives the negative transpose of the
$f$th coefficient column in every seed. Its entries have the same
signs as the entries of $B_{f,U}$. Under reddening the predecessor of a
sink-optimizing seed is source-optimizing for the same vertex
\cite[Lemma~1.5]{Fei}. In particular all the vectors $w_f$ of the
polynomial polarization lie in the support, since they are signed
coordinate vectors.

\begin{proposition}\label{I:frozen-comparison}
Let $\vartheta_g$ be the classical theta function with tropical label
$g$, and write $L_g(1)$ for the specialization of $L_g$ at $v=1$.
Under mutable reddening, every frozen divisor $f$ satisfies
\[ \nu_f(L_g)\leq\nu_f\bigl(L_g(1)\bigr)\leq\nu_f(\vartheta_g). \]
All three orders agree for every $g$ at a
divisor whose vector $w_f$ lies in the closure of the support of
$\mathcal F(B_\mu^{\mathsf T})$.
If all vectors $w_f$ with $f\in P$ lie in this closure, the labels of the
regular $R$-basis of $\mathcal U_{R,P}$ in each seed are the
integral points of the same closed rational polyhedral cone as the
regular theta labels.
\end{proposition}
\begin{proof}
With no mutable vertices the two functions are the same based monomial
and its specialization. Otherwise we expand $L_g$ in the quantum
theta basis $\{\Theta_h\}$ of the localized integral upper algebra
\cite[Theorem~1.3(1)]{TangTheta}. The order of a finite theta expansion
is the least order of its nonzero summands
\cite[Theorem~1.2]{TangTheta}. Moreover, the Newton polytope of
$\Theta_h$ equals that of its classical specialization $\vartheta_h$
\cite[\S8]{TangTheta}.
In the expansion $L_g=\sum_h c_{hg}\Theta_h$, pointedness gives
$c_{gg}=1$ and $h\preceq_t g$ whenever $c_{hg}\ne0$.
Specializing the Laurent expansion can only increase its least frozen
exponent. In the theta expansion the coefficient of $\vartheta_g$
remains one, and the order of the classical sum is again the minimum
over its nonzero summands \cite[Proposition~8.3]{TangTheta}.
This gives the two inequalities.
When $w_f$ lies in the stated closure, $L_g$ and $\vartheta_g$
have equal order by the comparison of compatibly pointed elements
in the proof of Proposition~\ref{I:regular}. The inequalities just
proved give the order of $L_g(1)$ as well.

Each frozen order of $\vartheta_g$ is the minimum
of finitely many rational linear forms in $g$
\cite[Lemma~3.1 and the proof of Corollary~3.2]{Fei}. Taking the inequalities
for just the divisors in $P$ gives the asserted cone in every seed.
Proposition~\ref{I:regular} gives the basis assertion.
\end{proof}

For an individual label $g$, convexity sometimes gives equality
beyond the closure of the fan. In the reference seed let
$P_\vartheta(g)$ and $P_L(g)$ be the convex hulls in $\mathbb R^U$
of the indices $r$ for which the exponent $g-rE$ occurs in
$\vartheta_g$ and $L_g$, respectively. For either polytope $Q$ write $m_Q(w)=\min_{r\in Q}rw^{\mathsf T}$.
The theta expansion above gives $P_\vartheta(g)\subseteq P_L(g)$
by valuative independence and independence of the rows of $E$.
Suppose $w_f=\sum_j t_j w_j$, where $t_j\geq0$, the $w_j$
lie in the fan closure, and they attain their minima at a common
point of $P_\vartheta(g)$. The proof gives equal minima for the two
polytopes at each $w_j$. Superadditivity then gives
$m_{P_L(g)}(w_f)\geq\sum_j t_jm_{P_\vartheta(g)}(w_j)
=m_{P_\vartheta(g)}(w_f)$, and containment gives the reverse inequality.
Adding the common frozen exponent $g_f$ gives
$\nu_f(L_g)=\nu_f(\vartheta_g)$.

\begin{corollary}\label{I:all-colors}
At a frozen vertex admitting an optimizing seed, the three orders in
Proposition~\ref{I:frozen-comparison} agree for every tropical label.
If every vertex of $P$ admits an optimizing seed, the regular lifts
form an $R$-basis of $\mathcal U_{R,P}$, and their specializations
at $v=1$ form a $\mathbb Z$-basis of the classical upper cluster
algebra with the same polynomial frozen variables. In each seed
their labels are the integral points of the regular theta cone for
$P$. In particular, these conclusions hold when every
polynomial frozen vertex is chosen as a reachable color.
\end{corollary}
\begin{proof}
If $f$ admits an optimizing seed, then $w_f$ lies in the support
of the fan, since the predecessor of a sink-optimizing seed is
source-optimizing. Proposition~\ref{I:frozen-comparison} therefore
gives equality of the three orders. In a seed $t$ where $f$ is a
source, the $f$th column of $E_t$ is nonpositive, so the pointed
monomial attains the least $f$-exponent in each expansion. At a sink
the column is nonnegative, so the copointed monomial attains the
minimum. Its coefficient is one by Lemma~\ref{I:predecessor},
applied to the pointed monomial in the predecessor.

If every $f\in P$ admits an optimizing seed,
Propositions~\ref{I:regular} and~\ref{I:frozen-comparison} give the
$R$-basis of regular lifts and identify its labels with the integral
points of the regular theta cone. By Proposition~\ref{B:fiber},
specializing this basis at $v=1$ gives a $\mathbb Z$-basis of the
classical upper cluster algebra with the same polynomial variables.
\end{proof}

Multiplication by frozen monomials translates the degrees of the lifts,
so we can test spanning with one frozen variable left polynomial at a
time. At a single divisor the order of a theta expansion is the least
order of its summands. This lets us express spanning in terms of the
orders of $L_g$ and $\vartheta_g$.

\begin{proposition}\label{I:spanning-reductions}
For any set $P$ of polynomial frozen vertices, the following are
equivalent.
\begin{enumerate}
\item The regular lifts $L_g$ span $\mathcal U_{R,P}$ over $R$.
\item They span $k\otimes_R\mathcal U_{R,P}$ over $k=\mathbb Q(v)$.
\item For each $f\in P$, the lifts $L_g$ regular at $f$
span $k\otimes_R\mathcal U_{R,\{f\}}$ over $k$.
\item For each $f\in P$ and every tropical label $g$,
$\nu_f(L_g)=\nu_f(\vartheta_g)$.
\end{enumerate}
\end{proposition}
\begin{proof}
Extending scalars preserves the coefficients in the basis of
Theorem~\ref{I:localized}. For $x\in\mathcal U_{R,P}$ those
coefficients already belong to $R$, and (2) forces them to vanish
at every nonregular label. This gives (1), and extension of scalars gives
the converse. The same argument applies to a single divisor.

Assume (1) and fix $f\in P$. For $x\in\mathcal U_{R,\{f\}}$,
multiply $x$ by a frozen monomial with sufficiently large exponents
at the vertices of $P\setminus\{f\}$ and zero exponent at $f$.
The product is regular along $P$. Multiplication by this monomial
takes each $L_g$ to a power of $v$ times the lift of translated degree
and preserves its order at $f$. Applying (1) to the product shows
that every summand of $x$ is regular at $f$, proving (3).
Conversely, apply (3) to the localized expansion of any
$x\in\mathcal U_{R,P}$, one divisor at a time. Each summand is
regular at every $f\in P$, giving (1).

To compare (3) and (4), fix $f$.
Under (3), the integral spanning assertion also holds for this single
divisor. For a nonzero localized element $x$ of order $d$, expand
$X_f^{-d}x$ in regular lifts. Multiplication by $X_f^d$ translates
these lifts up to powers of $v$, so every summand of $x$ has
order at least $d$. Apply this to the expansion of $\Theta_g$
in the $L_h$. Its coefficient on $L_g$ is one by pointedness, hence
$\nu_f(L_g)\geq\nu_f(\Theta_g)=\nu_f(\vartheta_g)$.
Proposition~\ref{I:frozen-comparison} gives (4).

Conversely, assume (4) and consider a nonzero finite sum of distinct
$L_g$ with nonzero scalar coefficients. Among the summands of least
order $d$, choose one with maximal pointed degree $g$ in the
dominance order. Its diagonal term $\Theta_g$ has order $d$.
A summand of higher order contains only theta terms of higher order.
Among the other summands of order $d$, only one of strictly larger
degree could contribute to $\Theta_g$. Thus its coefficient survives.
Theta valuative independence gives order $d$ for the sum. In particular
the localized expansion of any element regular at $f$ uses only lifts
regular at $f$, proving (3).
\end{proof}

In a fixed seed, one regularity test for each mutable degree suffices
for condition~(4). For $p\in\mathbb Z^U$ choose $a(p)\in\mathbb Z^F$
with $a(p)_f=-\nu_f(\vartheta_{(p,0)})$ for $f\in P$ and zero
elsewhere. Then the regular lifts span if and only if
$L_{(p,a(p))}$ is regular along $P$ for every $p$. Indeed, for both
families, translation by a frozen exponent $b$ adds $b_f$ to the
order at $f$. Thus $\vartheta_{(p,a(p))}$ has order zero at every
$f\in P$. By Proposition~\ref{I:frozen-comparison} the corresponding
lift has order at most zero there. Regularity forces equality, which
then holds at every frozen translate.

The first author's Newton polytope conjecture
\cite[Conjecture~6.11]{FeiTropical} predicts the equalities in (4)
for every $P$, as explained in Remark~\ref{B:spanning-expectation}.

\subsection{The distinguished residue basis and its balanced lifts}

A regular lift belongs to both the polynomial target and the
localized lattice. We first show that passage to the target preserves
its residue and the uniqueness of its lift. Spanning enters when we
ask whether these residues account for the whole target lattice.
For comparison with Kashiwara's definition we shall also use the
local ring $A_\infty=\mathbb Q[u]_{(u)}$ of rational functions regular
at $v=\infty$ \cite{Kashiwara,KashiwaraGlobal}.

\begin{corollary}\label{I:balance}
For any $P$, the residues of the regular $L_g$ are linearly
independent in $\mathcal L_P/u\mathcal L_P$, and $L_g$ is the
unique bar-invariant integral lift in $\mathcal L_P$ of its residue.
These residues form a basis if and only if the regular lifts span
$\mathcal U_{R,P}$. In fact, spanning the residue space
$(\mathcal L_P/u\mathcal L_P)\otimes_{\mathbb Z}\mathbb Q$ already
implies spanning over $R$.
When the regular lifts span, the triple
$(\mathcal U_{R,P},\mathcal L_P,\overline{\mathcal L_P})$ is
\emph{balanced} over $\mathbb Z$. In other words, reduction identifies
$\mathcal U_{R,P}\cap\mathcal L_P\cap\overline{\mathcal L_P}$
with $\mathcal L_P/u\mathcal L_P$.
\end{corollary}
\begin{proof}
Since $k\otimes_R\mathcal U_{R,P}$ is a $k$-subspace, we have
$\mathcal L_P\cap u\mathcal L^{\mathrm{loc}}=u\mathcal L_P$.
Thus inclusion induces an injection on residue modules. In the
localized lattice the residues of all $L_g$ form a basis by
Theorem~\ref{I:localized}. In the basis of distinguished lifts,
the three localized modules intersect in a direct sum of copies of
$R\cap A_0\cap\overline{A_0}=\mathbb Z$, since an integral Laurent
polynomial regular at zero and infinity is constant. Reduction sends
each copy isomorphically to $A_0/uA_0$. A bar-invariant integral
lift in $\mathcal L_P$ of a regular residue is also a lift in this
localized intersection and hence equals $L_g$.

If the regular lifts span $\mathcal U_{R,P}$, the expansion of a
target element uses only labels in $\mathcal B_P$. Thus
$\mathcal L_P=\bigoplus_{g\in\mathcal B_P}A_0L_g$, and the same
coordinate calculation identifies the triple intersection with the
$\mathbb Z$-span of the regular lifts. Reduction gives the asserted
isomorphism.

Conversely, suppose the regular lifts do not span. By
Proposition~\ref{I:spanning-reductions}, the target over $k$ contains
an element outside their $k$-span. Subtract its regular summands in
the localized basis. The remaining vector still belongs to the
target and has at least one nonzero coefficient. Multiply by a
common denominator and cancel the largest power of $u$ dividing
all the resulting coefficients in $\mathbb Z[u]$. This gives a
vector of $\mathcal L_P$ whose localized residue is nonzero and
supported entirely on irregular labels. Even over $\mathbb Q$ this
residue lies outside the span of the regular ones.
\end{proof}

\subsection{Recognition by highest residues}

Suppose the module already has an upper global basis. To identify
its vectors with the distinguished lifts, it suffices under the
hypotheses below to prove one lattice inclusion and to identify the
highest residues. The common bar forces the change of basis to have
constant coefficients, and the lowering operators determine the
other residues from the highest ones.
For a basis $\mathcal E$ of a $k$-vector space $M$, let
$\mathcal L(\mathcal E)$ be its $A_\infty$-span and
$B(\mathcal E)$ its set of residues modulo $u\mathcal L(\mathcal E)$.

\begin{proposition}\label{I:recognition}
Let $\mathcal D$ and $\mathcal G$ be bases of $M$ fixed elementwise
by the same semilinear involution over $v\mapsto v^{-1}$. Suppose
that every element of $\mathcal G$ has coordinates in
$\mathbb Q[v^{\pm1}]$ in the basis $\mathcal D$.
If $\mathcal L(\mathcal G)\subseteq\mathcal L(\mathcal D)$, the
two lattices coincide. If in addition they are bases of the same $R$-form,
their $A_0$-spans also coincide.

Under this inclusion, suppose a family of lowering operators on the
common residue space preserves $B(\mathcal D)\cup\{0\}$. Let
$H\subseteq B(\mathcal G)$ be a set from which every element of
$B(\mathcal G)$ can be obtained by successive applications of these
operators. Then $\mathcal G=\mathcal D$ if and only if
$H\subseteq B(\mathcal D)$.

There is also a sufficient condition for the inclusion. Suppose the
lowering operators on $M$ preserve $\mathcal L(\mathcal D)$ and
$\mathcal L(\mathcal G)=\bigoplus_\lambda T_\lambda$, where each
$T_\lambda$ is a finitely generated lattice stable under the same
operators. It suffices that each $T_\lambda$ contain a vector
$p_\lambda\in\mathcal L(\mathcal D)$ whose images under successive
lowerings, including $p_\lambda$ itself, span $T_\lambda/uT_\lambda$
on reduction.
\end{proposition}
\begin{proof}
Let $g\in\mathcal G$. By the lattice inclusion its coordinates in
$\mathcal D$ are regular at infinity. The common bar fixes each
coordinate, so each is also regular at zero. A Laurent polynomial
with both properties is constant. Thus the rational span of
$\mathcal G$ is contained in that of $\mathcal D$. The two spans
become equal on extension of scalars to $k$, since both families are
$k$-bases of $M$. They therefore already coincide over
$\mathbb Q$. This proves equality over $A_\infty$. If both bases
span the same $R$-form, the coordinates in either direction belong
to $R\cap\mathbb Q=\mathbb Z$, giving equality over $A_0$.

Suppose $H\subseteq B(\mathcal D)$. The lowering operators preserve
$B(\mathcal D)\cup\{0\}$, so every residue in $B(\mathcal G)$
belongs to $B(\mathcal D)$. Reduction modulo $u$ leaves constant
coordinates unchanged. Thus whenever $g\in\mathcal G$ and
$d\in\mathcal D$ have equal residues, $g=d$. This proves
$\mathcal G\subseteq\mathcal D$, and equality follows because
$\mathcal G$ spans $M$. The converse follows from
$H\subseteq B(\mathcal G)$.

For the last assertion, the vectors obtained by lowering
$p_\lambda$ lie in $\mathcal L(\mathcal D)$ and generate $T_\lambda$
over $A_\infty$ by Nakayama's lemma.
\end{proof}

For a highest weight crystal we may take one highest residue from each
component. The proposition requires equality of these residue vectors, including
their normalization.

\section{Quantum boundary operators and biperfectness}\label{B:boundary-action}

To recover the crystal arrows on the common lattice, we quantize
the classical boundary derivations and compute their divided powers
on the distinguished lifts. The term of greatest string height gives
the raised label. A bound on the degrees of all other coefficients
allows us to isolate this term on reduction modulo $u$.

\subsection{Sink and source operators}\label{B:section}

The first author constructed a boundary derivation for each reachable
sink \cite[\S2]{Fei}. In a seed where the chosen polynomial frozen
vertex is a sink, the derivation differentiates with respect to its
variable and multiplies by a monomial in the adjacent mutable variables.
To quantize it on the upper algebra of
Section~\ref{I:polynomial}, we express both operations as a single
commutator that is independent of the choice of seed.

The coefficients of this commutator depend on the symmetrizer at the
frozen vertex. We therefore extend the mutable rows $B_{U,V}$ to a
full valued exchange matrix $B=(b_{ab})_{a,b\in V}$ and fix positive
integral weights $D_V=\operatorname{diag}(d_a)_{a\in V}$ extending
$D$ such that $BD_V$ is skew-symmetric. Here $U$ is the mutable
vertex set and $V$ the full vertex set. Compatibility with the quantum
torus reads $B_{U,V}\Lambda=-(D\mid 0)$.
We write $q_i=v^{d_i}$ and use the symmetric quantum integers
$[n]_q=(q^n-q^{-n})/(q-q^{-1})$ and their factorials and binomial
coefficients.

Suppose $i\in P$ is a sink relative to the mutable vertices in some
seed. We call any such seed \emph{sink-optimizing} for $i$. Its
exchange entries satisfy $b_{ik}\leq 0$ for every $k\in U$.
The \emph{boundary monomial} $M_i=X^{m_i}$, with
$m_i=-\sum_{k\in U}b_{ik}e_k$, specializes to the monomial
that multiplies $\partial/\partial x_i$. Differentiation lowers the $i$th exponent by one,
so the full shift on monomials is $s_i=m_i-e_i$. The based monomial
$Z_i=X^{s_i}$ will enter the commutator. To express the degree of $Z_i$ in other seeds, let $\zeta_i$
denote the frozen part of the $i$th row of $B$ in this seed,
with zero entries at the mutable vertices. By \cite[Lemma~2.1]{Fei},
$m_i$ represents the same tropical point in every
optimizing seed, and $\zeta_i$ is independent of this choice.
Lemma~\ref{I:predecessor} shows that a localized cluster monomial is
determined by its tropical point. Thus both $M_i$ and its
normalized frozen translate $Z_i$ are well defined. The pointed
degree of $Z_i$ in a seed $t$ is
$\tau_i(t)=-B_{i,V}(t)-e_i+\zeta_i$.

\begin{lemma}\label{B:grading}
The functional $h_i(a)=d_i a_i-\Lambda(s_i,a)$ in a sink-optimizing
seed defines an integral grading. In a seed $t$, the degree of
$X^a$ is $d_i a_i-\Lambda(\tau_i(t),a)$.
\end{lemma}
\begin{proof}
For $k\in U$, compatibility gives
$\Lambda(s_i,B_{k,V})=d_k(s_i)_k=-d_k b_{ik}=d_i b_{ki}$.
Thus $h_i(B_{k,V})=0$, so both monomials in the exchange relation
at $k$ have the same degree.
To follow the grading under mutation, write $H=D_V+B\Lambda$.
Its mutable rows vanish. The functional in the statement is the sum
of the functional $a\mapsto H_{i,V}a^{\mathsf T}$ and
$\Lambda(e_i-\zeta_i,-)$. The
second summand measures commutation with a fixed frozen monomial. Under mutation at
$k$ write $B'=QBA$, where $A$ is the identity with its $k$th row
replaced by $-e_k+[B_{k,V}]_+$ and $Q$ is the identity
with its $k$th column replaced by $-e_k^{\mathsf T}+[-B_{V,k}]_+$. Here
$[\,\cdot\,]_+$ takes the positive part of each entry. Then $A^2=1$, $\Lambda'=A\Lambda A^{\mathsf T}$, and
$QD_VA^{\mathsf T}=D_V$ by skew-symmetrizability. We obtain
$H'=QHA^{\mathsf T}=HA^{\mathsf T}$, the transformation law for
the degrees of the seed variables.
\end{proof}

The grading compensates for the powers of $v$ introduced by
commuting $Z_i$ past a monomial. For a homogeneous element $x$
of degree $h_i(x)$, define the \emph{quantum boundary operator} by
\begin{equation}\label{B:comm}
 E_i(x)=\frac{v^{h_i(x)}Z_i x-v^{-h_i(x)}xZ_i}{q_i-q_i^{-1}}.
\end{equation}
Extend $E_i$ by linearity. If $\sigma_i(x)=v^{h_i(x)}x$, expanding
the commutator on a product gives the skew Leibniz rule
$E_i(xy)=E_i(x)\sigma_i(y)+\sigma_i^{-1}(x)E_i(y)$.

\begin{proposition}\label{B:operators}
In a sink-optimizing seed, $E_i X^a=[a_i]_{q_i}X^{a+s_i}$. Its
divided powers preserve $\mathcal U_{R,P}$ and commute with the
bar. On a monomial with $a_i\geq 0$ they are given by
\begin{equation}\label{B:monomials}
 E_i^{(j)}X^a=\begin{bmatrix}a_i\\j\end{bmatrix}_{q_i}X^{a+js_i},
 \qquad E_i^{(j)}=E_i^j/[j]_{q_i}!.
\end{equation}
In particular $E_i$ is locally nilpotent, and at $v=1$ it acts
as $M_i\,\partial/\partial x_i$.
\end{proposition}
\begin{proof}
On $X^a$, the two terms in the numerator of \eqref{B:comm} have the
same exponent $a+s_i$ and respective coefficients $v^{d_i a_i}$
and $v^{-d_i a_i}$. Each application of $E_i$ lowers the $i$th
exponent by one. This proves \eqref{B:monomials}, with the
binomial understood to vanish for $j>a_i$, and gives bar invariance
and the classical formula. Over $k$ the commutator belongs to the
localized upper algebra because $Z_i$ does. In the optimizing seed
\eqref{B:monomials} preserves nonnegative exponents at every frozen vertex in $P$. By seed independence of the frozen orders
(Section~\ref{I:polynomial}), this proves regularity along $P$.

Formula \eqref{B:monomials} gives integral coefficients in the optimizing
seed. Lemma~\ref{I:valuation-mutation} gives integrality in every seed.
\end{proof}

We obtain the operator $E_j^\star$ for a reachable source $j$ by
applying this construction to the opposite seed $(-B,-\Lambda)$.
The $R$-linear map $\iota:X^a\mapsto X_{\mathrm{op}}^a$ reverses
products because the commutation form changes sign. It intertwines
mutation because negating the mutable exchange row interchanges
the two based monomials in the exchange relation. Thus $\iota$
identifies the localized upper algebras and preserves every frozen
order. Transporting the opposite sink operator by $\iota$ gives
$E_j^\star$. In a source-optimizing seed the
shift in \eqref{B:monomials} is $\sum_{k\in U}b_{jk}e_k-e_j$.

The requirement $i\in P$ has a direct explanation. On $X_i^{-1}$,
the successive applications of $E_i$ introduce the factors
$[-1]_{q_i},[-2]_{q_i},\ldots$. Each is nonzero over $k$, so
inverting $X_i$ would destroy local nilpotence.

To determine the Cartan matrix for a family of sinks, we measure how
each derivation acts on the multiplier of the others. For distinct
reachable sinks $i,j$, let $a_{ij}$ be the greatest $i$th
exponent of $M_j$ in a seed optimizing $i$. The pointed term
attains this maximum because $b_{ki}\geq0$ for every $k\in U$
there. In particular
its value agrees before and after specialization. The boundary degree
formula gives $d_i a_{ij}=d_i(\zeta_j)_i+d_j(\zeta_i)_j=d_j a_{ji}$
\cite[\S2.3]{Fei}. The \emph{boundary Cartan matrix} has
entries $c_{ii}=2$ and $c_{ij}=-a_{ij}$ for $i\ne j$. For the same
vertices chosen as sources the matrix agrees whenever both families
are reachable \cite[Corollary~6.3]{Fei}. Indeed the source degrees
can be computed at $v=1$ from the pointed expansion in the opposite seed.
We write $U_v^+(C)$ for the positive half of the quantized enveloping
algebra of $C=(c_{ij})$ and use its integral form generated over $R$
by the divided powers of the simple generators \cite{Lusztig}.

\begin{theorem}\label{B:serre}
The sink operators give a representation of the divided-power
integral form of $U_v^+(C)$ on $\mathcal U_{R,P}$. The source operators likewise represent the divided-power integral form
for their boundary Cartan matrix.
\end{theorem}
\begin{proof}
It suffices to verify the quantum Serre relation over $k$.
Fix $i\ne j$, abbreviate $A=a_{ij}$, $q=q_i$, and expand
$Z_j=\sum_n z_nX^n$ in a seed optimizing $i$. For each term
$b=n_i$ satisfies $0\le b\le A$, since $M_j$ is regular at $i$.
Write $\gamma_i$ for the exchange degree of $Z_i$. Evaluating the
two functionals of Lemma~\ref{B:grading} on $\tau_j,\tau_i$ in the
same seed cancels their alternating pairings and gives
$h_i(\gamma_j)+h_j(\gamma_i)=d_i(\tau_j)_i+d_j(\tau_i)_j=d_i A$.
In the last equality the entries of $B$ cancel by
$d_i b_{ji}+d_j b_{ij}=0$, leaving the two frozen entries in
the boundary degree formula.
As $h_i(n)=h_i(\gamma_j)$, it follows that
$h_j(s_i)+\Lambda(n,s_i)=d_i(A-b)$.

To check the Serre relation, group monomials along lines parallel to
$s_i$. On the line $a=a^0-zs_i$ with $a^0_i=0$, identify
$X^a$ with $x^z$.
Then $E_i$ becomes the symmetric $q$-derivative
$\partial x^z=[z]_q x^{z-1}$. The summand $z_nX^n$ in
\eqref{B:comm} for $j$ sends this line to the one with origin
$a^0+n+b s_i$, taking level $z$ to $z+b$. With $T x^z=q^z x^z$,
the two maps it contributes are scalar multiples of $J=x^b T^p$
for $p=\pm(A-b)$, since the two powers of $v$ have exponents given by the
functionals $\pm(h_j+\Lambda(n,-))$ at $a^0-zs_i$.
Choose $\epsilon\in\{1,-1\}$ so that $p=\epsilon(A-b)$.
Evaluation on $x^z$ gives
\[ \partial J-q^{p+\epsilon b}J\partial=[b]_q x^{b-1}T^{p-\epsilon}. \]
Each application lowers $b$ by one and replaces $p$ by
$p-\epsilon$. Thus the differences with exponents
$\epsilon A,\epsilon(A-2),\ldots,\epsilon(A-2b)$ kill $J$, since
the last difference has coefficient $[0]_q$.
Left and right composition by $\partial$ commute. Thus the product of
all $A+1$ differences with exponents $A,A-2,\ldots,-A$ vanishes on
either map. Expanding the product by the quantum binomial formula
and dividing by $[A+1]_q!$ gives
$\sum_{r=0}^{A+1}(-1)^r E_i^{(A+1-r)}E_jE_i^{(r)}=0$.

Applying the same calculation in the opposite algebra proves the
source relations. The divided powers of both families preserve the
integral form by Proposition~\ref{B:operators}.
\end{proof}

To estimate the commutator on $L_g$, we need the coefficients of
multiplication by a cluster monomial at both the greatest and the least
exponent in the dominance order. The multiplication rule at the least
exponent follows from the rule at the pointed degree in the predecessor.
For the common triangular basis, these two rules and the passage to the
opposite seed are due to Qin
\cite[Proposition~4.9 and Theorem~4.10]{QinLeclerc}.

\begin{proposition}\label{B:opposite-family}
In every seed $t$ the distinguished lifts are copointed, with one
lift at each codegree. Let $L$ be the lift of codegree $\eta$
and let $X^a$ be a localized cluster monomial of $t$. In the
expansion of $v^{-\Lambda_t(\eta,a)}LX^a$, the element of codegree
$\eta+a$ has coefficient one. Every other term has strictly greater
codegree in the dominance order of $t$ and coefficient in
$u\mathbb Z[u]$. The map $\iota$ carries the distinguished
lifts onto those of the opposite seed and preserves regularity along
every polynomial frozen divisor.
\end{proposition}
\begin{proof}
Let $\psi$ be the order-reversing map from $t$ to its predecessor
in Lemma~\ref{I:predecessor}. If $L$ has pointed degree $b$ in
the predecessor, it has codegree $\psi^{-1}(b)$ in $t$, with
coefficient one. Each degree in the predecessor lattice occurs exactly once. In $LX^a$ the coefficient at the codegree
$\eta+a$ is $v^{\Lambda_t(\eta,a)}$. The lemma identifies this
coefficient with the pointed coefficient in the predecessor. There
$X^a$ is a localized monomial of the successor, so the right
multiplication formula of Proposition~\ref{I:multiplication} gives
the asserted expansion upon reversing the order by $\psi^{-1}$.

Under $\iota$ a codegree in $t$ becomes a degree in the opposite
seed. The preceding expansion gives its left multiplication rule,
since $\Lambda_{\mathrm{op}}(a,\eta)=\Lambda_t(\eta,a)$. The image
basis is bar invariant and contains all localized cluster monomials.
Thus it is triangular in the fixed opposite seed. The opposite of
$t[-1]$ is an injective successor, since the degrees $\psi^{-1}(e_j)$
of its mutable variables have mutable parts forming a negative
permutation matrix. The fixed-seed uniqueness
in the proof of Theorem~\ref{I:qin-comparison} identifies the image
basis with the distinguished basis there. Preservation of frozen
orders gives the last assertion.
\end{proof}

\subsection{Quantum biperfectness}\label{B:section-perfect}

To study a single operator at $i$, we may invert every other frozen
variable. The lifts regular at $i$ then form a basis by
Corollary~\ref{I:all-colors}. Let $a$ specify the vertex $i$ and a
choice of sink or source orientation. We write $E_a=E_i$ for a sink
and $E_a=E_i^\star$ for a source. In a seed optimizing $a$, take the
pointed degree of $L_g$ for a sink and of $\iota(L_g)$ in the opposite
seed for a source. Its $i$th coordinate is the greatest exponent of
$X_i$ in $L_g$. We call it the \emph{string height} $\rho_a(g)$.
When this height is positive, we add the shift of $E_a$ to the pointed
degree and transport the resulting label back to the reference seed,
using $\iota^{-1}$ in the source case. We denote the label by $r_a g$.

We first estimate the coefficients of the left product in the
commutator. Use full degrees in a seed optimizing a sink $i$ and
write $s=s_i$ for its boundary shift.

\begin{lemma}\label{B:product-bounds}
Let $r$ and $\ell$ be the greatest and least $i$th exponents
of $L_g$, and write $v^{h_i(g)}Z_iL_g=\sum_b A_b(v)L_b$ in
the localized upper algebra. The coefficient at $g+s$ is
$v^{d_i r}$. Every other label with $A_b\ne0$ satisfies
\[ b_i<r-1,\qquad d_i\ell\le t<d_i r\quad\text{for every exponent $t$ of $v$ in $A_b$}. \]
\end{lemma}
\begin{proof}
Let $\eta$ be the codegree of $L_g$ from
Proposition~\ref{B:opposite-family}. Every exponent of $L_g$ lies
between $\eta$ and $g$ in dominance order. Since
$b_{ki}\ge0$ for $k\in U$, we have $r=g_i$ and
$\ell=\eta_i$. The coefficient at $X^{g+s}$ in the product is
$v^{h_i(g)+\Lambda(s,g)}=v^{d_i r}$. As $s=m_i-e_i$ is
nonnegative on the mutable indices, the left multiplication formula
of Proposition~\ref{I:multiplication} gives coefficient one at
$L_{g+s}$ after division by this power of $v$. All other terms
have strictly lower degree and coefficient in $u\mathbb Z[u]$.
In particular $b_i\le r-1$ throughout the expansion and each
exponent in a nonleading $A_b$ is strictly less than $d_i r$.

We show that at height $r-1$ only the leading label occurs. For
each $b$, let $Q_b$ be the sum of the terms of $L_b$ with
$i$th exponent $b_i$. Each $Q_b$ is bar invariant and pointed at $b$. Hence they are linearly independent,
as one sees by choosing a dominance-maximal degree in any finite
linear relation. Multiplication by $Z_i=X^s$ translates exponents
by $s$. At an exponent $g-nB_{U,V}$ of height $r=g_i$,
we have $(nB_{U,V})_i=0$. By Lemma~\ref{B:grading},
$\Lambda(s,nB_{U,V})=d_i(nB_{U,V})_i=0$, so the factor from
multiplication by $Z_i$ is the same as at the pointed term.
Thus the terms of height $r-1$ in the normalized product are the
translate of $Q_g$, which is bar invariant. Their expansion in the
$Q_b$ has bar-invariant coefficients by independence. Every
coefficient off $g+s$ also belongs to $u\mathbb Z[u]$, so it vanishes.

Taking the bar gives the right product $v^{-h_i(g)}L_g Z_i$,
with coefficients $A_b(v^{-1})$. Its codegree is $\eta+s$.
Since $h_i(\eta)=h_i(g)$, its coefficient at $X^{\eta+s}$ is
$v^{-h_i(g)+\Lambda(\eta,s)}=v^{-d_i\ell}$.
Proposition~\ref{B:opposite-family} now bounds every exponent in
$A_b(v^{-1})$ above by $-d_i\ell$. This is the lower bound for
$A_b(v)$.
\end{proof}

\begin{theorem}\label{B:biperfect}
Suppose $i$ admits a seed optimizing the chosen orientation $a$ and
$L_g$ is regular at $i$. For $r=\rho_a(g)$, its expansion after
applying a divided power is
\begin{equation}\label{B:perfectformula}
 E_a^{(j)}L_g=
 \begin{bmatrix}r\\j\end{bmatrix}_{q_i} L_{r_a^j g}
 +\sum_{\rho_a(h)<r-j}c_{h,g}^{a,j}(v)L_h
 \qquad(0\le j\le r).
\end{equation}
Every label on the right is regular at $i$. Each nonzero coefficient
in the sum belongs to $\mathbb Z[v+v^{-1}]$ and its highest Laurent
degree satisfies
$\deg_v c_{h,g}^{a,j}<d_i j(r-j)$. In particular
$E_a^{(r)}L_g=L_{r_a^r g}$, and $E_a^{(j)}L_g=0$ for $j>r$.
If $L_g$ is regular at every vertex of $P\ni i$, then
$L_{r_a^r g}$ is regular at every vertex of $P$ as well. If the
regular lifts span $\mathcal U_{R,P}$, every term in
\eqref{B:perfectformula} is regular on this target.
\end{theorem}

We call this property for the two families \emph{quantum
biperfectness}. The unique term of greatest height gives the crystal
arrow at $v=1$. The strict degree bound makes all the other terms
vanish at $u=0$ after normalization.

\begin{proof}
Consider a sink $i$ and write $q=q_i$. With $A_b$ as in
Lemma~\ref{B:product-bounds}, the coefficient of $L_b$ in $E_i L_g$
is $(A_b(v)-A_b(v^{-1}))/(q-q^{-1})$ by \eqref{B:comm}.
Proposition~\ref{B:operators} and the regular basis for $P=\{i\}$ show
that these coefficients are Laurent polynomials over $\mathbb Z$
and vanish at labels that are irregular at $i$. Each is fixed by the bar.
For $r>0$ the coefficient of $L_{g+s_i}$ is $[r]_q$, so
$g+s_i$ is regular. Every other coefficient has degree strictly
less than $d_i(r-1)$, since regularity at $i$ gives
$\ell\ge0$. These terms have height at most $r-2$.

In \eqref{B:monomials}, distinct monomials that survive a fixed
divided power still have distinct exponents. Thus $r$ is the
largest $j$ for which $E_i^{(j)}L_g\ne0$. This characterizes the height
independently of the optimizing seed. The unique summand of height $r-1$ determines the raised label
in the same way. The raising map is injective
because it is a translation in any such seed.
Proposition~\ref{B:opposite-family} identifies the regular basis for
a source with that for the sink of the opposite seed, so the same
argument applies.

Following the leading term through $j$ applications of $E_i$
gives the coefficient
$[r]_{q_i}[r-1]_{q_i}\cdots[r-j+1]_{q_i}$. Every other contribution passes through a term of strictly smaller
height at some step. At that step it also has strictly smaller
degree, so its total degree is less than $d_i\sum_{\ell'=0}^{j-1}(r-\ell'-1)$.
Adding these contributions preserves the bound and division by
$[j]_{q_i}!$ subtracts $d_i j(j-1)/2$ from the degree. This proves
\eqref{B:perfectformula} and the asserted strict estimate. When $j=r$
no lower terms have nonnegative height.
For a larger $P$, Proposition~\ref{B:operators} gives
$E_a^{(r)}L_g\in\mathcal U_{R,P}$ whenever $L_g\in\mathcal U_{R,P}$.
The equality $E_a^{(r)}L_g=L_{r_a^r g}$ proves regularity of this
lift on the whole target. If the regular lifts span the target, uniqueness of expansion in
the localized basis gives regularity of every summand for each $j$.
\end{proof}

The residue calculation works for any $P$. For $g\in\mathcal B_P$
and a boundary $a$ at $i$, let $r=\rho_a(g)$ and $0\le j\le r$.
The leading quantum binomial in \eqref{B:perfectformula} is monic
of degree $D=d_i j(r-j)$ in $v$. The estimate places
$u^D E_a^{(j)}L_g$ in the localized lattice, and it belongs to the
target by Proposition~\ref{B:operators}. Under the injection
$\mathcal L_P/u\mathcal L_P\hookrightarrow
\mathcal L^{\mathrm{loc}}/u\mathcal L^{\mathrm{loc}}$ of
Corollary~\ref{I:balance}, its residue therefore maps to that
of $L_{r_a^j g}$. Whenever this lift is regular on the target, we
obtain the congruence in the target lattice
\[ u^D E_a^{(j)} L_g \equiv L_{r_a^j g} \pmod{u\mathcal L_P}. \]
If the regular lifts span, this holds for every $j$ by
Theorem~\ref{B:biperfect}.

In the polynomial polarization every frozen vertex admits an
optimizing seed. In its
reference seed the two frozen vertices of weight $d_k$ with columns
$e_k^{\mathsf T}$ and $-e_k^{\mathsf T}$ have mutable rows
$-e_k$ and $e_k$, respectively. The first gives a sink operator
and the second a source operator. In both cases the boundary monomial
is $X_k$.

We next specialize the algebra and its boundary operators at $v=1$.
Let $\mathcal U^{\mathrm{cl}}_{\mathbb Z,P}$ denote the classical
upper cluster algebra over $\mathbb Z$ with the same mutable
exchange rows and the same set $P$ of polynomial frozen variables.
When all frozen variables are polynomial and each admits an optimizing
seed, Oya, Qin and Yakimov prove that a compatibly pointed basis and
its regular part specialize to classical bases
\cite[Corollaries~3.21--3.22]{OQY}.

\begin{proposition}\label{B:fiber}
For every $P$, the regular quantum theta functions form an $R$-basis
of $\mathcal U_{R,P}$. Specialization in any seed takes this basis
to the regular classical theta basis with the same labels and induces
an isomorphism
\[ \mathcal U_{R,P}/(v-1)\mathcal U_{R,P}\ \xrightarrow{\ \sim\ }\ \mathcal U^{\mathrm{cl}}_{\mathbb Z,P}. \]
\end{proposition}
\begin{proof}
In the localization, the quantum theta functions form an $R$-basis
and specialize to the corresponding classical theta basis over
$\mathbb Z$ \cite[Proposition~3.2, Theorem~9.3 and Corollary~11.14]{TangTheta}.
At each frozen divisor the order of a nonzero theta expansion is the
least order of its summands, on both the quantum and the classical
side. Moreover, a quantum theta function and its specialization have
the same Newton polytope in each seed
\cite[Theorem~1.2 and Propositions~8.2--8.3]{TangTheta}.
Thus imposing regularity at the divisors in $P$ selects the same
subset of labels in both bases. Specialization sends one resulting
basis bijectively onto the other, proving the isomorphism.
\end{proof}

\begin{corollary}\label{B:specialization}
Suppose the regular lifts span $\mathcal U_{R,P}$.
In the classical polynomial upper algebra over a field of
characteristic zero, their specializations form a perfect basis
for each family of boundary derivations. If each family is equipped
with an exchange grading giving the weights of its crystal, the basis is
\emph{BK-biperfect} in the sense of \cite[Definition~8.9]{FeiCrystal}.
Its elements are also compatibly pointed, so it is biperfect in the
sense of \cite[Definition~9.7]{FeiCrystal}.
\end{corollary}
\begin{proof}
Proposition~\ref{B:fiber} identifies the specialized basis with a
basis of the classical target. The binomial in
\eqref{B:perfectformula} specializes to $\binom rj$.
For $j=1$ the leading coefficient is $r$, the other heights
are at most $r-2$, and the leading map is injective. The unit is
$L_0$, each $L_g$ is homogeneous for any exchange
grading, and its pointed coefficient remains one in every seed.
\end{proof}

A projective presentation weight is \emph{$\mu$-supported} if the
cokernel of a general presentation of that weight is supported on the
mutable vertices. Conjecture~8.14 of \cite{FeiCrystal} asks, in
part, whether the specialized triangular elements are BK-biperfect
whenever they form a basis of the unlocalized upper algebra indexed
by these weights. We prove this implication under mutable reddening and
compatible quantization in the following quiver setting. Here a frozen
vertex is \emph{dually reachable} if a sequence of mutations takes its
dual boundary representation to the simple at that vertex
\cite[Definition~4.7]{FeiCrystal}.

\begin{corollary}\label{B:presentation-biperfect}
Let $(\Delta,\mathcal S)$ be a skew-symmetric ice quiver with
potential over $\mathbb C$ as in \cite[Construction~2.5]{FeiCrystal}.
Assume that the full Jacobian algebra is finite dimensional, the
mutable part admits a reddening sequence, and the extended exchange
matrix admits a compatible quantization. For one such quantization,
suppose that the classical specializations $L_{-\delta}(1)$,
indexed by the $\mu$-supported projective presentation weights
$\delta$, form a basis of the unlocalized
upper cluster algebra over $\mathbb C$. Then
the regular theta basis has the same labels. Both bases are BK-biperfect for
the reachable ordinary and dually reachable starred colors with
their compatible gradings. If every frozen vertex admits an optimizing
seed, the basis hypothesis holds for every compatible quantization.
\end{corollary}
\begin{proof}
By Corollary~\ref{I:symmetrizer}, we may compute the specialized
elements in a quantization with a single multiplier $d$ on the
entire ice quiver. To construct one, use the full row rank of $E$
to prescribe $\omega(rE,a)=-r a_{\rm mut}^{\mathsf T}$ on its
rational row space against $\mathbb Q^V$. This prescription is
alternating on pairs of rows since $B_\mu$ is skew-symmetric. Extend
it to an alternating form on $\mathbb Q^V$ and clear denominators.
The resulting integral form satisfies $E\Lambda=-(dI\mid0)$.
We take $D_V=dI$ for the original full skew-symmetric matrix
$B(\Delta)$. Thus the exponent $m_i$ uses the frozen row of the
given quiver, and $[a_i]_{v^d}$ specializes to $a_i$.
At $v=1$ the sink operator therefore agrees with the root derivation
in \cite[(8.1)]{FeiCrystal}. The same calculation in the opposite
quiver identifies the starred derivation.

Fix one boundary at $i$ and first take $P=\{i\}$. By
Corollary~\ref{I:all-colors}, in an optimizing seed the least
$i$th exponent of $L_g$ occurs at its pointed or copointed monomial
with coefficient one. Thus every $L_{-\delta}(1)$ in the assumed
classical basis has a quantum lift regular at $i$. Specializing
\eqref{B:perfectformula} for that lift gives its boundary derivative.
The derivative belongs to the unlocalized classical algebra
by \cite[Theorem~2.4]{Fei}, so it has an expansion in the assumed basis.
The specialized localized elements are linearly independent by
pointedness, so the two expansions agree. In particular, at positive
height the leading term survives with coefficient $r$, and every
other term has height at most $r-2$. It remains to identify the labels.
At an ordinary color, our height and translation in an optimizing seed
are the theta height and translation for the same $g=-\delta$.
At a starred color we use the opposite seed. The pointed degree
there is the original codegree by Proposition~\ref{B:opposite-family},
and Lemma~\ref{I:predecessor} recovers this codegree from the tropical
point.
The triangular and theta labels therefore agree in both orientations.
These arrows and heights agree with those on $\mu$-supported
presentation weights \cite[Lemma~4.7 and Corollary~4.9]{Fei}.
In the opposite quiver the label is $-\check\delta$, where
$\check\delta$ is the injective presentation weight of the same
general decorated representation.

It remains to check that no regular theta label is missing. Write
$\mathcal T$ for the negatives of the $\mu$-supported weights.
Each $h\in\mathcal T$ is a regular theta label
\cite[Lemma~4.7]{Fei}. Conversely, expand a regular $\vartheta_g$
in the assumed basis $\{L_h(1):h\in\mathcal T\}$. In a fixed seed,
the monomial of each maximal degree $h$ in this finite expansion
survives. Pointedness of $\vartheta_g$ therefore gives $h\preceq g$
for every summand. The coefficient at $g$ can now come only from
$L_g(1)$, so $g\in\mathcal T$. Compatibility gives full row rank,
and hence spanning by the regular theta functions
\cite[Proposition~1.3]{Fei}. The regular theta basis is
BK-biperfect for the same presentation crystals
\cite[Corollary~4.9]{Fei}. The triangular basis is integral as well.
Indeed, all $L_g(1)$ form a $\mathbb Z$-basis of the localized
classical algebra by Theorem~\ref{I:localized} and
Proposition~\ref{B:fiber} with $P=\varnothing$. An integral
element in the complex span of the indicated subset therefore has
integral coefficients in that subset.

Finally, suppose every frozen vertex admits an optimizing seed.
The frozen-order formula identifies $\mathcal T$ with the regular
theta labels \cite[Proposition~4.8]{Fei}. By Corollary~\ref{I:all-colors}
the $L_g$ for $g\in\mathcal T$ form a basis of the unlocalized
integral form. Their specializations therefore give the required
classical basis by Proposition~\ref{B:fiber}.
\end{proof}

\subsection{One-sided strings on the common lattice}\label{B:section-strings}

For a fixed boundary $i$, we seek a basis of the common lattice
with the same residues as the distinguished lifts. On each string
the boundary operator should send the vector of height $j>0$ to
$[j]_{q_i}$ times the preceding vector. Assume the regular lifts span
$\mathcal U_{R,P}$. By Proposition~\ref{I:spanning-reductions}, it
suffices to check spanning over $k$ in $M=k\otimes_R\mathcal U_{R,P}$.
The string heights determine the kernels of all powers of the operator.

\begin{lemma}\label{B:filtrations}
For any family of boundaries $a$ and integers $j_a\ge0$,
\[ \bigcap_a\ker_M E_a^{j_a+1}=\bigoplus_{\{g\in\mathcal B_P:\,\rho_a(g)\le j_a\ {\rm for\ all}\ a\}}kL_g. \]
Intersecting with $\mathcal U_{R,P}$ or $\mathcal L_P$ gives the
corresponding direct sum over $R$ or $A_0$.
\end{lemma}
\begin{proof}
This is the quantum counterpart of the classical kernel formula
\cite[Corollary~8.11]{FeiCrystal}. For one boundary, take the greatest height $r$ in the basis
expansion of a nonzero vector. By Theorem~\ref{B:biperfect},
$E_a^{(r)}$ sends the terms of height $r$ to distinct basis
elements and kills every other term. Thus the vector has height
exactly $r$. Each kernel is consequently spanned by a subset of
the same basis, and taking intersections gives the assertions.
\end{proof}

To construct a string for one boundary operator $E$, we correct
each lift by vectors of smaller height. Write
$M_{\le j}=\ker E^{j+1}$ and $M_{\le-1}=0$ for the resulting
filtration.

The raising map for $E$ is injective on the labels of positive
height and lowers height by one. Iterating it from any label ends
at a unique label of height zero. The strings in the optimizing
seed therefore have labels $p-js$,
where $p$ has height zero and $j$ ranges over an initial interval
of $\mathbb N$. For either orientation we work in the corresponding
sink chart.

\begin{theorem}\label{B:strings}
For each boundary $i$ there is an $A_0$-basis of $\mathcal L_P$
consisting of strings indexed by the $i$-strings of regular
labels, with the same length at each highest label. The vector $w_j$ at height
$j$ in the string through $L_p$ satisfies
\[ Ew_0=0,\qquad Ew_j=[j]_{q_i}w_{j-1},\qquad w_j-L_{p-js}\in u(\mathcal L_P\cap M_{\le j-1}). \]
Define two linear operators on the string basis by
$w_j\mapsto w_{j-1}$ and $w_j\mapsto w_{j+1}$, taking the value
zero whenever the indicated vector is absent.
The construction below gives the same vectors in every optimizing
seed. Both shifts preserve $\mathcal L_P$. On its residue basis
they give the raising arrow and its partial inverse for that boundary.
\end{theorem}
\begin{proof}
Write $q=q_i$ and $b_j=L_{p-js}$. If the string has last height
$N$, take $w_j=E^{(N-j)}b_N/\left[\begin{smallmatrix}N\\j\end{smallmatrix}\right]_q$.
The maximal raise is $w_0=b_0$, and the divided-power identity
gives $Ew_j=[j]_q w_{j-1}$ for $j>0$.
The binomial in the denominator has highest $v$-degree
$d_i j(N-j)$ with coefficient one. After division by this binomial, each lower coefficient belongs to
$uA_0$ by Theorem~\ref{B:biperfect}, which proves the congruence.

For an infinite string we construct the successive vectors by
multiplication. The required multiplier is $S=L_{-s}$. In the optimizing
seed $-s$ has mutable part $-m_i|_U$. Thus $S$ is the
injective monomial with multiplicities $m_i|_U$, translated in
the frozen directions. Let $\phi$ be the tropical transformation
to the successor. Since $b_j$ is regular, the exponent
$\phi(p-js)_f$ of its pointed term is nonnegative for every
$f\in P$. Dividing by $j$ and taking the limit gives
$\phi(-s)_f\ge0$ by homogeneity and continuity. In the successor
$S$ is a single based monomial with these frozen exponents.
It is therefore regular, has height one, and satisfies $ES=1$.
Since raising is injective on the labels of height one, these
properties determine $S$ among the distinguished lifts.

Take $w_j=v^{j\Lambda(p,s)}L_p S^j$. Since $S^j$ is a based
monomial in the successor, Proposition~\ref{I:multiplication} gives
leading term $b_j$ and lower coefficients in $u\mathbb Z[u]$.
The product is regular, so target spanning gives $w_j-b_j\in
u\mathcal L_P$. Since $h_i(S)=d_i$ and $ES=1$, the skew Leibniz
rule gives $E(S^j)=[j]_q S^{j-1}$. Also $EL_p=0$, and
$h_i(p)=\Lambda(p,s)$ because $p_i=0$. The skew Leibniz rule therefore gives $Ew_j=[j]_q w_{j-1}$, with
the required normalization $v^{(j-1)\Lambda(p,s)}$. Thus $w_j$ depends
only on $L_p$, $S$, and the grading $h_i(p)$. For a finite
string the last lift and its successive divided powers determine all
the $w_j$. This proves independence of the optimizing seed in both
cases. On an infinite string we have
$E^{(j)}w_j=L_p=E^{(j)}b_j$, so the difference has height at most
$j-1$.

In both cases the change $b_j\mapsto w_j$ is the identity modulo
$u$ and its difference from the identity strictly lowers height.
The change of basis thus has an inverse given on each vector by a
finite geometric sum. The string vectors form a basis of the original
lattice, and their two shifts induce the asserted arrows on its
residues.
\end{proof}

\begin{corollary}\label{B:string-lengths}
Fix one boundary with shift $s$, and let $S$ be the distinguished
lift of degree $-s$ in its optimizing seed. For a
source we use the opposite seed. The following
conditions are equivalent.
\begin{enumerate}
\item[\textup{(i)}] Some string for this boundary has infinite length.
\item[\textup{(ii)}] Every string for this boundary has infinite length.
\item[\textup{(iii)}] The element $S$ belongs to $\mathcal U_{R,P}$.
\end{enumerate}
Condition \textup{(iii)} can be checked on the frozen exponents of
the single monomial $S$ in the successor of the optimizing seed.
When these conditions hold, the string vectors also form an $R$-basis
of $\mathcal U_{R,P}$, and both shifts preserve the integral form.
\end{corollary}
\begin{proof}
An infinite string implies the regularity of $S$ by the proof of
Theorem~\ref{B:strings}. Conversely, suppose $S$ is regular and
let $p$ have height zero. For every $j\ge0$ the product $L_p S^j$
belongs to $\mathcal U_{R,P}$. Its pointed degree is $p-js$, so
$L_{p-js}$ occurs in its expansion in the localized upper algebra. By target spanning this label is regular. Thus the string through $p$ contains a label of every height. Finally, the unit
has height zero, so there is at least one string. In the infinite
case all $w_j$ are normalized products $L_p S^j$. Their
coefficients in the distinguished basis belong to $R$, and the
difference $w_j-L_{p-js}$ has smaller height. Induction on height
therefore inverts the change of basis over $R$.
\end{proof}

When the positive action is jointly locally nilpotent, successive
maximal raises determine both the joint highest label and the candidate
upper global vector in the copy indexed by that label
(Sections~\ref{G:positive}--\ref{S:section}).
For the full quantum group, Section~\ref{X:extension} compares the
length of each string with the prescribed coroot weight at its highest
label to decide whether the positive action has an integrable extension.
\section{The positive action on the distinguished lifts}\label{G:positive}

Theorem~\ref{B:strings} describes the crystal arrows on the residues
of the regular lifts when these lifts span. We now ask when a
positive-module embedding identifies the lifts themselves with upper
global vectors in copies of the restricted dual of the negative half.
Choose a nonempty set $I$ of polynomial frozen vertices, each admitting
a sink-optimizing seed. Their boundary operators preserve
$M_R=\mathcal U_{R,P}$ and give an action of $U_v^+(C)$ on
$M=k\otimes_R M_R$, where
$C=(c_{ij})$ is the boundary Cartan matrix (Theorem~\ref{B:serre}).
For sources we use the sink operators in the opposite seed.

The joint local nilpotence test below applies to every polynomial target.
For the basis comparison we assume that the regular lifts span $M_R$.
Corollary~\ref{I:all-colors} gives spanning when $I=P$. In general it
suffices to check spanning over $k$ at the divisors in $P\setminus I$
(Proposition~\ref{I:spanning-reductions}).
We allow a separate copy of the restricted dual for each joint highest
lift. Under joint local nilpotence, successive maximal divided powers
determine both the copy and the upper global vector to which each
$L_g$ must map.

Let $Q=\bigoplus_{i\in I}\mathbb Z\alpha_i$ be the root lattice of $C$,
and write $Q_+=\sum_i\mathbb N\alpha_i$. The height of
$\sum_i n_i\alpha_i\in Q_+$ is $\sum_i n_i$.
The exchange-degree shift $\gamma_i$ of $E_i$ defines a homomorphism
$\gamma:Q\to\Gamma$ by $\gamma(\alpha_i)=\gamma_i$.
The positive action requires no choice of the coroot weights used in
Section~\ref{X:extension}.

To describe the action on the dual, take $U^-=U_u^-(C)$, with simple generators
$f_i$. Its crystal $B(\infty)$ indexes Kashiwara's lower global
basis $\{G^{\mathrm{low}}(b):b\in B(\infty)\}$ \cite{Kashiwara}. We
grade $U^-$ by \emph{root content}, assigning content $\alpha_i$ to
$f_i$, and denote the content of $G^{\mathrm{low}}(b)$ by $\nu_b$.
The basis dual to the lower global basis is the upper global basis
$\{a_b:b\in B(\infty)\}$ of the \emph{restricted graded dual}
$\mathcal A=\bigoplus_{\nu\in Q_+}(U^-_\nu)^*$. The simple generator $E_i$ acts by
$(E_i\varphi)(y)=\varphi(f_i y)$. Transposition reverses the
monomials in each Serre relation. Reversing their order in the sum
recovers the same relation up to an overall sign, and replacing $v$
by $u=v^{-1}$ fixes its quantum binomial coefficients.
The bar is $\overline\varphi(y)=\overline{\varphi(\overline y)}$,
and the $R$-span and $A_0$-span of the $a_b$ are its integral
form and upper crystal lattice, respectively. These conventions agree with the divided-power formulas of \cite[\S\S 2.2 and
3.2]{Kimura}.

A \emph{coordinate realization} of the distinguished lifts is an
injective $U_v^+(C)$-module map from $M$ to a direct sum of copies
of $\mathcal A$ that sends each $L_g$ to a basis vector $a_b$ in
one copy. When the regular lifts span $M_R$, their images form a
basis of the image with its induced integral form, bar and $A_0$-lattice. Theorem~\ref{S:shifts} determines when both ambient
modified Kashiwara operators of a given color preserve this submodule.

\subsection{Joint local nilpotence}

On the restricted dual, every application of a generator lowers root
height by one. A necessary condition for a coordinate realization is
therefore that sufficiently long words in the $E_i$ kill each
vector, even when the colors vary. For $x\ne0$, let $N(x)$ be the
supremum of the lengths of the words acting nontrivially on $x$,
including the empty word. The action is \emph{jointly locally
nilpotent} if $N(x)<\infty$ for every such $x$.
Proposition~\ref{B:operators} establishes local nilpotence in each
color separately. To test joint local nilpotence we will use only the
centered product rule and the embedding of the target in a seed torus.

Choose an order on $I$ and repeat it periodically. On any initial
segment $\boldsymbol i=(i_1,\ldots,i_h)$ and for $t\in\mathbb N^h$,
write $E_{\boldsymbol i}^{(t)}
=E_{i_h}^{(t_h)}\cdots E_{i_1}^{(t_1)}$, allowing zero exponents.
Write $\boldsymbol i_s$ for the first $s$ periods and
$n_s(x)=\max\{|t|:E_{\boldsymbol i_s}^{(t)}x\ne0\}$,
where $|t|=\sum_jt_j$. Every word occurs as a subsequence of some
$\boldsymbol i_s$, so $N(x)=\sup_s n_s(x)$. For a fixed segment,
the lexicographically greatest tuple with nonzero image is the
\emph{string parameter} of $x$ along that segment.
Its entries record the maximal divided powers taken at successive
letters. Let $Z$ be the set of unchosen frozen vertices, $m=|V|$ the
rank of a seed torus, and $d=m-|Z|$. Write $\operatorname{top}(x)$ for
the result of these successive maximal raises along $\boldsymbol i_d$.
We call a vector \emph{joint highest} if it belongs to the
\emph{joint kernel} $H=\bigcap_{i\in I}\ker E_i$.

\begin{proposition}\label{G:period}
For $x\ne0$, we have $N(x)<\infty$ if and only if every word of length
$n_d(x)+1$ kills $x$. In that case $N(x)=n_d(x)$ and
$\operatorname{top}(x)\in H$. It suffices to check the condition on
the variables of one extended cluster to decide joint local nilpotence on $M$.
The conclusions hold with any nonnegative integer
$d\ge\lfloor m-\operatorname{GKdim}H\rfloor$, where
$\operatorname{GKdim}$ denotes Gelfand--Kirillov dimension.
In particular, one may use $d=m-|J|$ for any set $J$ of
joint highest variables in a single extended cluster.
\end{proposition}

\begin{proof}
For homogeneous $x,y$ the centered product rule gives
\begin{equation}\label{G:leibniz}
 E_i^{(r)}(xy)=\sum_{a+b=r}
 v^{a h_i(y)-b h_i(x)}(E_i^{(a)}x)(E_i^{(b)}y).
\end{equation}
Indeed $E_i$ decreases $h_i$ by $d_i$, and the induction step
uses $q_i^{-b}[a]_{q_i}+q_i^a[b]_{q_i}=[a+b]_{q_i}$.
For a fixed finite word, local nilpotence at each successive letter
leaves only finitely many tuples that act nontrivially on a given
vector. We use two orders on these tuples, the lexicographic order
and the order comparing first the total sum and then the tuple
lexicographically.
Iterating \eqref{G:leibniz} expresses
$E_{\boldsymbol i}^{(t)}(xy)$ as a sum of nonzero scalar multiples of $(E_{\boldsymbol
i}^{(a)}x)(E_{\boldsymbol i}^{(b)}y)$, with $a+b=t$. For either order, translation invariance implies that the sum of the
greatest tuples for $x$ and $y$ occurs only once. Its coefficient
is nonzero because the seed torus is a domain. Thus the greatest
tuple for $xy$ is the sum of those for $x$ and $y$. The greatest total sum on finite
words cannot decrease when letters are inserted. Taking the supremum over the periodic initial segments,
which contain every word as a subsequence, gives $N(xy)=N(x)+N(y)$, also for infinite values.

We next bound the number of periods needed when $N(x)$ is finite
and $x$ is homogeneous. In any tuple giving a nonzero image of $x$,
move the nonzero powers, in order, to the earliest available occurrences
of their colors. Their product and total length stay the same, and the
tuple can only increase lexicographically. All nonzero powers now lie
in the first $N(x)$ periods. Hence a greatest tuple exists on the
whole sequence for either order. We extend it by zeros after its last
nonzero entry.

Fix either order and choose homogeneous $x_1,\ldots,x_r$ of finite
$N$ with linearly independent greatest tuples over $\mathbb Q$.
The ordered products $x^b=\prod_j x_j^{b_j}$, $0\le b_j\le l$,
have distinct greatest tuples. We can multiply each of them by
independent elements of $H$ without losing linear independence.
Indeed, let $T$ be a finite-dimensional graded subspace of $H$
containing $1$ and let $y$ range over a homogeneous basis of
$T^l$. In a relation among the $yx^b$, apply the divided word for
the greatest tuple present. Only one exponent vector $b$ survives,
and its images are nonzero scalar multiples of the $y$ times a
common nonzero vector, by \eqref{G:leibniz}. They are independent
since the torus is a domain.

In a fixed seed, the Laurent supports of these products lie in a
box in $\mathbb Z^m$ of side lengths $O(l)$. Hence
$(l+1)^r\dim_k T^l=O(l^m)$. Every finite-dimensional subspace
of $H$ is contained in such a graded $T$, since each $E_i$
is homogeneous. The definition of Gelfand--Kirillov dimension gives
$r\le m-\operatorname{GKdim}H$. A compatible set $J$ of joint
highest variables has independent ordered monomials in its seed, so
$\operatorname{GKdim}H\ge |J|$. This applies in particular to
$J=Z$.

Suppose the last nonzero entry of a greatest tuple occurs in period
$r$. For $0\le j<r$, apply the prescribed powers through the
first $j$ periods. The greatest tuple for the resulting vector
is the remaining suffix, since a greater suffix would improve the
original tuple. These $r$ suffixes end in different periods and
are therefore linearly independent. The dimension bound gives
$r\le\lfloor m-\operatorname{GKdim}H\rfloor\le d$. For the order by total sum this yields $N(x)=n_d(x)$.
For lexicographic order it shows that maximal raising reaches the
joint kernel within $d$ periods.

For an inhomogeneous $x$ of finite $N(x)$, a fixed word shifts
all exchange degrees by the same amount. It therefore survives on
$x$ exactly when it survives on some homogeneous component. For
either order the greatest tuple for $x$ is the greatest of the
tuples for its components, so the same bound applies. Conversely,
if every word of length $n_d(x)+1$ kills $x$, so does every longer word.

If the variables in one extended cluster have finite $N$, every
homogeneous $x$ has finite $N$ by clearing its Laurent denominators
with a monomial and applying product additivity. The homogeneous components
of an arbitrary element give the same conclusion.
\end{proof}

\subsection{Evaluation on the joint kernel}

Each copy of the restricted dual has a one-dimensional joint kernel.
Thus the joint highest lifts distinguish the copies that occur in a
coordinate realization. For the rest of this section, assume that
the regular lifts span $M_R$ and that the positive action is jointly
locally nilpotent. By Lemma~\ref{B:filtrations}, the lifts
$h_p=L_p$ of height zero in every color form an $R$-basis of
$H\cap M_R$.

Take the maximal divided power at each letter of $\boldsymbol i_d$
on $L_g$. Let $t(g)$ record these powers. Each intermediate
vector is a distinguished lift by Theorem~\ref{B:biperfect}, and the
last is a joint highest lift $h_{p(g)}$ by Proposition~\ref{G:period}.
Let $\widetilde e_i,\widetilde f_i$ denote the crystal operators on
$B(\infty)$, and let $\varepsilon_i(b)$ be the greatest $r$
for which $\widetilde e_i^{\,r}b\ne0$. In the restricted dual, if
$\varepsilon_i(b)=r$, then $E_i^{(r)}a_b=a_{\widetilde e_i^{\,r}b}$ and
$E_i^{(r+1)}a_b=0$. This follows by transposing the multiplication formula for
$f_i^{(r)}G^{\mathrm{low}}(b')$ \cite[Theorem~2.27]{Kimura}. Distinct indices of fixed
$i$-height have distinct images under their maximal raise. Let $b_\circ$ be the crystal index
of $1\in U^-$.

To distinguish vectors with the same joint highest lift we record
the highest components of all their raises. For $y\in U^-$, write
$\partial_y$ for the operator on $M$ obtained by sending $f_i$
to $E_i$ and reversing the order of multiplication. The Serre relations make this a well-defined
$k$-linear antihomomorphism.
Let $\rho_0:M\to H$ fix the joint highest elements $h_p$ and send
every other $L_g$ to zero. We view $\mathcal A\otimes H$ as the
space of $H$-valued linear maps on $U^-$ supported on finitely
many homogeneous components. Define $J:M\to\mathcal A\otimes H$ by
\begin{equation}\label{G:evaluation}
 J(x)(y)=\rho_0(\partial_y x)\qquad(y\in U^-).
\end{equation}
Only contents of height at most $N(x)$ contribute. Evaluating at
$f_i y$ gives $J(E_i x)(y)$, so $J$ intertwines the positive
action. It is the unique such map whose value at $1$ is $\rho_0$.
Indeed, for any intertwiner $T$ we have $T(x)(y)=T(\partial_y x)(1)$.
In particular $J(h)=a_{b_\circ}\otimes h$ for $h\in H$.
For $x\ne0$, choose a word $y$ of length $N(x)$ such that
$\partial_y x\ne0$. Its value belongs to $H$, so
$J(x)(y)=\partial_y x\ne0$. Thus $J$ is injective.
It also commutes with the bars. Indeed $\rho_0$ is projection in
a bar-invariant basis and each $E_i$ commutes with the bar. Finally
$J(M_R)$ lies in the integral restricted dual tensored with
$H\cap M_R$, since every word in divided powers preserves $M_R$
and $\rho_0(M_R)=H\cap M_R$.

\begin{lemma}\label{G:indices}
For each $g$ there is a unique $b(g)\in B(\infty)$ with string parameter $t(g)$ along
$\boldsymbol i_d$ and final index $b_\circ$. The pair $(b(g),p(g))$ determines $g$.
\end{lemma}
\begin{proof}
Expand $J(L_g)$ in the $a_b$ with coefficients in $H$. Raising
maximally along the chosen word kills all summands except those
having the lexicographically greatest tuple. It preserves the
coefficients of the surviving summands. At
each step the crystal raising map is injective on indices of the given height. Since $J$ is an
injective intertwiner, this tuple is $t(g)$ and the result is $a_{b_\circ}\otimes h_{p(g)}$.
The terminal index and the tuple reconstruct the initial index by successive inverse crystal
arrows, so the leading index is unique. To recover $g$, undo its
maximal raises in reverse order, subtracting the recorded multiple of
the shift in each optimizing chart. This is the inverse on tropical
labels of \cite[Proposition~10.1]{Fei}.
\end{proof}

Define $T:M\to\mathcal A\otimes H$ by
$T(L_g)=a_{b(g)}\otimes h_{p(g)}$. Lemma~\ref{G:indices} gives
injectivity. The image of $L_g$ has exchange degree
$\deg h_{p(g)}-\gamma(\nu_{b(g)})$, since its maximal raises
have total content $\nu_{b(g)}$.

\begin{proposition}\label{G:prescribed}
A coordinate realization exists if and only if $T(E_iL_g)=E_iT(L_g)$
for every $i\in I$ and $g\in\mathcal B_P$. With the copies indexed
by their joint highest lifts, the realization is the evaluation map
$J=T$. It identifies each boundary raising arrow with the corresponding
crystal arrow in $B(\infty)$.
\end{proposition}
\begin{proof}
Index the copies by their joint highest lifts. By Lemma~\ref{G:indices},
successive maximal raises force any realization to agree with $T$ on
each $L_g$. Conversely, the stated identities make $T$ an intertwiner
because the $E_i$ generate the positive half. Evaluation of $T$ at $1$
is $\rho_0$, so the uniqueness in \eqref{G:evaluation} gives $T=J$.

Since $T$ is an injective intertwiner, $\rho_i(g)=\varepsilon_i(b(g))$. When this
height is positive, compare the unique terms of height $\rho_i(g)-1$
in $T(E_iL_g)$ and $E_iT(L_g)$ using Theorem~\ref{B:biperfect} and the
transposed multiplication formula \cite[Theorem~2.27]{Kimura}. This
gives $b(r_i g)=\widetilde e_i b(g)$ and $p(r_i g)=p(g)$.
\end{proof}

\section{Recognition of a quantum Schubert cell}\label{C:section}

The joint kernel of the raising operators in the restricted dual is
spanned by the functional of root content zero. Suppose the boundary
action is jointly locally nilpotent and its joint kernel consists of
the constants. Evaluation then embeds the algebra in a single copy
of the dual. We identify its image as a quantum Schubert cell and its
regular lifts as the upper global vectors in that cell. Choose a nonempty set $I$ of polynomial frozen vertices
admitting sink-optimizing seeds, and allow any polynomial frozen set $P\supseteq I$. As before,
$M=k\otimes_R\mathcal U_{R,P}$ and $m=|V|$. To specify the Weyl
group element of the cell, we use a realization of the boundary
Cartan datum with linearly independent simple roots and fundamental
weights $\varpi_i$ dual to the simple coroots. Write $\mathsf{s}_i$ for the reflection in
$\alpha_i$ and $\ell(w)$ for the length of $w$ in these simple
reflections.

\subsection{A one-dimensional joint kernel}

\begin{theorem}\label{C:cell}
Suppose that the joint kernel of the $E_i$ is $k1$ and that their
action is jointly locally nilpotent. Then every frozen vertex belongs
to $I$, and $\gamma:Q\to\Gamma$ is an isomorphism. After a twist
of multiplication by the root grading, evaluation
\eqref{G:evaluation} at the constants identifies $M$ with the quantum
Schubert cell of a Weyl group element $w$ of length $m$. Every
simple reflection occurs in a reduced expression for $w$, and $w$
is determined by
\[ \varpi_i-w\varpi_i=-\gamma^{-1}(\deg X_i)\qquad(i\in I). \]
Under evaluation, the regular distinguished lifts $L_g$ are exactly
the upper global vectors of the cell. Thus evaluation identifies
the integral forms, bars and crystal lattices.
Conversely, any embedding of positive modules into a single copy of
the restricted dual implies both hypotheses.
\end{theorem}

The flag minors of the Schubert cell remain polynomial. We use the
coideal convention in which the transposes of left multiplication
act on the cell. To recognize the algebra, we show that the image of
evaluation is a graded coideal subalgebra of finite growth. To
recognize its basis, we find a cluster whose variables map to chamber
minors. Inverting the flag minors in this cluster lets us compare the
two bases. Regularity along the same minors then singles out the
upper global vectors of the polynomial cell.

The hypothesis on the joint kernel already forces every frozen
vertex to belong to $I$. Indeed, each unchosen frozen variable
is transcendental over $k$, and \eqref{B:monomials} in the respective
optimizing seeds shows that all the $E_i$ kill it. For the rest of the
proof we may therefore take every frozen vertex to be polynomial and
chosen. Corollary~\ref{I:all-colors} then gives target spanning and a
rational polyhedral cone $\mathcal C$ of regular labels.

\begin{lemma}\label{C:bound}
Under the hypotheses of Theorem~\ref{C:cell}, the root degree map
$\gamma:Q\to\Gamma$ is an isomorphism. Every regular lift reaches
$1$ by successive maximal divided raises, of total root content
$-\gamma^{-1}(\deg L_g)$. The degree-zero subspace of $M$ is $k1$.
\end{lemma}
\begin{proof}
By Proposition~\ref{G:period} and Theorem~\ref{B:biperfect}, successive
maximal divided raises take $L_g$ to a distinguished lift in $k1$.
Its pointed coefficient is one, so this lift is $1$. If the raises
have total content $\beta\in Q_+$, then $\deg L_g=-\gamma(\beta)$.
In particular the image of $\gamma$ contains the degrees of the
initial variables, which generate $\Gamma$. Compatibility makes
the mutable exchange rows linearly independent. Since all frozen
vertices are chosen, $\Gamma$ and $Q$ both have rank $|I|$.
The kernel of this surjection has rank zero in the free abelian group $Q$,
so $\gamma$ is an isomorphism. A lift of degree zero has
$\beta=0$ and is already $1$. Since the lifts are homogeneous and span $M$, the entire degree-zero
subspace is $k1$.
\end{proof}

When every frozen vertex is chosen, we can test both hypotheses of
Theorem~\ref{C:cell} on the regular cone $\mathcal C$. By
Lemma~\ref{B:filtrations}, the joint kernel is $k1$ precisely when
the simple heights vanish simultaneously only at the origin.
To check this, take a common rational subdivision of $\mathcal C$
on which the heights are linear. On each cone of the subdivision
their common zero set is rational polyhedral, so a nonzero solution
gives a nonzero integral one.

Assume this kernel condition. Any real linear functional $\lambda$
on $\Gamma_{\mathbb R}$ with $\lambda(\gamma_i)\ge1$ for all $i$
and $\lambda(\deg g)\le0$ on $\mathcal C$ gives joint local nilpotence.
Indeed, every nonzero homogeneous vector has the degree of a regular
lift, and each raise increases $\lambda$ by at least one. Conversely,
joint local nilpotence forces $\gamma$ to be an isomorphism by
Lemma~\ref{C:bound}. A longest nonvanishing word on a homogeneous
$x\ne0$ takes it to a scalar, so its length is
$N(x)=-\operatorname{ht}(\gamma^{-1}(\deg x))$.
Choose integral vectors generating $\mathcal C$ as a real cone.
It therefore suffices to test whether $\gamma$ is an isomorphism
and whether $\operatorname{ht}(\gamma^{-1}(\deg g))\le0$ for each
of these generators $g$.
When the test holds, the equations $\omega_i(\gamma_j)=c_{ij}$
also determine the adapted integral coroots.

\subsection{The coideal and its twisted multiplication}

The evaluation image will be a graded coideal subalgebra of a
quantum half. We use the reflection theorem for Nichols algebras
to show that such a subalgebra is a Schubert cell whenever its
Gelfand--Kirillov dimension is finite. We work in characteristic zero with indeterminate
quantum parameter. Write $(\alpha_i,\alpha_j)=d_i c_{ij}$ for the
symmetric form on $Q$, and $\Phi^+$ and $\Phi^-$ for the positive
and negative real roots. The \emph{inversion set}
$\Phi(w)=\Phi^+\cap w\Phi^-$ has cardinality $\ell(w)$.

\begin{proposition}\label{C:coideal}
Let $A$ be a unital, root-graded one-sided coideal subalgebra of a
symmetrizable quantum half. If $A$ has finite Gelfand--Kirillov
dimension, then $A$ is the quantum Schubert cell of a unique Weyl
group element $w$, in the corresponding coideal convention, and
$\operatorname{GKdim}A=\ell(w)$.
\end{proposition}
\begin{proof}
Let $U$ denote the quantum half and $r$ its braided coproduct.
We first treat right coideals, so $r(A)\subseteq A\otimes U$.
For the Cartan braiding, the quantum half is the Nichols algebra of
the span of its simple generators. Indeed Lusztig's nondegenerate form is adjoint to
the coproduct, so a primitive element of height greater than one
would be orthogonal to all words of its degree. A nonzero homogeneous
Hopf ideal disjoint from degree one would contain such a primitive
element in its least positive degree. Moreover, a root-graded
subspace is stable under the grading torus, which acts on each
homogeneous component by a character. Thus the root-graded coideals
are objects of the braided category used in \cite{HS}. A nontrivial
graded coideal $A$ contains a simple generator $x_i$
\cite[Lemma~5.1]{HS}. For $j\ne i$ the braided adjoint string
through $x_j$ has length $-c_{ij}$. The quantum Serre relation
kills its next term, and no defining relation occurs in degree
$r\alpha_i+\alpha_j$ for $0\le r\le-c_{ij}$. Thus the Nichols
reflections act on degrees as the simple Weyl reflections. The
reflected braiding has the same matrix by Weyl invariance of the
symmetric form. Let $F$ be the intersection of $A$ with
the coinvariants for projection to $k[x_i]$. Multiplication gives
a root-graded vector space isomorphism $F\otimes k[x_i]\simeq A$
\cite[Lemma~2.2]{HS}. The reflection identifies $F$ with the
reflected coideal $A'$ as an algebra, carrying $F_\beta$ to
$A'_{\mathsf{s}_i\beta}$. In particular $A'_{\alpha_i}=0$, since
$F_{-\alpha_i}=0$ \cite[Theorem~5.6 and~(4.6)]{HS}. For any finite-dimensional
subspace $V\subset F$ containing $1$, the subspaces
$V^n x_i^j$, $0\le j\le n$, are independent and contained in
$(V+kx_i)^{2n}$. Taking growth exponents and then the supremum over
$V$ gives $\operatorname{GKdim} A\ge1+\operatorname{GKdim} A'$.
Successive reflections therefore reach $k$ after finitely many steps.

Reading the reflections in reverse reconstructs $A$ from $k$
\cite[Definition~6.3]{HS}. We show by induction, starting with the
empty word, that the resulting word for $w$ is reduced and that the
coideal has Hilbert series $\prod_{\beta\in\Phi(w)}(1-e^\beta)^{-1}$.
Suppose this holds for the reflected coideal $A'$ and the element $w'$.
Since $A'_{\alpha_i}=0$, we have $\alpha_i\notin\Phi(w')$ and
$\ell(\mathsf{s}_iw')=\ell(w')+1$. In the graded factorization of $A$,
the factor $F$ has degrees obtained from those of $A'$ by
$\mathsf{s}_i$, and $k[x_i]$ contributes the factor
$(1-e^{\alpha_i})^{-1}$. The resulting root factors are indexed by
$\{\alpha_i\}\cup \mathsf{s}_i\Phi(w')=\Phi(\mathsf{s}_iw')$.
For a coideal reconstructed from a reduced word in this way, its
Hilbert series distinguishes it from every other graded right coideal
\cite[Lemma~6.4]{HS}. The Poincar\'e--Birkhoff--Witt (PBW) basis
of the quantum Schubert cell has this series, and its braided
coproduct has the required coideal property by \cite[Propositions~4.37--4.38]{Kimura}.
Thus the two subalgebras coincide. The factors $1-e^\beta$ for
distinct positive real roots are nonassociate irreducibles in
$\Q[Q]$, since real roots are primitive. The Hilbert series
therefore determines the inversion set and hence $w$.
In the filtration by root height, the PBW monomials of height at most
$n$ are indexed by tuples $(a_\beta)_{\beta\in\Phi(w)}$ with
$\sum_\beta a_\beta\operatorname{ht}(\beta)\le n$. Their number grows
on the order of $n^{\ell(w)}$. The finitely many root vectors have
positive heights and generate the cell, so this filtration gives
$\operatorname{GKdim}A=\ell(w)$. Word reversal exchanges the two
coproduct conventions.
\end{proof}

We now return to the cluster algebra under the hypotheses of
Theorem~\ref{C:cell}. By Lemma~\ref{C:bound}, assigning content
$\beta=-\gamma^{-1}(\deg x)$ to each nonzero homogeneous $x$
gives $M$ a $Q_+$-grading whose degree-zero part is $k1$. Thus
$\rho_0$ in \eqref{G:evaluation} is projection onto content zero.
We identify $\mathcal A\otimes k1$ with $\mathcal A$. For the algebra
comparison we use $U_v^-$ to realize this dual. Replacing
$u=v^{-1}$ by $v$ preserves the divided-power presentation under
the linear identification of words in the $f_i$. It also preserves
each lower global vector under this identification, by bar invariance.
Consequently the paired vectors $a_b$ and the transposed left
actions are the ones already fixed in Section~\ref{G:positive}.

\begin{lemma}\label{C:twist}
After twisting multiplication on $M$ by the root grading, evaluation
identifies $M$, via Kashiwara's form, with the graded coideal
$J(M)=U_v^-(w,+1)$ for a unique Weyl group element $w$. Its length is $m$, and
$J(M_R)=J(M)\cap\bigoplus_b R a_b$.
\end{lemma}
\begin{proof}
We choose the twist so that the boundary product rule agrees with
the one on the restricted dual. Define the bilinear form $\chi$ by
$\chi(\alpha_i,\alpha_j)=h_i(\gamma_j)$. The boundary Cartan
calculation in the proof of Theorem~\ref{B:serre} gives
$\chi+\chi^{\mathsf T}=-(\, ,\,)$, including the diagonal identity
$h_i(\gamma_i)=-d_i$. For elements of contents $\beta,\eta$
define $x\star y=v^{\chi(\beta,\eta)}xy$. The centered product rule
then becomes
\[ E_i(x\star y)=(E_i x)\star y+v^{-(\alpha_i,\beta)}x\star E_i y. \]
This is the product rule for the transpose of left multiplication
in the dual of the braided coproduct of $U_v^-$
\cite[\S2.1]{Kimura}. The projection onto content zero is
multiplicative. Induction on total content now shows
$J(x\star y)=J(x)J(y)$. Indeed the difference is killed by
every transpose of left multiplication by $f_i$, and these
operators have common kernel $k1$. Kashiwara's form differs from
Lusztig's coproduct-adjoint form by a character of the root grading
\cite[Lemma~2.12]{Kimura}. It is therefore also adjoint to the
braided coproduct and identifies convolution with multiplication
in $U_v^-$. Contracting the first factor of the coproduct of
$J(x)$ against a word in the $f_i$ gives the transpose of left
multiplication by that word applied to $J(x)$, which belongs to
$J(M)$ by equivariance. Words span every homogeneous component,
and these components are finite-dimensional. Thus the coproduct takes $J(M)$ into
$U_v^-\otimes J(M)$. This is Kimura's $+1$ convention
\cite[Proposition~4.38]{Kimura}.

The growth estimate in the proof of Proposition~\ref{G:period}
gives $\operatorname{GKdim}J(M)\le m$. Indeed, the twist only rescales
products of homogeneous elements and thus preserves their Laurent
supports. The ordered monomials in the $m$ initial variables are
linearly independent, giving the reverse inequality.
Proposition~\ref{C:coideal} now identifies $J(M)$ with the cell
of a unique $w$ of length $m$.

Evaluation \eqref{G:evaluation} takes $M_R$ into the integral
restricted dual. Conversely, let $\sum_g c_gJ(L_g)$ be an integral
functional in the image. Choose a summand with lexicographically
greatest string parameter $t(g)$. By Lemma~\ref{G:indices} this
parameter determines $g$, since every string ends at $1$. Evaluate
the functional on the corresponding word in the divided powers of
the $f_i$. The word kills every summand with smaller parameter
and takes $L_g$ to $1$, so its value is $c_g$. Thus $c_g\in R$.
Subtracting this summand and repeating proves the reverse inclusion.
\end{proof}

\subsection{The chamber seed and the upper global basis}\label{C:chamber-seed}

By Lemma~\ref{C:twist}, evaluation identifies the target integral
form with the integral form of the cell. Kimura's theorem gives a
basis of this form consisting of the ambient upper global vectors
that lie in the cell \cite[Theorem~4.22]{Kimura}. To find their
preimages, we identify a seed in the target whose variables map to
chamber minors.

We will recognize the chamber variables by their extremal weight
vectors. On the integrable highest weight module $V(\varpi_i)$ use
the contravariant form for which $e_j$ and $f_j$ are adjoint and
the highest vector has norm one. Let $v_{t\varpi_i}$ be the extremal
vector obtained from the highest vector by successive divided powers along
a reduced expression for $t$. Each extremal vector also has norm
one, since $e_j^{(r)}f_j^{(r)}$ is the identity at the top of
a $j$-string of length $r$. For $t\le w$ in Bruhat order, define
$D_i(t)(y)=(v_{w\varpi_i},y^*v_{t\varpi_i})$ on the negative half,
where $y^*$ reverses words in the $f_i$. Thus $D_i(t)$
is the image under word reversal of the usual unipotent quantum
minor. The extremal crystal element of weight $t\varpi_i$ lies
in the Demazure crystal $B_w(\varpi_i)$, as one sees by using
a reduced subword for $t$ in a reduced expression for $w$.
The minor is therefore a nonzero upper global vector, and word
reversal preserves the upper basis
\cite[Propositions~3.19, 3.34 and 3.46]{KimuraOya}. With this
convention $E_j$ acts by inserting $f_j$ on $v_{t\varpi_i}$.
For every symmetrizable Cartan datum, Goodearl and Yakimov identify
the integral quantum unipotent subgroup over $\mathbb Z[v^{\pm1/2}]$
with both the cluster algebra and the upper cluster algebra in which
the frozen flag minors are polynomial. Their chamber variables are
normalized quantum minors \cite[Theorem~7.3(a),(b)]{GYIntegral}.

Before comparing chamber variables we must account for the
normalizations in the two realizations of the cell. For a homogeneous
functional $d$, let $d_K$ and $d_L$ denote its representatives
under Kashiwara's and Lusztig's forms, respectively. Choose a weight $\rho$ with
$\langle\alpha_i^\vee,\rho\rangle=1$ for every $i$. For
$\beta=\sum_i n_i\alpha_i$ this gives $(\beta,\rho)=\sum_i d_i n_i$.

\begin{lemma}\label{C:normalization}
Evaluation identifies $M_R$, with the twist of Lemma~\ref{C:twist},
with the integral cell in the restricted dual. After applying
$\Q\otimes_{R,\,v\mapsto1}(-)$, evaluation identifies the
classical polynomial upper cluster algebra with the coordinate
algebra of the split unipotent subgroup.
Over $\Q(v^{1/2})$, the map on functionals of content $\beta$ given by
\[ d\longmapsto v^{-(\beta,\beta)/4-(\beta,\rho)/2}\,d_L \]
identifies the original target product with an alternating twist of
the cell product. It takes each upper global functional of the cell
to the corresponding rescaled dual canonical vector in Qin's normalization.
\end{lemma}
\begin{proof}
By Lemma~\ref{C:twist}, $J(M_R)=J(M)\cap\bigoplus_b R a_b$.
This intersection has as an $R$-basis the upper global functionals
in the cell \cite[Theorem~4.22]{Kimura}. Each takes values in $R$
on the divided-power form. Evaluation at $v=1$ thus gives
functionals on the rational Chevalley form. Over $\mathbb C$ these
give Kimura's isomorphism with the coordinate algebra of the
unipotent subgroup \cite[Theorem~4.39]{Kimura}. Both the subgroup
and the evaluation map are defined over $\Q$, so faithful flatness
gives the isomorphism over $\Q$.

The twist is trivial at $v=1$, and Proposition~\ref{B:fiber}
identifies the specialized target over $\Q$ with the classical
upper algebra with the same polynomial frozen variables.

In the realization by Kashiwara's form, the original target product
has twist $v^{f(\beta,\eta)}$ with $f=-\chi$. Rescaling content
$\beta$ by $v^{-f(\beta,\beta)/2}$ changes this exponent to
$(f(\beta,\eta)-f(\eta,\beta))/2$, as one sees by evaluating on
a homogeneous pair. Here $f(\beta,\beta)=(\beta,\beta)/2$.
For $\beta=\sum_i n_i\alpha_i$, the equality of functionals
$(d_K,-)_K=(d_L,-)_L$ gives
$d_L=\prod_i(1-v^{2d_i})^{n_i}d_K$ \cite[Lemma~2.12]{Kimura}.
This last factor and $v^{-(\beta,\rho)/2}$ are characters of the
root grading, so both preserve twisted multiplication. Combining the three
factors gives the stated map. In weight $-\beta$ it is Qin's rescaling
$\operatorname{cor}$ \cite[end of \S7]{QinDCB}.
\end{proof}

For positive content the factor relating the two bilinear forms
vanishes at $v=1$. We therefore specialize a rescaled upper
vector through its integral functional $d$, whose value on a
specialized divided-power element $y$ is $d(y)|_{v=1}$.
The remaining rescalings are powers of $v^{1/2}$ and specialize
to one. In this integral form the standard cluster variables
generate the polynomial cell \cite[Remark~6.1 and Theorem~7.3(a)]{GYIntegral}.
Their specializations therefore generate its classical coordinate
algebra, whose elements are Laurent in every standard chart.

To distinguish a chamber variable from a product, work in a
classical upper cluster algebra with any choice of polynomial
frozen variables. Fix a seed with nonzero mutable exchange rows.
If the product of two elements of the upper algebra is a monomial
in this seed, both factors are units in its Laurent torus, hence
scalar Laurent monomials. Their exponents at polynomial
frozen vertices are nonnegative. A negative exponent at a mutable
vertex would give a nonunit exchange binomial in the denominator in
the adjacent chart, with a monomial numerator. Thus the only variables
with possibly negative exponents are the inverted frozen ones. In
particular these variables generate the units up to scalars, and
each other seed variable is irreducible.

In the setting of Section~\ref{B:boundary-action}, fix a polynomial
frozen vertex $i$ admitting a sink-optimizing seed. We claim that
the maximal divided raises at $i$ of all based monomials with
nonnegative exponents in one extended cluster are based monomials
in a common seed. For such a monomial $X$, write $x=X|_{v=1}$.
In a chart optimizing $i$, the pointed term has coefficient one
and attains the greatest $i$th exponent $r$. Thus $X$ and $x$ have the same height. Write
$R_i^{(r)}=R_i^r/r!$ for the divided power of the classical boundary
derivation $R_i$. All the classical maximal raises $R_i^{(r)}x$
are monomials with nonnegative exponents in a common seed. From
the chosen optimizing seed, one reaches this common seed without
mutating at any vertex adjacent to $i$ in the optimizing seed
\cite[Theorem~10.2]{Fei}. Theorem~\ref{B:biperfect} gives a
distinguished lift $E_i^{(r)}X$, and \eqref{B:monomials} specializes
it to $R_i^{(r)}x$. Its pointed coefficient survives in the common
seed, so its label there is the exponent of this monomial.
Lemma~\ref{I:predecessor} identifies the lift with the based quantum
cluster monomial of that degree. For a word of maximal divided raises, apply
the same argument in the common seed obtained at the preceding step.
This gives a common seed at each step of the word.

\begin{proposition}\label{C:chamber}
We have $J(X_i)=D_i(1)$ for every color. For any reduced expression
$w=\mathsf{s}_{j_1}\cdots \mathsf{s}_{j_m}$, write $w_{\le a}=\mathsf{s}_{j_1}\cdots \mathsf{s}_{j_a}$
with $w_{\le 0}=1$. The $m$ functionals
$D_{j_a}(w_{\le a-1})$ are the images of the individual variables of
one extended cluster of the target.
\end{proposition}
\begin{proof}
We identify the frozen variables first, using their principal ideals.
Each $X_i$ is normal, and compatibility implies that its commutation
character in every chart factors through the exchange grading. Its principal
ideal consists precisely of the elements of strictly positive
order at $i$. Indeed dividing such an element by $X_i$ preserves
regularity at all polynomial frozen variables and in all seed tori.
The order is additive on nonzero products, so the quotient is a
domain. Normality persists under the twist of Lemma~\ref{C:twist}. To see
that the twisted quotient is still a domain, choose an additive order
on $Q$ and take the largest nonzero homogeneous components of two
nonzero elements. Their twisted product is
a nonzero scalar multiple of their original product. The frozen
variables give distinct prime ideals, as one also sees from their
orders in the initial chart.

The homogeneous prime elements of the cell, up to scalars, are its
fundamental flag minors, indexed by $\operatorname{supp}(w)$
\cite[Corollary~6.6]{GYIntegral}. Passing to the opposite product
by word reversal gives the list $D_i(1)$ in our convention.
There are $|I|$ distinct frozen primes. Thus $w$ has full
support and $J$ matches the two lists. The action determines which minor corresponds to each variable. Only
$E_i$ can act nontrivially on $D_i(1)$, since $f_jv_{\varpi_i}=0$
for $j\ne i$. On the target,
$E_iX_i=M_i\ne0$, so $J(X_i)$ is proportional to $D_i(1)$.
Successive maximal divided raises take both vectors to $1$ with
coefficient one, by Lemma~\ref{G:indices} and its proof. This
determines the scalar, and comparison of contents gives the weight
formula in Theorem~\ref{C:cell}. The sum of the fundamental weights
has trivial Weyl stabilizer, so the formula determines $w$.

We recover the rest of the chamber by induction along the reduced
word, starting at its right end. At $j_a$, apply the maximal divided
power of $E_{j_a}$ to every minor already obtained, and adjoin
$X_{j_a}$. A minor from position $b>a$ has the form $D_{j_b}(t)$
with $t=\mathsf{s}_{j_{a+1}}\cdots \mathsf{s}_{j_{b-1}}$.
The word $\mathsf{s}_{j_a}t$ is reduced. For
$r=\langle\alpha_{j_a}^\vee,t\varpi_{j_b}\rangle\ge0$ the
extremal vector identity gives
$f_{j_a}^{(r)}v_{t\varpi_{j_b}}=v_{\mathsf{s}_{j_a}t\varpi_{j_b}}$, and
the next divided power vanishes. The resulting minor is nonzero
because $\mathsf{s}_{j_a}t\le w$ by the subword criterion. Thus the
maximal raise is exactly $D_{j_b}(\mathsf{s}_{j_a}t)$. Induction gives the stated family of minors.

At each step the maximal images are based monomials of one seed by
the preceding observation. We may adjoin $X_{j_a}$ in that seed
because it is frozen.

To see that these monomials are variables, use the seed of
\cite[Theorem~7.3(b)]{GYIntegral} for the permutation reversing
$1,\ldots,m$. Index its variables by the left endpoints $a$ of
the intervals $[a,m]$.
If $k$ is the last occurrence of $j_a$, then
$w_{\le k}\varpi_{j_a}=w\varpi_{j_a}$, since the remaining
reflections fix $\varpi_{j_a}$. Thus its chamber variable is
precisely the normalized minor associated to $D_{j_a}(w_{\le a-1})$
in our word-reversal convention. Under Lemma~\ref{C:normalization}
these minors specialize to distinct, algebraically independent
variables. Each is irreducible in the classical cell algebra by the
factorization argument above. All frozen
variables are polynomial, so each corresponding cluster monomial in
the target is a single variable. The $m$ minors thus occupy distinct
positions in one extended cluster.
\end{proof}

\begin{proof}[Completion of the proof of Theorem~\ref{C:cell}]
By Lemma~\ref{C:normalization} and Proposition~\ref{C:chamber}, the normalized
target chamber variables agree with the standard chamber minors over
$\Q(v^{1/2})$, with the cell product twisted by an alternating
form. Such a twist preserves based Laurent monomials in the given
generators. Indeed, the scalar in their ordered product cancels
the change in its torus normalization. Each mutable exchange row has
content zero, so normalized mutation and the normalized products in the
triangular-basis axioms are preserved as well. Word reversal sends
the usual seed to its opposite, negating its exchange and commutation
forms. Proposition~\ref{B:opposite-family} transports the distinguished
lifts to this opposite seed.

It remains to recover the mutable exchange rows from the common chamber.
Its frozen variables are the fundamental minors on both sides.
By Lemma~\ref{C:normalization} and the Laurent property above,
every function in the common classical coordinate algebra is
Laurent in every target and standard chart. Fix a mutable variable
$x$ in the common chamber and write $P/x,Q/x$ for its two
neighbors. In the chart containing $y=P/x$, Laurent regularity
of $Q/x=(Q/P)y$ implies that $Q/P$ is Laurent in the unchanged
variables. Interchanging the charts shows that $Q/P$ is a Laurent
unit. Each numerator has minimum exponent zero in every variable,
so that unit is a scalar. Both monomials in each numerator have
coefficient one, giving $P=Q$. Their exponent difference determines
the exchange row at $x$ up to sign. The quantum commutation forms
agree on the common frame after the twists above. At position $k$, compatibility requires the exchange row
to pair negatively with $e_k$ under $\Lambda$. This fixes the
sign, since the alternating twist leaves this pairing unchanged
on a row of content zero. Thus the two exchange rows coincide.

We now invert the frozen flag minors. For each position $a$ in the
reduced word, let $a_{\max}$ be the last occurrence of its letter.
After word reversal and the normalization of Lemma~\ref{C:normalization},
our chamber variables are Qin's interval minors indexed by $[a,a_{\max}]$.
They form his shifted seed $t_0[1]$, which is joined by mutation to
the initial seed $t_0$ \cite[\S7.4 and Theorem~8.2.2]{QinDCB}.
Since the frames and exchange rows agree in the chamber, the reverse
sequence gives the same normalized exchanges in both algebras.
Theorem~\ref{I:localized} therefore lets us compare the lifts in $t_0$.

Kimura's multiplication formula for an extremal minor translates the
PBW parameter of a dual canonical vector \cite[Theorem~6.19]{Kimura}.
Qin proved that the normalized frozen translates of the rescaled dual
canonical vectors form a triangular basis in $t_0$
\cite[(9.3), Proposition~9.3.1 and Theorem~9.3.2]{QinDCB}. By the fixed-seed uniqueness argument in the proof
of Theorem~\ref{I:qin-comparison}, its element of degree $g$ has the
distinguished residue of that degree and is therefore its unique
integral bar-invariant lift $L_g$. In the presence of half-powers
the same argument uses $u^{1/2}$ as local parameter.

Undoing the normalization preserves each homogeneous basis line,
and word reversal permutes the upper global vectors. Thus each
upper global functional of the polynomial cell spans a distinct
localized basis line. These functionals form a basis of the cell by
Kimura's theorem, so their lines are exactly the localized basis
lines contained in the cell. By Lemma~\ref{C:twist}, the images
of the regular lifts span this same space and hence occupy all these
lines. Write $J(L_g)=c\,a_b$ on one of them. By the maximal
divided-power calculation of Lemma~\ref{G:indices}, the leading
coefficient of $J(L_g)$ is one, so $c=1$. These equalities hold
over $k$, giving the claimed identifications of integral forms,
crystal lattices and bars. Conversely, in a single copy of
$\mathcal A$ all sufficiently long words annihilate each vector
and the joint kernel is one-dimensional. Its nonzero submodule
contains a joint highest vector, which proves necessity.
\end{proof}

Theorem~\ref{S:shifts} determines which upper Kashiwara operators
restrict to this cell and to its copies in a coordinate realization.

\section{The frozen-product criterion for the positive half}\label{S:section}

We test the positive action on a single frozen monomial. Choose a
nonempty set $I\subseteq P$ of sink colors and let $F$ be the based
product of their variables with multiplicities $k_i>0$. The remaining
frozen variables, indexed by $Z$, generate a quantum torus
$\mathcal T_k(Z)$ in the localized algebra. Every $E_i$ vanishes
on this torus.

We use the periodic word of Section~\ref{G:positive}. Then $n_s(F)$
is the greatest total length over $s$ periods, and $\operatorname{top}(F)$
is obtained by raising maximally at each letter over $d=|V|-|Z|$
periods. Each maximum exists by local nilpotence of the individual
generator.

\begin{theorem}\label{S:single-product}
The positive module $M$ embeds over $k$ into a direct sum of copies
of the restricted dual by a map taking each regular distinguished
lift $L_g$ to an upper global vector if and only if $N(F)<\infty$
and $U_v^+F$ meets $\mathcal T_k(Z)$ nontrivially. If $Z$ is empty,
the intersection condition is $1\in U_v^+F$. One can choose the
embedding to identify $\mathcal U_{R,P}$, its bar and its common
lattice with those induced on the image. The regular lifts then
map bijectively onto the ambient upper global vectors lying in the
image. These vectors form a basis of the image precisely when the
regular lifts span $\mathcal U_{R,P}$.

The two conditions on $F$ are equivalent to $n_{d+1}(F)=n_d(F)$
and the requirement that $\operatorname{top}(F)$ be a
based Laurent monomial in the variables in $Z$. The same criterion
holds after inverting the variables at $P\setminus I$.
\end{theorem}

\begin{proof}
By Theorem~\ref{S:frozen-product}, joint local nilpotence is
equivalent both to $N(F)<\infty$ and to $n_{d+1}(F)=n_d(F)$.
An embedding forces this property. It holds on the original target
and on the localization at $P\setminus I$ under the same condition,
since the orbit of $F$ is unchanged.

Assume joint local nilpotence and let $A=\operatorname{top}(F)$.
By Theorem~\ref{S:minimal}, the localized target admits a coordinate
realization precisely when $A$ has no mutable factor. Such an $A$
is a frozen monomial killed by every $E_i$, so its exponent at each
chosen vertex is zero by \eqref{B:monomials}. Thus the condition is
$A\in\mathcal T_k(Z)$, or equivalently
$U_v^+F\cap\mathcal T_k(Z)\ne0$ by Corollary~\ref{S:frozen-orbit}.
Proposition~\ref{G:localization} transfers the realization and the
induced structures to $M$. Its regular lifts form a basis of the
image exactly when they span over $k$, or equivalently over $R$ by
Proposition~\ref{I:spanning-reductions}.
\end{proof}

\subsection{Detection of joint local nilpotence}

\begin{theorem}\label{S:frozen-product}
The following conditions are equivalent.
\begin{enumerate}
\item The positive action on $M$ is jointly locally nilpotent.
\item $N(F)<\infty$.
\item $n_{d+1}(F)=n_d(F)$.
\end{enumerate}
If $n_{s+1}(F)=n_s(F)$ for some $s\ge0$, then $N(x)=n_s(x)$
for every nonzero $x\in M$. Whenever $n_{s+1}(F)>n_s(F)$, the
difference is at least $\min_{i\in I}k_i$.
The equivalence still holds with $d=m-|J|$ for any set $J$ of
joint highest variables in one extended cluster.
\end{theorem}

\begin{corollary}\label{S:tropical-test}
Compute $n_d(F)$ and $n_{d+1}(F)$ by the backward covector calculation
below. Starting from the exponent of $F$, read $\boldsymbol i_d$ from
left to right. At each letter $i$, transport the current degree by
tropical mutation to a seed where $i$ is a sink and replace the
resulting exponent $a$ by $a+a_i s_i$. Return the final degree to the
initial seed. The embedding of Theorem~\ref{S:single-product} exists
precisely when $n_{d+1}(F)=n_d(F)$ and this degree is supported on $Z$.
For fixed valued exchange data and colors, the test is independent
of the compatible quantum form. Both calculations also apply to the
ordinary divided boundary derivatives at $v=1$.
\end{corollary}

\begin{proof}[Proof of Theorem~\ref{S:frozen-product} and Corollary~\ref{S:tropical-test}]
Let $d_\ell(x)$ be the maximum of $b\ell$ over the exponents $b$
in the Laurent support of $x$ in a given seed. Since
$d_\ell=-\nu_{-\ell}$ for the trivial scalar valuation,
Lemma~\ref{I:valuation-mutation} transports this maximum across
mutation at $j$ by replacing $\ell_j$ with
$\max\{[B_{j,V}]_+\ell,[-B_{j,V}]_+\ell\}-\ell_j$.
In a seed optimizing $i$, the exponents surviving $E_i^{(r)}$
are exactly those with $b_i\ge r$ and they have distinct images
$b+r s_i$. Hence
\[ \max_{r:\,E_i^{(r)}x\ne0}\{r+d_\ell(E_i^{(r)}x)\}
   =d_{\ell+\delta e_i^*}(x),\qquad \delta=[1+s_i\ell]_+. \]
Starting with the zero covector, read $\boldsymbol i_s$ backwards.
At each letter $i$, transport the covector to a seed optimizing $i$
and add $\delta e_i^*$. The resulting covector $\ell_s$ in the
initial seed satisfies $n_s(x)=d_{\ell_s}(x)$.

Transport preserves the value on each frozen variable. One further
period therefore increases $n_s(F)$ by the sum of the $k_i\delta$
over its letters. Each term is nonnegative. If the sum is positive,
it is at least $\min_i k_i$. If it vanishes, every update is zero,
and transport back to the initial seed returns $\ell_s$ by
Lemma~\ref{I:valuation-mutation}. All subsequent periods fix this
covector, so $N(x)=n_s(x)$ for every $x\ne0$. Conversely, if
$N(F)<\infty$, Proposition~\ref{G:period} gives $N(F)=n_d(F)$,
hence $n_{d+1}(F)=n_d(F)$. The same proposition permits the bound
$m-|J|$ and shows that $\operatorname{top}(x)$ is joint highest
under these equivalent conditions.

For the forward calculation, $F$ is the distinguished lift of its
degree. At a sink $i$ the maximal divided power takes a regular lift
of degree $a$ to the lift of degree $a+a_i s_i$
(Theorem~\ref{B:biperfect}). Thus the final degree is that of
$\operatorname{top}(F)$ and is supported on $Z$ precisely when this
lift is the corresponding frozen Laurent monomial.
Theorem~\ref{S:single-product} now gives the criterion. The backward
calculation uses only exchange exponents and boundary shifts.
In the forward calculation each maximal raise has pointed coefficient
one and therefore keeps its degree and height upon specialization.
This gives the same two calculations at $v=1$.
\end{proof}

\subsection{Localization at joint highest cluster variables}

Inverting the unchosen polynomial variables lets us work in an
algebra in which the regular lifts span. Indeed, the resulting
algebra is $\widehat M_R=\mathcal U_{R,I}$, so
Corollary~\ref{I:all-colors} applies. The following proposition
allows us to test existence of the embedding in this localization.

\begin{proposition}\label{G:localization}
There is a positive-module embedding of $M$ over $k$ that sends
each regular $L_g$ to an upper global vector in a copy of the
restricted dual if and only if the distinguished basis of
$\widehat M_R$ admits a coordinate realization.
\end{proposition}
\begin{proof}
For a frozen exponent $a$ supported on $P\setminus I$, let $\tau_a$
be normalized multiplication by $X^a$. On each exchange-homogeneous
component it acts in every chart by $X^b\mapsto X^{b+a}$, since
$a$ pairs trivially with the mutable exchange rows. Thus
$\tau_a L_g=L_{g+a}$, and $\tau_a$ commutes with every $E_i^{(j)}$
by \eqref{B:monomials} in a seed optimizing $i$.

Suppose the embedding on $M$ exists. A sufficiently large
translation takes any finite family of localized lifts into $M_R$.
Each localized lift therefore has finite raising length, and these
lifts span the localization. Lemma~\ref{G:indices} assigns each
one a standard index $b(g)$ and a terminal highest label $p(g)$.
Translation fixes its string parameter and sends $p(g)$ to $p(g)+a$.

For a fixed $g$ and $i$, translate $L_g$ and all summands of $E_iL_g$
into $M_R$ at once. The embedding sends these translated lifts to
upper global vectors. The string parameter determines the index
in the restricted dual, and the terminal highest lift determines the
copy. Distinct highest lifts give distinct copies by injectivity.
Removing the common translation from the highest labels therefore
gives the intertwining identity of Proposition~\ref{G:prescribed}
for $E_iL_g$. This constructs the localized realization.

Conversely, restrict a localized realization to $M$. Each ambient
upper global vector in the localized image has a single lift as
preimage, and this lift belongs to $M$ exactly when it is regular
along $P$. Since scalar denominators do not affect frozen orders,
$M_R=M\cap\widehat M_R$. Likewise $\mathcal L_P=M\cap\mathcal L_I$,
and the localized bar restricts to $M$. The localized realization
identifies its integral form, bar and lattice by mapping its basis
to ambient upper global vectors. Intersecting with the image of $M$
gives the asserted identifications on the original target.
\end{proof}

To construct the required localizations, assume that the positive
action on $\widehat M_R$ is jointly locally nilpotent.
Proposition~\ref{G:period} gives $\operatorname{top}(x)\in
H=\bigcap_{i\in I}\ker E_i$ for each nonzero $x$. If $x$ is a
regular distinguished lift, its maximal raise is again a regular
distinguished lift by Theorem~\ref{B:biperfect}.

We seek cluster variables whose Laurent monomials span the joint
kernel after localization. Start with an extended cluster
$x_1,\ldots,x_m$ and raise its members maximally
along $\boldsymbol i_d$. By Section~\ref{C:chamber-seed}, the resulting elements
$\operatorname{top}(x_1),\ldots,\operatorname{top}(x_m)$ are
based cluster monomials in a single seed.
Denote by $J$ the set of variables occurring in these $m$
monomials, and choose a seed containing $J$. Each member is joint
highest, by additivity of the simple heights on products. In
particular $Z\subseteq J$, since each unchosen frozen variable is
its own maximal raise. Every member of $J\setminus Z$ is mutable,
since $E_iX_i=M_i\ne0$ for $i\in I$.

We may freeze and invert any compatible set of mutable variables
killed by all the $E_i$. We call such a choice an \emph{admissible
freezing}. Let $K$ consist of $Z$ and the newly frozen variables.
Write $M_R(K)$ for the resulting integral upper algebra and
$M(K)=k\otimes_R M_R(K)$. Write $\mathcal T_R(K)$ for the
based quantum torus generated by the variables in $K$.
Thus $M_R(Z)=\widehat M_R$.

\smallskip\noindent\emph{Reddening after freezing.}\phantomsection\label{S:freeze-reddening}
Take principal coefficients in a seed containing the variables to be
frozen, and let $D$ be the mutable positions that remain. The
successor has $g$-vectors $-e_j$ up to permutation. By the enough
$g$-pairs theorem, a sequence of mutations in $D$ gives a cluster
whose projected $g$-cone contains their projections, hence the
negative orthant in $\mathbb R^D$
\cite[Proposition~4.6 and Theorem~4.8]{CaoLi}. Along this sequence the
untouched $g$-vectors project to zero and the $c$-vectors at positions
in $D$ are supported on $D$. The projected vectors at $D$ are
therefore the $g$-vectors of the restricted principal pattern by
the mutation recurrences. In each coordinate these vectors have the
same weak sign, so their cone equals the negative orthant. They
form an integral basis, hence are the negative coordinate vectors up
to permutation. This gives a successor after freezing.

\begin{lemma}\label{S:invariants}
For every admissible freezing $K$, the original operators extend to
$M_R(K)$, and every chosen color has a seed optimizing it both before
and after freezing. Along any word in the colors, the maximal divided
raises of all based cluster monomials of a given seed of the frozen
pattern are based monomials of a single seed of that pattern.
Here the exponents in $K$ may be arbitrary integers.

The joint kernel in $\widehat M_R$ is
$\widehat M_R\cap\mathcal T_R(J)$. In $M_R(J)$,
the joint kernel of the original operators is $\mathcal T_R(J)$
and the action is jointly locally nilpotent. For each admissible
freezing $K$, the algebra $M_R(K)$ is both the left and the right
localization of $\widehat M_R$ at the variables in $K\setminus Z$.

More generally, take any subset of the variables of a cluster of
$M_R(K)$ whose maximal raises along $\boldsymbol i_d$ are
Laurent monomials in $K$. Every joint highest element of $M_R(K)$
that is Laurent in these variables and $K$ belongs to $\mathcal T_R(K)$.
\end{lemma}

\begin{proof}
First choose a seed containing the joint highest variables to be
frozen. Apply the common-seed construction preceding
Proposition~\ref{C:chamber} in color $i$. Each of these variables
survives its maximal raise. Its classical specialization is therefore
a monomial in the resulting seed. With all original frozen variables
polynomial, the surviving variable is irreducible and all seed
variables are nonunits by the factorization argument in
Section~\ref{C:chamber-seed}. Thus it occurs as a single variable,
also in the quantum seed by its label. The common seed is obtained
from a seed optimizing $i$ without mutating its mutable neighbors,
so it still optimizes $i$. Connectivity of the seeds containing
the chosen variables \cite[Theorem~6.2]{CaoLi} places it in the
frozen pattern.

Since freezing restricts the allowed mutations,
$\widehat M_R\subseteq M_R(K)$. For $x\in M_R(K)$ choose a monomial
$b$ in the newly inverted variables that clears their negative
exponents in a chart containing $K$. Both $bx$ and $xb$ are
Laurent in each adjacent chart away from $K$. At a newly frozen
position they are polynomial in the variable being replaced, so
substitution of its exchange relation gives Laurent regularity in
that adjacent chart as well. The quantum upper bound theorem
\cite[Theorem~5.1]{BZ} places both products in the original localized
upper algebra. Their orders at the polynomial colors are nonnegative
in the chosen chart, so they lie in $\widehat M_R$. This gives both
fraction expressions. The skew product rule extends the original
operators since they kill the denominators. On a homogeneous
fraction, clearing the denominator commutes with a divided-power
word up to a power of $v$. Thus joint local nilpotence passes to
$M(K)$, and the same number of periods reaches its joint kernel.

In a shared optimizing seed, let $c_i$ be the part of $m_i$ on
$K\setminus Z$. The frozen boundary has shift $s_i-c_i$.
Formula~\eqref{B:monomials} gives the same binomial for both operators
and hence
\[ E_i^{(r)}=\tau_{r c_i}(E_i^{\rm fr})^{(r)}, \]
where $\tau_c$ is based translation in the Laurent variables.
Both operators commute with these translations and have the same
height. The frozen pattern admits reddening, so the common-seed
argument of Section~\ref{C:chamber-seed} applies there
\cite[Theorem~10.2]{Fei}. Clearing negative Laurent exponents by
translation and then restoring them proves the assertion about
maximal raises, one letter at a time.

Finally, let $Y$ be a subset of an extended cluster whose variables
have maximal raises in the units of a quantum torus $\mathcal T$
on joint highest variables. If a homogeneous joint highest element
$z$ is Laurent in $Y$, choose a monomial $b$ in $Y$ such that $bz$
is polynomial in $Y$.
In $bz$ the greatest string parameter of its monomials survives.
Indeed, among the monomials with that parameter the terminal
labels are distinct, since each maximal raise of fixed height is
injective on labels. By \eqref{G:leibniz} their terminal monomials
are products of the maximal raises of the variables and lie in
$\mathcal T$. Thus $\operatorname{top}(bz)\in\mathcal T$. The same
product rule gives $\operatorname{top}(bz)=v^c\operatorname{top}(b)z$
for some $c\in\mathbb Z$, and $\operatorname{top}(b)$ is a unit
of $\mathcal T$. Hence $z\in\mathcal T$.

For the original cluster take $\mathcal T=\mathcal T_R(J)$. All
its elements are joint highest, so intersection with $\widehat M_R$
gives the claimed kernel. Clearing denominators in $J$ gives the
kernel after localization at $J$. For the indicated subset of a
cluster of $M_R(K)$, include the variables of $K$ in $Y$ and take
$\mathcal T=\mathcal T_R(K)$ to obtain
the last assertion. The kernels respect the exchange grading, so
the conclusions hold for arbitrary $z$.
\end{proof}

\smallskip\noindent\emph{Deleting Laurent variables.}\phantomsection\label{S:delete-Laurent}
Suppose that a seed admits reddening and its frozen variables consist of the polynomial
colors $I$ and a set $Y$ of invertible joint highest variables.
Allow the original shift at each color to differ from the canonical
shift in an optimizing seed by an exponent in $\Z^Y$, as in
Lemma~\ref{S:invariants}. Suppose that each variable of one extended
cluster has a nonzero iterated divided raise in the torus on $Y$.
We show that the upper algebra and its distinguished lifts decompose
into based Laurent translates of those of the seed with $Y$ deleted.

Project the exponent lattice of the cluster away from $Y$ and
quotient by the integral span of the projected exchange rows.
Write $\bar\Gamma$ for the resulting group and $\eta_i$ for
the projected exchange degree of the original shift $s_i$.
The raise of each cluster variable in the $Y$-torus has degree zero
in $\bar\Gamma$. Thus the $\eta_i$ generate $\bar\Gamma$.
If $n$ coordinates remain after projection, there are $n-|I|$
exchange rows, so $\bar\Gamma$ has rank at least $|I|$.
The surjection $\Z^I\to\bar\Gamma$ is then an isomorphism.
The projected exchange rows together with representatives of the
$\eta_i$ form an integral basis of $\Z^n$.

Consequently there is a unique grading by $\Z^Y$ that assigns to
each element of $Y$ its coordinate degree and vanishes on the
exchange rows and the original shifts. Indeed, these conditions
prescribe the degrees of the remaining coordinates on the integral
basis just obtained. In each chart, replace every variable $x$
outside $Y$ by its degree-zero translate $\tau_{-\deg_Y x}x$.
The two exchange terms have the same degree and project to the two
terms of the exchange after deleting $Y$. The induced skew form
is compatible because the exchange rows pair trivially with the
coordinate vectors at $Y$. Hence in every seed the degree-zero
chart is the chart of the deleted seed. Intersecting these charts,
with polynomial exponents at $I$, recovers its integral upper
algebra. Every other homogeneous component is a based Laurent
translate of degree zero. Theorem~\ref{I:localized} gives identical
normalized coefficients for the two frozen patterns, and every
projected exponent has a unique degree-zero lift. The translates
therefore identify the distinguished lifts, the bars and the common
lattices.

The original operators preserve this grading. In degree zero their
shifts project to the canonical shifts of the deleted seed, so
\eqref{B:monomials} gives precisely its boundary operators. The
original boundary monomial $M_j$ is a based Laurent translate of
the degree-zero monomial, and
translation in $Y$ commutes with each $E_i$. Thus the heights
$\rho_i(M_j)=-c_{ij}$ for $i\ne j$ agree, and the symmetrizers
$d_i$ are unchanged. This identifies the Cartan data.

Once the variables in $J$ are invertible, their exponents separate
the copies of the restricted dual. The degree-zero component has only
scalars in its joint kernel, so the recognition theorem for Schubert
cells applies to it. The growth bound in Proposition~\ref{G:period}
has a geometric interpretation here. It equals the length of the
Weyl group element of this cell. Write $H_M=\bigcap_{i\in I}\ker(E_i|_M)$
for the joint kernel on the original target.

\begin{theorem}\label{S:generic}
The distinguished basis of $M_R(J)$ has a coordinate realization for the original positive
action. Its image consists of copies of a single quantum Schubert
cell indexed by the Laurent exponents in $\Z^J$. Its Weyl group
element has full support and length
$m-|J|=m-\operatorname{GKdim}H_M$.
Conversely, if the
upper algebra obtained by freezing and inverting a compatible set of
joint highest variables embeds as a positive module into a direct sum of copies of the restricted dual,
the original action is jointly locally nilpotent.
\end{theorem}

\begin{proof}
Each member of $J$ is a variable of an extended cluster of the
original seed pattern. In that chart its orders at all polynomial
frozen divisors are nonnegative, so it belongs to $M$.
The joint kernel on $M$ therefore contains the quantum
affine space on these variables. By Lemma~\ref{S:invariants} it is
contained in their quantum torus. Both algebras have
Gelfand--Kirillov dimension $|J|$, as one sees by counting based
monomials. Thus $|J|$ is intrinsic to the original action.

Apply the grading construction above with $Y=J$. The maximal raises
of every variable in a chart containing $J$ belong to its torus
by Lemma~\ref{S:invariants}. The degree-zero algebra is therefore
the upper algebra of the seed with $J$ deleted, which inherits
reddening.
Lemma~\ref{S:invariants} gives joint local nilpotence and a scalar
joint kernel on this algebra. We may now apply Theorem~\ref{C:cell}
to the degree-zero algebra. It identifies its distinguished lifts
with the upper global basis of a cell of full support. Its length
is the number of variables in the deleted chart,
$m-|J|$. Each Laurent translate gives the same realization for
the original action. Conversely, suppose an admissible freezing
$K$ gives the embedding in the theorem. The restricted dual is
jointly locally nilpotent, and this property passes to its direct
sums and submodules. The inclusion of the original target in $M(K)$
intertwines the original operators, so the same property holds there.
\end{proof}

Write $A=\operatorname{top}(F)$ and $A_i=\operatorname{top}(X_i)$.
By \eqref{G:leibniz}, at the maximal order only the product of the
maximal divided powers survives. Thus $A$ is the based product of
the $A_i^{k_i}$ in their common cluster.

\begin{corollary}\label{S:frozen-orbit}
Let $w$ be the Weyl group element of Theorem~\ref{S:generic} and
$\beta$ the total simple-root content of the maximal divided powers
taking $F$ to $A$, and let
$\lambda=\sum_{i\in I} k_i\varpi_i$ in a full weight realization
of the boundary Cartan matrix. Then $w\lambda=\lambda-\beta$,
which determines $w$. The joint highest subspace of $U_v^+F$
is $kA$. In particular $U_v^+F$ meets $\mathcal T_k(Z)$
nontrivially if and only if $A\in\mathcal T_k(Z)$.
\end{corollary}

\begin{proof}
The raises of $F$ stay in its Laurent degree, whose joint kernel is
one-dimensional by Theorem~\ref{S:generic}. They include $A$, so
the joint highest subspace of $U_v^+F$ is $kA$. This gives the
intersection criterion because every element of $\mathcal T_k(Z)$
is joint highest.

In degree zero, the translate of $X_i$ is the fundamental minor
$D_i(1)$ of content $\varpi_i-w\varpi_i$ by
Proposition~\ref{C:chamber}. The translate of $F$ therefore has
content $\lambda-w\lambda$. Translation commutes with the raises
and takes $A$ to a scalar, so $\beta=\lambda-w\lambda$. This
determines $w$ because $\lambda$ is strictly dominant.
\end{proof}

\subsection{The smallest localization}

Inverting all variables of $J$ gives a sum of copies of a Schubert
cell. To find the smallest localization with a coordinate
realization, it suffices to examine the maximal raises of the chosen
frozen variables.

\begin{lemma}\label{S:necessary}
Suppose $K$ defines an admissible freezing and the distinguished
basis of $M_R(K)$ has a coordinate realization. Then $\operatorname{top}(X_i)$ is a Laurent monomial in $K$
for each chosen frozen variable $X_i$.
\end{lemma}

\begin{proof}
The realization assigns to each distinguished lift both its root
content and its copy in the ambient restricted dual. Multiplication
on either side by $X_i$ takes each lift to a
nonzero scalar multiple of another one. String parameters add on
these products by \eqref{G:leibniz}. Their total simple-root
contents therefore add as well, since each string in a copy ends
at its unique highest vector. Thus both multiplication operators
for $X_i$ are homogeneous for root content. For any
exchange-homogeneous $y$ with this property, the product rule
shows that both multipliers for $E_j y$ are homogeneous of
content lower by $\alpha_j$. Indeed $E_j$ lowers content and
the exchange grading operators $\sigma_j^{\pm1}$ are diagonal on
the distinguished basis. Iteration gives root content zero for the
multipliers of $a=\operatorname{top}(X_i)$.

The element $a$ quasi-commutes with every joint highest lift
$h$. Apply the word taking $X_i$ to $a$ to the
quasi-commutation of $X_i$ with $h$. All powers on $h$
other than zero vanish. Write $ah=v^c ha$, and let $x$ be a
distinguished lift of content $\nu=\sum_j n_j\alpha_j$ in the copy
through $h$. For every raising word $D$ of content $\nu$, we have
$D(x)\in kh$. Since $a$ is joint highest, the centered product rule
gives $D(ax)=v^{c-2\sum_j n_j h_j(a)}D(xa)$. The exponent is
independent of the order of the letters in $D$. Both products have
content $\nu$, and evaluations on the words spanning $U^-_\nu$
separate vectors of this content in each copy. Thus $a$ quasi-commutes
with $x$.

It remains to determine which variables can occur in $a$.
By Lemma~\ref{S:invariants}, maximal raises for the original action
give a cluster monomial of the frozen pattern, up to a Laurent
monomial in $K$.
In a chart for this monomial, quasi-commutation
with the adjacent cluster variable at any mutable position $j$
forces its two exchange terms to have the same commutation factor
with $a$. Compatibility gives $d_j a_j=0$, where $a_j$
is the exponent of $a$ there. Thus only frozen variables occur.
In a sink-optimizing seed for a chosen vertex $i$, the $i$-height
of this frozen monomial is its exponent at $X_i$. Since every $E_i$
kills $a$, all these exponents are zero. Its based coefficient is one.
\end{proof}

The factor $A_i=\operatorname{top}(X_i)$ is independent of the
order of the colors. Indeed, fix the Laurent grading of
Theorem~\ref{S:generic}. Every sequence of maximal raises from
$X_i$ stays in its component, whose joint kernel is one-dimensional.
Any two sequences ending in that kernel give proportional
distinguished lifts, and their pointed coefficients are one.
Since every $k_i$ is positive, the mutable factors of
$\operatorname{top}(F)$ are precisely those occurring in the $A_i$.

\begin{theorem}\label{S:minimal}
Let $T\subseteq J$ consist of $Z$ and the mutable factors of
$\operatorname{top}(F)$. For an admissible freezing $K$, the
distinguished basis of $M_R(K)$ has a coordinate realization if and
only if $T\subseteq K$. Its image is a direct sum of copies of the
cell of Theorem~\ref{S:generic}.
\end{theorem}

To prove sufficiency, we must recover the chamber variables after inverting
only $T$. First specialize at $v=1$ in a chart containing $J$.
By Proposition~\ref{B:fiber}, the specializations for $T$ and $J$
are classical upper algebras in the same rational function field.
The former lies in the latter because $T\subseteq J$.
The degree-zero algebra of Theorem~\ref{S:generic} specializes to
the coordinate algebra $D_w$ of the unipotent subgroup over $\Q$,
since the twist of multiplication disappears at $v=1$. We continue to write $X_i$ and $A_i$
for their specializations, so that the fundamental flag minors are
$X_i/A_i$.

\begin{lemma}\label{S:flag-generation}
The fundamental flag minors and all their iterated derivatives under
the simple raising operators generate $D_w$ as a $\Q$-algebra.
\end{lemma}
\begin{proof}
Let $D_0$ be the subalgebra generated by these derivatives. It is
stable under each simple derivation. By Proposition~\ref{C:chamber}, the
chamber minors for a reduced expression of $w$ are obtained from
the fundamental flag minors by maximal divided raises. They belong
to $D_0$ and generate the fraction field of $D_w$.

For $y\in D_w$, let $V_y$ be the span of $y$ and all its
iterated raising derivatives. Each derivation decreases root height
by one, so $V_y$ is finite-dimensional. Since $D_0$ and $D_w$
have the same fraction field, clearing denominators for a basis of
$V_y$ gives a nonzero element of the ideal
$\mathfrak a_y=\{a\in D_0:aV_y\subseteq D_0\}$.
For a simple derivation $\partial$, $a\in\mathfrak a_y$ and
$z\in V_y$, both terms in
$(\partial a)z=\partial(az)-a\,\partial z$ belong to $D_0$.
Thus $\mathfrak a_y$ is stable under all the simple derivations.
Apply a raising word of greatest nonvanishing length to any nonzero
element of this ideal. Its image is a nonzero joint highest element.
In the classical restricted dual such an element is constant, since
it vanishes on every word of positive length in the lowering
generators. Hence $1\in\mathfrak a_y$ and $y\in D_0$.
\end{proof}

\begin{lemma}\label{S:chamber-descent}
For each reduced word for $w$, its chamber minors agree, up to
Laurent monomials in $T$, with distinct variables of one extended
cluster in the classical specialization of $M_R(T)$. Their positions form a set $S$
containing the chosen frozen vertices. Every mutable neighbor of a
vertex in $S$ also lies in $S$.
\end{lemma}
\begin{proof}
Each $A_i$ is a unit in the classical specialization of $M_R(T)$,
so the fundamental flag minors $X_i/A_i$ belong to it. The simple
derivations preserve this algebra, which therefore contains $D_w$ by
Lemma~\ref{S:flag-generation}.

Fix a reduced word. The construction in Proposition~\ref{C:chamber}
successively raises the current family in one color and adjoins
the next fundamental flag minor $X_i/A_i$. This minor is a frozen
variable times a Laurent unit in $T$, so it can be adjoined in
any seed. At each step Lemma~\ref{S:invariants}, specialized at
$v=1$, gives one seed for the maximal raises of the whole family.
Induction places the chamber minors in a common chart of the
specialization of $M_R(T)$, up to Laurent monomials in $T$.

In this chart factor a chamber minor into ordinary variables and
a Laurent unit in $T$. After localization at $J$ it is a single
chamber variable times a Laurent unit, by Proposition~\ref{C:chamber}
and Theorem~\ref{S:generic}, and is irreducible there by the
factorization argument of Section~\ref{C:chamber-seed}.
All but one of its factors are therefore units after localization,
hence Laurent monomials in $J$ and joint highest. Successive maximal
divided powers take the minor to $1$. By the product rule, each
joint highest factor divides $1$ already in the target for $T$,
whose units are scalar Laurent monomials in $T$. This leaves a
single variable for each minor. Distinct minors occupy distinct
positions by algebraic independence over $\Q(J)$. The first
occurrence of $i$ gives $X_i/A_i$, since the preceding reflections
fix $\varpi_i$. Thus $I\subseteq S$.

At a mutable chamber position $j$, denote the chamber coordinates
by $y_a$ and the exchange variable in the cell by $Q_j/y_j$.
This variable belongs to $D_w$, hence to this classical specialization. Write
$y'_j=P_j/y_j$ for the target exchange at the same position,
absorbing the Laurent rescalings in $T$ into $P_j$.
In the adjacent target chart the cell variable is
$(Q_j/P_j)y'_j$. Its Laurent expansion therefore gives
$P_j\mid Q_j$ in the Laurent ring on the unchanged variables.
For each mutable position outside $S$, the difference between
the largest and smallest exponents in $Q_j$ is zero. These
differences add on products, so the same holds for $P_j$,
where the difference is the absolute value of the exchange entry.
Thus every mutable neighbor of $j$ lies in $S$.

Let $i$ be a chosen frozen vertex. By Lemma~\ref{S:invariants}
we can reach a seed optimizing $i$ by mutations away from $T$.
Throughout these mutations, the chamber and exterior mutable
positions remain disconnected. Each set therefore generates the
same rational function field over the frozen variables as before.
In the optimizing seed, $M_i=A_i\partial_i(X_i/A_i)$ lies in the
chamber field. Indeed, $D_w$ is stable under $\partial_i$, and
its fraction field is generated by the chamber minors, which differ
from the variables in $S$ by units in $T$. The exterior mutable variables are algebraically
independent over that field, so none occurs in $M_i$.
Every mutable neighbor of $i$ contributes a factor to $M_i$, so
$i$ has no neighbor among the exterior mutable vertices in this
seed. This remains true on reversing the mutations, since the two
mutable components have no exchange entries between them.
\end{proof}

\begin{proof}[Proof of Theorem~\ref{S:minimal}]
If $M_R(K)$ has a coordinate realization, each $A_i$ is a Laurent
unit in $K$ by Lemma~\ref{S:necessary}. After specialization at
$v=1$, every mutable cluster factor $x$ of $A_i$ is therefore a
unit in the classical upper algebra for $K$. Its units are scalar
Laurent monomials in $K$ by the factorization argument in
Section~\ref{C:chamber-seed}. In the original algebra with all
frozen variables polynomial, $x$ is regular. In a chart containing
$K$, its monomial has nonnegative frozen exponents by regularity.
Its mutable exponents are nonnegative as well, since a negative
one introduces an exchange denominator in the adjacent chart.
Irreducibility of $x$ now makes it a single variable of $K\setminus Z$,
with coefficient one by pointedness. Thus $T\subseteq K$.

Conversely, start with $K=T$. Restrict the chart of
Lemma~\ref{S:chamber-descent} to the chamber variables $S$ and
the Laurent variables $T$. Each mutable exchange row of the chamber
is supported on $S\cup T$, so the restricted form is compatible.
The chamber admits reddening by the
\hyperref[S:freeze-reddening]{freezing argument}. Its upper algebra
embeds in $M(T)$, since exterior mutations fix the chamber chart
and its exchanges. In a shared optimizing chart the original boundary
shift is supported on $S\cup T$ as well. Thus the chamber algebra
is stable under the original operators, and joint local nilpotence
passes to it from $M(J)$.

Every chamber variable raises maximally to a Laurent unit in $T$:
its classical specialization is a chamber minor times such a unit,
and Lemma~\ref{S:invariants} identifies the quantum maximal raises
by their labels. Apply the argument of Theorem~\ref{S:generic}
to this chamber with Laurent variables $T$. Lemma~\ref{S:invariants}
identifies its joint kernel with the $T$-torus. The
\hyperref[S:delete-Laurent]{grading obtained by deleting these variables}
splits it into translates of a degree-zero algebra with scalar joint
kernel. Theorem~\ref{C:cell} identifies the latter with a Schubert
cell, and Corollary~\ref{S:frozen-orbit} gives the same Weyl group
element $w$.

It remains to include the exterior mutable components. Their seed
variables may also be incident with the Laurent variables $T$.
There are no exchange entries between the two mutable sets. By the
factorization of the distinguished lifts proved in
Section~\ref{I:mutation}, the localized basis consists of normalized
products of the component lifts. Assign each free Laurent exponent
to the chamber factor. The exterior component has no mutable
neighbor of a chosen color by Lemma~\ref{S:chamber-descent}, so
regularity at the colors involves only the chamber factor.
In a seed optimizing a chosen color, the exterior factor has exponent
zero at that color and the boundary shift is supported in $S\cup T$. Formula
\eqref{B:monomials} therefore acts on each normalized product
through its chamber factor alone. Each exterior basis element indexes one copy of the cell,
with its induced integral form, bar and common lattice.

Finally, an admissible enlargement of $T$ freezes only exterior
variables in $M_R(T)$. In any cluster, a variable of the
chamber component remains nonconstant after projecting off $T$,
whereas its joint highest elements all lie in the $T$-torus.
Freezing exterior variables changes only the indexing of the
copies. This proves sufficiency.
\end{proof}

\subsection{Restriction of the upper Kashiwara operators}

For the evaluation pairing in each root degree, the upper operators
on the restricted dual are $\widetilde e_i^{\rm up}
=(\widetilde f_i^{\rm low})^t$ and $\widetilde f_i^{\rm up}
=(\widetilde e_i^{\rm low})^t$, where the lower operators act on
$U_u^-$. On an admissible freezing, Proposition~\ref{G:prescribed}
identifies the raising arrows on the distinguished residues.
To compare the linear operators we must also identify the normalized
string vectors and decide when the ambient lowering operator
preserves the image.

On an admissible freezing we use the original boundary action.
In a shared optimizing seed write $s_i=s_i^{\rm fr}+c_i$ as in
Lemma~\ref{S:invariants}. To recover the original string through
a highest lift, translate the level-$j$ vector of its frozen
string by $-jc_i$. Indeed, Laurent translation preserves regularity
and $E_i^{(r)}=\tau_{rc_i}(E_i^{\rm fr})^{(r)}$. For an infinite
string this replaces the multiplier $L_{-s_i^{\rm fr}}$ by $L_{-s_i}$.

\begin{theorem}\label{S:shifts}
Suppose $K$ defines an admissible freezing and the distinguished
basis of $M_R(K)$ has a coordinate realization for the original positive action. For the
Weyl group element $w$ of Theorem~\ref{S:generic} and a fixed
color $i$, both linear shifts of Theorem~\ref{B:strings} agree
with the ambient upper Kashiwara operators on $M(K)$ precisely when
$\ell(\mathsf{s}_iw)<\ell(w)$.

Suppose now that the original target satisfies the frozen-product
criterion of Theorem~\ref{S:single-product}, and use its embedding
by restriction from $\widehat M_R$. Both ambient upper
operators of color $i$ preserve $M$ if and only if $S_i=L_{-s_i}$
in an optimizing seed is regular along $P$. They then preserve
$M_R$ and $\mathcal L_P$. If the regular lifts span $M_R$, their
distinguished string shifts are these restrictions.

If $S_i$ is regular along $P$ for every color, $w$ is the longest
element of the boundary Weyl group. In this case the regular
lifts span $M_R$ if and only if the regular joint highest lifts
span $\bigcap_{i\in I}\ker E_i$ in $M$ over $k$. The image is
then a direct sum of whole ambient copies of the restricted dual.
\end{theorem}
\begin{proof}
In the ambient copy through $1$, the upper lowering operator of
color $i$ sends $1$ to the functional $a_i$ of content $\alpha_i$
with $a_i(f_i)=1$. If $a_i$ belongs to the image, its preimage is a
distinguished lift of height one whose raise is $1$. In an optimizing
seed a maximal raise at height one adds $s_i$ to the degree, so this
lift is $S_i=L_{-s_i}$. Conversely, whenever $S_i$ is regular, its
height is one and its maximal raise is $E_iS_i=1$. Its image
is therefore $a_i$ in the copy through $1$.

For an admissible freezing, the cell through $1$ is the cell of $w$
by Theorems~\ref{C:cell} and~\ref{S:minimal}. Its PBW degrees include
$\alpha_i$ precisely when $\alpha_i\in\Phi(w)$, since a simple root
cannot be a sum of two positive roots. This is the left descent
condition $\ell(\mathsf{s}_iw)<\ell(w)$. When it holds, $S_i$ has
Laurent degree zero since $E_i$ preserves that degree and $E_iS_i=1$.
All $i$-strings in the cell are then infinite by
Corollary~\ref{B:string-lengths}.

For a homogeneous $i$-highest vector $y$ define
$W_j(y)=v^{j h_i(y)}yS_i^j$. Since $E_iS_i=1$ and $h_i(S_i)=d_i$,
the skew product rule gives $E_i(S_i^j)=[j]_{q_i}S_i^{j-1}$ and
$E_iW_j(y)=[j]_{q_i}W_{j-1}(y)$. At an $i$-highest lift $y=L_p$,
$h_i(y)=\Lambda(p,s_i)$, so $w_j(p)=W_j(L_p)$ by
Theorem~\ref{B:strings}.

To identify these linear shifts in every copy at once, define
$D_i x=v^{h_i(x)}xS_i$ on homogeneous elements. This operator
increases root content by $\alpha_i$. In the copy through $1$ this
follows from the twisted cell product and the content of $S_i$.
Each other copy is obtained by normalized left multiplication by an
exterior lift and a Laurent monomial (Theorem~\ref{S:minimal}), so
right multiplication by $S_i$ has the same content shift there.
Since $S_i$ represents $a_i$, we have $E_jS_i=\delta_{ij}$. The
skew product rule and
$h_i(\gamma_j)+h_j(\gamma_i)=-(\alpha_i,\alpha_j)$ give
\[ E_jD_i=u^{-(\alpha_i,\alpha_j)}D_i E_j+\delta_{ij}\operatorname{id}. \]
Left convolution by $a_i$ on the restricted dual at $u$ satisfies
the same relation, by the braided coproduct. The two operators
therefore agree under the realization. Indeed, by induction on
the input root height the difference of their values is joint highest, whereas
the ambient joint kernel has content zero and the output has
positive root height.

Now $W_j(y)=u_i^{j(j-1)/2}D_i^j y$, where $u_i=u^{d_i}$.
For $z$ in the kernel of Kashiwara's derivation $e_i'$ on $U_u^-$,
the quantum boson relation gives
$e_i'(f_i^{(r)}z)=u_i^{1-r}f_i^{(r-1)}z$ for $r>0$.
Under Kashiwara's form this derivation is adjoint to left
multiplication by $f_i$. If $\lambda$ is the component of $y$
in an ambient copy, the corresponding component of $W_j(y)$ is
$u_i^{j(j-1)/2}\lambda\circ(e_i')^j$. Its value on $f_i^{(r)}z$
is $\delta_{jr}\lambda(z)$, using $E_i y=0$ when $r>j$.
The lower string decomposition is the direct sum of
$f_i^{(r)}\ker e_i'$ \cite[\S2.2.4]{Kimura}. Thus the transposed
lower shifts send $W_j(y)$ to its two neighboring string vectors,
with zero at the top for the raising operator.

For the original target, the ambient lowering operator at $1$
again forces $S_i$ to be regular by Theorem~\ref{S:single-product}.
Conversely, if $S_i\in M_R$, it also belongs to $\widehat M_R$,
where the preceding calculation identifies both shifts on the
strings $W_j(y)$. For a homogeneous $x\in M_R$ of height $n$,
let $y=E_i^{(n)}x$. Then $y$ is integral and highest, and
$E_i^{(n)}W_n(y)=y$. Subtracting $W_n(y)$ lowers the height.
Induction therefore gives
\[ M_R=\bigoplus_{j\ge0} W_j(M_R\cap\ker E_i), \]
where $W_j$ is extended linearly over the exchange grading.
Directness follows by applying the maximal divided power. Both
shifts preserve $M_R$ and hence $M$. They preserve the common
lattice by $\mathcal L_P=M\cap\mathcal L_I$. When the regular
lifts span, Theorem~\ref{B:strings} identifies these with the
distinguished string shifts on $M$.

If all $S_i$ are regular, $w^{-1}$ sends every simple root, and
hence every positive real root, to a negative root. Thus $w$ is
the longest element. Suppose the regular joint highest lifts span
the joint kernel, and let $V'$ be the sum of the ambient copies
through their images. The image of $M$ is stable under left
multiplication by each $a_i$, which acts as $D_i$. The $a_i$
generate the restricted dual over $k$, since Kashiwara's form
identifies them with the generators $f_i$ of $U_u^-$. Hence the
image contains $V'$. Subtract from any vector in the image its
component in $V'$. If the difference were nonzero, a longest
nonvanishing raising word on it would give a nonzero joint highest
vector in the complementary copies, a contradiction. Hence the
image is $V'$, whose induced integral form has the ambient upper
global vectors as a basis. By Theorem~\ref{S:single-product}
these are exactly the regular lifts. Conversely, when the regular
lifts span, their images form a subset of the ambient upper basis.
Its span meets the joint kernel in the span of the highest vectors
in this subset.
\end{proof}
\section{Integrable extensions and their upper global bases}\label{X:extension}

Once the coroot weights are fixed, the positive boundary action
determines the lowering operators of any integrable extension.
We construct them one color at a time from the divided powers that
identify opposite weight spaces. These maps can be tested at
$v=1$. A finite test asks whether the reflection prescribed by
each coroot preserves the regular theta cone.
The upper Kashiwara operators preserve the common lattice for
every target admitting an integrable extension. When the regular
lifts span, their residues give its upper crystal basis and the
lifts themselves are the upper global vectors.

\subsection{The string-length criterion}

Choose a set $I\subseteq P$ of frozen vertices, each of which is a
sink in some seed, and let $C$ be their boundary Cartan matrix.
Write $M_R=\mathcal U_{R,P}$ and $M=k\otimes_R M_R$.
Recall that $E_j$ shifts the exchange degree by $\gamma_j=\deg Z_j$.
We prescribe adapted integral weights $\omega_i:\Gamma\to\mathbb Z$,
so that $\omega_i(\gamma_j)=c_{ij}$ for the entries of $C$
\cite[\S3.4]{Fei}. On the component of exchange degree $\lambda$,
the toral generator $K_i$ acts by $q_i^{\omega_i(\lambda)}$.
We write $U_v(\mathfrak g'(C))$ for the quantum group of the derived
Kac--Moody algebra with Cartan matrix $C$. An \emph{integrable
type-one action} means that $M$ has the indicated weight
decomposition and that all simple raising and lowering operators
are locally nilpotent. We seek lowering operators on $M$ with
these properties. When the chosen sinks have exact source partners,
Section~\ref{B:section-exact} identifies the lowering operators with
the corresponding source operators.

For a fixed color $i$, write $M_{i,m}$ for the subspace of
$i$-th coroot weight $m$. When the regular lifts span $M_R$,
the crystal arrows of color $i$ partition their labels into strings. This spanning
holds in particular when $I=P$, by Corollary~\ref{I:all-colors}.
We call $p$ a \emph{highest label} for $i$ if $E_iL_p=0$
and denote the number of edges in its string by $\ell_i(p)$. This
number can be infinite. By Corollary~\ref{B:string-lengths}, the
strings of color $i$ are infinite exactly when the distinguished
lift of degree $-s_i$ in a seed optimizing $i$ is regular.

\begin{theorem}\label{X:criterion}
The positive generators $E_i$ and the prescribed operators $K_i$ extend
to an integrable type-one $U_v(\mathfrak g'(C))$-action on $M$
if and only if the maps
\begin{equation}\label{X:opposite-weights} E_i^{(r)}:M_{i,-r}\longrightarrow M_{i,r}\qquad(i\in I,\ r\ge1) \end{equation}
are bijective. The extension is unique, its lowering operators
commute with the bar, and all simple divided powers preserve $M_R$.
Both upper Kashiwara operators of each color preserve $\mathcal L_P$,
and $A_\infty\otimes_{A_0}\mathcal L_P$ is an upper crystal lattice
of $M$.

Suppose the regular lifts span $M_R$. Then the maps
in \eqref{X:opposite-weights} are bijective if and only if
every highest label $p$ of each color $i$ satisfies
\begin{equation}\label{X:length} \ell_i(p)=\omega_i(\deg L_p). \end{equation}
Under this condition the residues of the regular $L_g$ give an
upper crystal basis of $M$. Its upper
global basis for the given integral form and bar consists of the
regular lifts. On these lifts the negative divided powers satisfy the
height and Laurent-degree bounds of \eqref{B:perfectformula} along
the reversed strings.
\end{theorem}

\begin{proof}
Fix $i$ and write $E=E_i$, $q=q_i$, and $M_m=M_{i,m}$.
An integrable type-one $U_q(\mathfrak{sl}_2)$-module is a direct
sum of finite simple strings \cite[Propositions~6.3.4 and~6.3.6]{Lusztig}. On each string
$E^{(r)}$ identifies weights $-r$ and $r$, proving necessity.

Conversely, assume \eqref{X:opposite-weights}. Injectivity at
negative weights forces every homogeneous highest vector $p$ to have
nonnegative weight $N$. Let $z$ be its unique preimage of weight
$-N$ under $E^{(N)}$. Define
\[ w_j(p)=\begin{bmatrix}N\\j\end{bmatrix}_q^{-1}E^{(N-j)}z\qquad(0\le j\le N). \]
Then $w_0(p)=p$ and $Ew_j(p)=[j]_q w_{j-1}(p)$ for $j>0$.
For $0\ne x\in M_m$ of height $r$, the vector $p=E^{(r)}x$ is
highest of weight $N=m+2r$. We have $r\le N$, since otherwise
$-m>r$ and $E^{(-m)}x=0$ contradicts injectivity. Subtracting
$w_r(p)$ lowers the height of $x$. This gives a unique expansion
in strings, since the divided power at the greatest index of a
relation recovers the highest vector at that index.

On each string define $F_iw_j(p)=[N-j]_q w_{j+1}(p)$, with
$w_{N+1}(p)=0$. These are the simple rank-one modules. In any
integrable extension $F_i^{(N)}p$ must equal $z$, so the same
formulas force $F_i$. The construction also shows that $F_i$ has
exchange degree $-\gamma_i$ and commutes with the bar, since
$E^{(N)}$ is homogeneous and commutes with the bar. A weight-preserving
map intertwining $E_i$ carries the unique preimage $z$ to the
corresponding preimage of the image of $p$. It therefore commutes
with $F_i$. This functoriality holds at $q=1$ as well.

For the relations between different colors, fix $i\ne j$ and write
$E,F,K$ for $E_i,F_i,K_i$ and $n=-c_{ij}$. For any finite
$i$-submodule $V_0$, the spaces $E_jV_0$ and $F_jV_0$ lie in a
finite $i$-submodule $V_1$. With the coproduct
$\Delta E=E\otimes1+K\otimes E$, $\Delta F=F\otimes K^{-1}+1\otimes F$,
the adjoint action on $T\in\operatorname{Hom}_k(V_0,V_1)$ of
$i$-weight $a$ is
\[ (\operatorname{ad}E)T=ET-q^aTE,\qquad (\operatorname{ad}F)T=(FT-TF)K. \]
Take $T=E_j|_{V_0}$ of weight $-n$. The successive commutators
in $(\operatorname{ad}E)^{n+1}T$ have factors
$q^{-n},q^{-n+2},\ldots,q^n$, so the positive Serre relation
makes this vector zero. In a finite simple string, a vector of
weight $-n$ killed by $E^{n+1}$ is lowest. Thus
$(\operatorname{ad}F)T=0$, proving $[F_i,E_j]=0$ on $V_0$
and hence on $M$, for every ordered pair $i\ne j$.

For $T=F_j|_{V_0}$ use the opposite coproduct. Now
$(\operatorname{ad}E)T=(ET-TE)K^{-1}=0$ by the mixed relation.
This highest vector has weight $n$, so $(\operatorname{ad}F)^{n+1}T=0$.
Since $(\operatorname{ad}F)T=FT-q^{-a}TF$ at weight $a$, this
is the negative Serre relation. The exchange degree of each $F_i$
gives the Cartan relations, completing the integrable action.
For stability of $M_R$, use its regular quantum theta basis from
Proposition~\ref{B:fiber}. For each color, the maximal divided power sends theta functions
of a fixed height to distinct theta functions with coefficient one
(Lemma~\ref{X:theta-top}). Proposition~\ref{X:integral-criterion}
below therefore gives stability under the negative divided powers
as well.

Suppose now that the regular lifts span $M_R$. By
Theorem~\ref{B:biperfect}, the maximal raise at height $j$ sends each
such lift to the highest lift of its string, with coefficient one.
The basis criterion of Proposition~\ref{X:integral-criterion} below
shows that a highest lift of weight $N$ occurs among these images
exactly for $0\le j\le N$. Thus integrability gives \eqref{X:length}.
Conversely, under \eqref{X:length} the string basis of
Theorem~\ref{B:strings} has weights $N,N-2,\ldots,-N$ on each
string, so every map in \eqref{X:opposite-weights} is bijective.

These string vectors agree with $w_j(L_p)$ by the uniqueness of
the preimage at weight $-N$. They are bar invariant by their
construction in Theorem~\ref{B:strings}. In each simple string
they satisfy
$w_j=\left[\begin{smallmatrix}N\\j\end{smallmatrix}\right]_{q_i}^{-1}F_i^{(j)}w_0$.
Kashiwara's upper operators shift these vectors by one place
\cite[\S3.3.1]{Kimura}. Each $w_j$ is congruent modulo $u\mathcal L_P$ to the
distinguished lift with the same label. Thus the upper operators
preserve the $A_0$-span of the string basis, which is $\mathcal L_P$,
and act on its residues by the distinguished string arrows. This
gives the upper crystal basis over $A_\infty$. Its global lifts are
the $L_g$ by the balanced-triple calculation of
Corollary~\ref{I:balance}, since
$\mathbb Q[v^{\pm1}]\cap A_\infty\cap\overline{A_\infty}=\mathbb Q$.

It remains to compute the negative divided powers in this basis.
Fix an $i$-string with highest label $p$ of weight $N$, and let
$g=p-js_i$ in a seed optimizing $i$. Set $r=N-j$. On the string
basis, $F_i^{(a)}w_j(p)=\left[\begin{smallmatrix}r\\a\end{smallmatrix}\right]_{q_i}w_{j+a}(p)$
for $0\le a\le r$, and higher divided powers vanish. At a fixed
coroot weight $m$, a string vector of $E_i$-height $h$ has
$F_i$-height $m+h$. Both changes between the string basis and the
distinguished basis are unitriangular in $h$, with lower coefficients
in $uA_0$, by Theorem~\ref{B:strings}. Consequently $L_g$ has
$F_i$-height $r$. For $0\le a\le r$, the term of greatest
$F_i$-height in $F_i^{(a)}L_g$ is
$\left[\begin{smallmatrix}r\\a\end{smallmatrix}\right]_{q_i}L_{p-(j+a)s_i}$.
Every other term has height less than $r-a$. Its coefficient is a sum
of products, each containing a factor in $uA_0$ and a quantum binomial
of degree at most $d_i a(r-a)$. Multiplying it by $u^{d_i a(r-a)}$
therefore gives an element of $uA_0$. Integrality of $F_i^{(a)}$
and bar invariance give the asserted strict bound on its Laurent
polynomial coefficient. For an arbitrary polynomial target, the proof of lattice stability
is completed after Theorem~\ref{X:partial-scale} by localization
at the remaining unchosen divisors of weight zero.
\end{proof}

The normalized string vectors $w_j$ need not be integral. To prove
integrality of the lowering operators, we compare $w_j$ with a basis
vector having the same maximal divided power. The quotient in this
comparison can have $R$-torsion. In the following lemma, a weight module over $R$
is a direct sum of integer weight spaces, with the toral quantum
binomials acting on each space by evaluation at its weight.

\begin{lemma}\label{X:descent}
Let $S_R$ and $Y$ be weight modules for the divided-power integral
form of $U_q(\mathfrak{sl}_2)$, with $S_R$ finitely generated and
projective over $R$.
Suppose that for each $y\in Y$ we have $F^{(a)}y=0$ for all
sufficiently large $a$. A weight-preserving $R$-linear map
$S_R\to Y$ commuting with all positive divided powers commutes
with all negative divided powers.
\end{lemma}
\begin{proof}
Finite projectivity identifies the map with a weight-zero tensor
$z\in Y\otimes_R S_R^*$ killed by $E^{(a)}$ for $a>0$, where the
dual action uses the antipode. Also $F^{(b)}z=0$ for large $b$ by
the coproduct formula, since $S_R^*$ has finite weight support and
the finitely many factors from $Y$ are killed by all sufficiently
high divided powers. For $r>0$ and sufficiently large $s$, the
commutation formula \cite[\S23.1.3]{Lusztig} gives
\[ 0=E^{(s)}F^{(s+r)}z=(-1)^s\begin{bmatrix}s+r-1\\r-1\end{bmatrix}_q F^{(r)}z. \]
Write $A_j(s)=\left[\begin{smallmatrix}s+j\\j\end{smallmatrix}\right]_q$.
By quantum Pascal,
$A_j(s)-q^jA_j(s-1)=q^{-s}A_{j-1}(s)$. Since $q$ is a unit,
annihilation by $A_{r-1}(s)$ for all large $s$ successively gives
annihilation by $A_0(s)=1$. Thus $F^{(r)}z=0$ for $r>0$, which
means that the map intertwines the negative divided powers.
\end{proof}

\begin{proposition}\label{X:integral-criterion}
Let $V$ be an integrable type-one $U_q(\mathfrak{sl}_2)$-module
over $k$, with $q=v^d$ for a positive integer $d$.
Let $V_R=\bigoplus_m V_{R,m}$ be an $R$-form graded by weight
and stable under the positive divided powers. Suppose every highest
vector $p\in V_{R,N}$ belongs to $E^{(j)}V_{R,N-2j}$ for
$0\le j\le N$. Then every negative divided power preserves $V_R$.

This hypothesis holds if $V_R$ has a weight-homogeneous basis
$\mathcal B$ with the following property. At each height $j$, the maximal divided power $E^{(j)}$ sends the
basis vectors of that height to unit multiples of distinct basis
vectors. In that case a highest basis vector of weight $N$ occurs
among these images, up to a unit, precisely for $0\le j\le N$.
\end{proposition}
\begin{proof}
For a nonzero homogeneous $x\in V_R$, let $N$ be the largest
highest weight of a simple summand occurring in $x$. Induct on $N$.
Let $W$ be the sum of the simple summands of highest weight less
than $N$. By induction $W_R=W\cap V_R$ is stable under all divided
powers, so $V/W_R$ is a module for the integral form.

Write the weight of $x$ as $N-2j$ and let $p=E^{(j)}x$, a nonzero
highest vector of weight $N$. Choose $b_t\in V_{R,N-2t}$ with
$E^{(t)}b_t=p$ for $0\le t\le N$, taking $b_j=x$.
The normalized string vectors $w_t$ of the simple module of highest
weight $N$ span a free $R$-module stable under the divided powers.
The map $w_t\mapsto b_t\bmod W_R$ intertwines the positive
divided powers. Indeed, for $a\le t$ the difference
$E^{(a)}b_t-\left[\begin{smallmatrix}t\\a\end{smallmatrix}\right]_q b_{t-a}$
is integral and killed by $E^{(t-a)}$. At weight $N-2(t-a)$ this
kernel lies in $W$, by the simple-string decomposition. For $a>t$
both images vanish. Lemma~\ref{X:descent} now gives
\[ F^{(a)}x-\begin{bmatrix}N-j\\a\end{bmatrix}_q b_{j+a}\in W_R, \]
where $b_{j+a}=0$ for $j+a>N$. This proves stability.

For the basis criterion, the height of a nonzero linear combination
is the greatest height of its summands, since their maximal images
at that height are independent. Thus the height-zero basis vectors
span the integral highest space. Over $k$ the highest space of weight
$N$ lies in $E^{(j)}(\ker E^{j+1})$ for $0\le j\le N$ and
has zero intersection with it for $j>N$, by the simple-string
decomposition. This image is spanned by the maximal images of the
height-$j$ basis vectors. Hence each highest basis vector of
weight $N$ occurs, up to an $R$-unit, among these images exactly
for $0\le j\le N$. This gives the required integral preimages.
\end{proof}

\subsection{Specialization of the divided powers}

We use the regular theta basis of Proposition~\ref{B:fiber} to
compare the opposite-weight maps at generic $v$ and at $v=1$.
For a single color $i$ we need only an integral exchange grading
with $\omega_i(\gamma_i)=2$.
In a seed where $i$ is a sink write $\Theta_g$ for the
regular quantum theta function of pointed degree $g$, and $s_i$
for the shift of $E_i$. We take all theta expansions at the same
point of the seed chamber, chosen outside every proper rational linear
hyperplane as in \cite[\S2]{TangTheta}.
Let $R_i$ be the classical boundary derivation at $i$.
For the comparison of weight spaces we use the classical upper
algebra $M^{\mathrm{cl}}=\mathbb Q\otimes_{\mathbb Z}
\mathcal U^{\mathrm{cl}}_{\mathbb Z,P}$, with the same coroot gradings
as $M$.

\begin{lemma}\label{X:theta-top}
Let $h=g_i$. For $0\le r\le h$, the expansion of $E_i^{(r)}\Theta_g$
in the regular quantum theta basis has a unique term of height $h-r$,
namely $\left[\begin{smallmatrix}h\\r\end{smallmatrix}\right]_{q_i}\Theta_{g+r s_i}$.
All other terms have smaller height. In particular
$E_i^{(h)}\Theta_g=\Theta_{g+h s_i}$, and $E_i^{(r)}\Theta_g=0$ for $r>h$.
The corresponding formulas at $v=1$ hold with $R_i^r/r!$ and
ordinary binomial coefficients.
\end{lemma}
\begin{proof}
If there are no mutable vertices, the claim is the monomial formula
\eqref{B:monomials}. Otherwise each exponent in $\Theta_g$ has the
form $g-nB_{U,V}$ with $n\ge0$. Since $B_{U,i}\ge0$, its
$i$-coordinate equals $h$ precisely when $n$ vanishes on the mutable
support of $m_i$. Freezing that support selects exactly these terms,
giving the theta function of degree $g$ in the frozen seed
\cite[Lemma~9.5]{TangTheta}. By \eqref{B:monomials}, $E_i^{(h)}$
translates these terms by $h s_i$ and kills the rest. This translation
is frozen in the new seed, so the result is its theta function of
degree $g+h s_i$ \cite[Corollary~6.7]{TangTheta}. The result lies
in the original target by Proposition~\ref{B:operators} and hence
has a finite expansion in its Laurent-polynomial theta functions
by Proposition~\ref{B:fiber}. Tang's comparison after freezing
therefore identifies it with $\Theta_{g+h s_i}$ in the original seed
\cite[Theorem~9.6]{TangTheta}. Higher divided powers vanish by
\eqref{B:monomials}.

Distinct theta labels of height $h$ have distinct maximal raises,
since their labels are translated by $h s_i$. Hence the height of a
finite theta expansion is the greatest height of its summands.
For $r\le h$, applying $E_i^{(h-r)}$ to $E_i^{(r)}\Theta_g$ gives
$\left[\begin{smallmatrix}h\\r\end{smallmatrix}\right]_{q_i}\Theta_{g+h s_i}$,
and one more raise kills the result. The only label of height $h-r$
with this maximal raise is $g+r s_i$. It therefore occurs with the
stated coefficient. Proposition~\ref{B:fiber} specializes this
expansion to the classical one, with the leading term still nonzero.
\end{proof}

\begin{theorem}\label{X:classical}
For every $i$, $r\ge1$ and coroot weight $m$, consider the maps
\[ E_i^{(r)}:M_{i,m}\longrightarrow M_{i,m+2r},\qquad R_i^r/r!:M^{\mathrm{cl}}_{i,m}\longrightarrow M^{\mathrm{cl}}_{i,m+2r}. \]
In a seed optimizing $i$, the regular source theta labels $g$ with
$g_i<r$ index a basis of the kernel. The classes of the regular
target theta functions with labels $b$ for which $b-rs_i$ is not
regular form a basis of the cokernel. These descriptions hold over
$k$ and at $v=1$ over $\mathbb Q$. In particular, injectivity and
surjectivity can each be tested after specialization.

For adapted integral coroot gradings, the positive action therefore
admits an integrable type-one quantum extension if and only if each
classical boundary derivation $R_i$ extends to an integrable
$\mathfrak{sl}_2$-action with the prescribed coroot. Equivalently,
$R_i^r$ must identify the spaces of weights $-r$ and $r$ for
every $i$ and $r\ge1$. Every such quantum extension preserves
$M_R$ under all simple divided powers, and its lowering operators
specialize to the classical ones. If the regular distinguished lifts
span $M_R$, they form the upper global basis for the distinguished
residues.
\end{theorem}
\begin{proof}
Work over either field. The regular theta labels agree by
Proposition~\ref{B:fiber}, and the leading binomials of
Lemma~\ref{X:theta-top} are nonzero. The source vectors of height
less than $r$ are killed. Among all other source labels $g$, the
greatest-height terms in their images have distinct labels
$g+r s_i$. Taking the greatest height in a finite sum proves the
kernel assertion, also over $R$. In the target, a term with label
$b$ having regular preimage $b-rs_i$ can be eliminated by subtracting
a multiple of $E_i^{(r)}\Theta_{b-rs_i}$, using the ordinary divided
power at $v=1$. Only smaller heights are introduced. This gives
spanning in the cokernel over either field. Independence follows
because every nonzero image has a greatest-height term with a
regular preimage.

The opposite-weight maps are therefore bijective at generic $v$
exactly when they are bijective at $v=1$. At either parameter the
rank-one construction of Theorem~\ref{X:criterion} identifies
these bijections with integrable rank-one actions. For adapted
weights the same theorem gives the full quantum action and its
integral divided powers. The specialized lowering operators remain
locally nilpotent and satisfy the classical rank-one relations,
so uniqueness identifies them with the prescribed classical ones.
The upper global basis assertion follows from Theorem~\ref{X:criterion}.
\end{proof}

Specialization restricts the possible boundary Cartan matrices.
It also determines when an integrable extension can pass the joint
local nilpotence test of Theorem~\ref{S:single-product}.

\begin{corollary}\label{X:indefinite}
For adapted integral weights, suppose that each classical boundary
derivation extends to an integrable $\mathfrak{sl}_2$-action with
the prescribed coroot on the classical polynomial upper algebra.
Then every connected component of the boundary Cartan matrix has
finite or affine type. On an affine component $J$ the action has
\emph{level zero}, that is, $\sum_{i\in J}a_i^\vee\omega_i=0$ on
$\Gamma$, where the $a_i^\vee$ are the coefficients of the
canonical central element in the simple coroots.

In the quantum extension the positive action is jointly locally
nilpotent if and only if every component has finite type. More
precisely, on an affine component $J$ the vectors killed by all
sufficiently long raising words in the colors of $J$ form exactly
$\bigcap_{i\in J}\ker E_i$.
\end{corollary}

\begin{proof}
The integral divided powers in Theorem~\ref{X:classical} specialize
to an action of $\mathfrak g'(C)$ with the prescribed coroots.
Multiplication intertwines the classical raising derivations and
their coroots, hence also the lowering operators by rank-one
functoriality. Thus the Chevalley generators act by derivations on
a domain whose fraction field has finite transcendence degree over
$\mathbb Q$. Each component acts nontrivially, since
$R_i(x_i)=x^{m_i}\ne0$. By \cite[Theorem~8.2 and Corollary~8.3]{Fei},
each component is finite or affine and each affine central element
acts by zero. The seed degrees generate $\Gamma$, so the
level-zero identity holds on $\Gamma$ and also in the quantum
extension with the same coroots.

Let $J$ be affine. A weight vector killed by all $E_i$ for $i\in J$
has nonnegative $i$-weight for each such $i$, by rank-one integrability.
Since all $a_i^\vee$ are positive, every one of these weights is zero.
In an integrable rank-one module a highest vector of weight zero lies
in a trivial string. Thus the joint kernel meets $E_iM$ trivially for
each $i\in J$, as one sees on each weight space. If a longest
nonvanishing raising word in $J$ on a vector had positive length,
its image would belong both to this kernel and to the image of its
last generator. This proves the assertion about vectors of finite
raising length. In particular $X_i$ has infinite raising length in
$J$, since $E_iX_i=M_i\ne0$.

If every component has finite type, each vector lies in a
finite-dimensional integrable submodule
\cite[Propositions~6.3.4 and~6.3.6]{Lusztig}. There are only finitely
many differences of weights in this submodule. Since the simple roots
are independent, the height of such a difference bounds the length of
a raising word between the two weight spaces.
\end{proof}

\subsection{A finite test on the regular theta cone}

In each seed the regular theta labels are the integral points of a
rational polyhedral cone $\mathcal C_P$. Indeed, for each $f\in P$
the Newton polytope of the dual theta function of degree $e_f^*$
gives finitely many inequalities expressing regularity at $f$
\cite[Corollary~3.2 and Lemma~3.3]{Fei}. Intersecting these
inequalities over $P$ describes the specified polynomial target
by Proposition~\ref{B:fiber}.

Fix $i$ and work in a seed where it is a sink. Write the
inequalities as $\lambda_a(g)\ge0$, and let
$H=\mathcal C_P\cap\{p:p_i=0\}$. This face parametrizes the
highest theta functions. Write $\omega_i(g)$ for the coroot
weight of degree $g$. Since $\omega_i(s_i)=2$, the map
$g\mapsto g-\omega_i(g)s_i$ is a linear reflection.

\begin{theorem}\label{X:polyhedral}
For each color $i$, the following are equivalent.
\begin{enumerate}
\item $E_i$ extends to an integrable rank-one action with coroot
$\omega_i$.
\item This reflection preserves the regular theta cone
$\mathcal C_P$.
\item The set $J=\{a:\lambda_a(s_i)>0\}$ is nonempty. Each linear
form $\lambda_a/\lambda_a(s_i)-\omega_i$, for $a\in J$, is
nonnegative on $H$, and at least one vanishes on all of $H$.
\end{enumerate}
For adapted integral weights, the full integrable extension exists
if and only if these conditions hold for every chosen color.
\end{theorem}
\begin{proof}
Write $\sigma(g)=g-\omega_i(g)s_i$. At weight $-r$ with $r>0$,
$\sigma(g)=g+r s_i$ is regular whenever $g_i\ge r$ by
Lemma~\ref{X:theta-top}. For $g_i<r$ its $i$-coordinate is negative,
so it is irregular. Thus $\sigma$ preserves regularity from weight
$-r$ to weight $r$ exactly when the map in
\eqref{X:opposite-weights} is injective. From weight $r$ to weight
$-r$, regularity of $\sigma(b)=b-rs_i$ for every regular $b$ is
exactly surjectivity by Theorem~\ref{X:classical}. Weight zero is
fixed. Since the rational cone is generated by its integral points,
Theorem~\ref{X:criterion} applied to the single color proves the
equivalence of (1) and (2).

For each $g\in\mathcal C_P$, the point $p=g+g_i s_i$ belongs
to $H$. On integral points this follows from the maximal raise
in Lemma~\ref{X:theta-top}, and on the cone it follows by rationality. The fibre through $p$
consists of $p-t s_i$ for
$0\le t\le \min_{a\in J}\lambda_a(p)/\lambda_a(s_i)$, with the
empty minimum equal to $+\infty$. Indeed, $(s_i)_i=-1$ and the
defining inequalities with $\lambda_a(s_i)>0$ give precisely the
upper bounds on $t$. Reflection sends $t$ to $\omega_i(p)-t$,
so it preserves this interval exactly when the minimum is finite
and equals $\omega_i(p)$.

Equality for every $p\in H$ means that each difference in (3) is
nonnegative and at least one vanishes at every $p$. A difference
vanishing at a relative interior point of $H$ vanishes on all
of $H$, proving (3). Conversely (3) gives the required minimum.
The full-action assertion follows from Theorem~\ref{X:criterion}.
\end{proof}

Both tests are finite. One may check
$\lambda_a(g-\omega_i(g)s_i)\ge0$ on generators of $\mathcal C_P$
for each defining inequality, since the reflection is an involution.
Alternatively, one may test the forms in the third condition on
generators of $H$. In either test one can take a lift of each
extremal ray generator modulo the lineality space, together with
both signs of a basis for that space.
When the regular distinguished lifts span, equality of frozen orders
in Proposition~\ref{I:spanning-reductions} identifies their labels
with the integral points of the same cone.

The cone also determines whether any choice of integral coroot
gradings permits an extension. The seed variables other than the
chosen sink are highest in a seed optimizing it. Their string lengths
determine all the coroot weights in that seed.

\begin{corollary}\label{X:reconstruct-weights}
For each $i$, work in a seed where $i$ is a sink. For $a\ne i$
let $N_{ia}=\sup\{r\in\mathbb N:e_a-rs_i\in\mathcal C_P\}$.
If these numbers are finite, define the integral covector $w_i$ by
\[ w_i(e_a)=N_{ia}\quad(a\ne i),\qquad w_i(e_i)=\sum_{u\in U}(m_i)_u N_{iu}-2. \]
An integrable type-one extension whose coroots are integral exchange
gradings exists if and only if all $N_{ia}$ are finite, each $w_i$
descends to an adapted coroot grading, and the reflection
$g\mapsto g-w_i(g)s_i$ preserves $\mathcal C_P$ in this seed for
every $i$. When the extension exists, its coroots and lowering
operators are unique.
\end{corollary}
\begin{proof}
For $a\ne i$, the exponent $e_a$ lies on the highest face $H$.
In an integrable extension the fibre through $e_a$ has endpoint
$\omega_i(e_a)$ in the parameter $t$ of Theorem~\ref{X:polyhedral}.
This is an integer, so $N_{ia}=\omega_i(e_a)$. The relation
$s_i=m_i-e_i$ and $\omega_i(s_i)=2$ give $w_i(e_i)$.
Since the seed degrees generate $\Gamma$, these values determine
the coroot. Conversely, if the covectors descend and are adapted,
the reflection test of Theorem~\ref{X:polyhedral} gives the
extension, and Theorem~\ref{X:criterion} gives its uniqueness.
\end{proof}

The positive-half test of Theorem~\ref{S:single-product} permits
inverting any unchosen frozen variable. For an integrable extension,
the effect of this localization is determined by the coroot weights.

\begin{corollary}\label{X:localization}
Suppose the action on the target for $P$ has an integrable
extension with adapted integral coroots $\omega_i$, and let
$T\subseteq P\setminus I$. The positive action on the target for
$P\setminus T$ admits an integrable extension with adapted integral
coroots if and only if $\omega_i(\deg X_f)=0$ for all $i\in I$
and $f\in T$. In that case the extension has the original coroots
and restricts to the given action. If an inverted variable has positive
$i$-weight, even the localized $E_i$ alone admits no integrable
type-one rank-one extension.
Spanning by the regular lifts passes from the original integral form
to the localized one. When spanning holds and the indicated weights
vanish, inclusion sends each original upper global vector to an upper
global vector of the localized target.
\end{corollary}
\begin{proof}
Every $X_f$ for $f\in T$ is highest for every chosen color, so its
original coroot weights are nonnegative. Suppose its $i$-weight is
positive. In the original integrable module, $X_f=E_i y$ for some
homogeneous $y$ on its $i$-string. Since $E_i X_f^{-1}=0$, the
product rule makes $E_i(yX_f^{-1})$ a nonzero scalar. Rescale
$yX_f^{-1}$ to $z$ with $E_i z=1$. Its exchange degree is
$-\gamma_i$, so $h_i(z)=d_i$. The product rule gives
$E_i(z^n)=[n]_{q_i}z^{n-1}$, hence $E_i^{(n)}z^n=1$ for all
$n\ge1$. In a direct sum of finite simple $i$-strings, no nonzero
vector belongs to the image of every divided power. This rules
out any integrable type-one rank-one extension.

Conversely, suppose the indicated weights vanish. The localized
target is the union of the subspaces $MD^{-1}$, where $D$ ranges
over monomials in the variables at $T$. Right multiplication by
$D^{-1}$ preserves coroot weights, and
$E_i(xD^{-1})=v^{-h_i(D)}(E_i x)D^{-1}$. Hence $E_i^{(r)}$
is an isomorphism from weight $-r$ to weight $r$ on each
$MD^{-1}$, and therefore on their union. Theorem~\ref{X:criterion}
gives the extension with the original coroots. These are unique by
Corollary~\ref{X:reconstruct-weights}. Inclusion also intertwines
the lowering operators by rank-one uniqueness.

Spanning passes to $P\setminus T$ by
Proposition~\ref{I:spanning-reductions}, since fewer frozen orders
need to be tested. Theorem~\ref{X:criterion} then identifies the
global vectors on both targets with their regular lifts.
\end{proof}

\subsection{Compatibility with multiplication}
For the target multiplication to intertwine the quantum group
action, its product rule for $E_i$ must agree with the coproduct.
We use the symmetric coproduct
$\Delta E_i=E_i\otimes K_i^{-1/2}+K_i^{1/2}\otimes E_i$, with the
same formula for $F_i$, adjoining $v^{1/2}$ if necessary.
On each homogeneous component we can compare the powers of $v$
in this formula with those in the boundary product rule.

\begin{proposition}\label{X:product}
The centered generators satisfy the symmetric product rule on
the target if and only if $2h_i=-d_i\omega_i$ on $\Gamma$ for every $i\in I$.
If the action extends to an integrable type-one module with the prescribed $K_i$, the same equalities are necessary and sufficient for
$M$ to be a module algebra for the full action with this coproduct.
\end{proposition}
\begin{proof}
The stated equality identifies $\sigma_i$ in the boundary product
rule with $K_i^{-1/2}$. For necessity, choose a seed where $i$ is a sink.
For any seed variable $X_a$ with $a\ne i$ we have $E_iX_a=0$
and $E_iX_i=M_i$. Applying the two product rules to $X_iX_a$
and cancelling $M_iX_a$ gives the required equality on $\deg X_a$.
It also holds on $\gamma_i$, since $h_i(\gamma_i)=-d_i$ and
$\omega_i(\gamma_i)=2$. The relation
$\deg X_i=\deg M_i-\gamma_i$ gives the equality on the last seed
degree. The seed degrees generate $\Gamma$.

Under these equalities, multiplication $M\otimes M\to M$ and
the unit $k\to M$ preserve weights and intertwine every $E_i$.
In an integrable extension they therefore intertwine each $F_i$
by the rank-one argument in the proof of Theorem~\ref{X:criterion}.
\end{proof}

Under integrability the two conditions of the frozen-product
criterion reduce to one. The reason is that the coroot weights of
a joint highest frozen monomial detect whether a component can be
affine.

\begin{corollary}\label{X:simultaneous}
Suppose $I\ne\varnothing$ and the positive action has an integrable
type-one extension. For the frozen monomial $F$ of
Theorem~\ref{S:single-product}, the embedding of that theorem
exists precisely when $\operatorname{top}(F)$ is a based Laurent
monomial in the unchosen frozen variables. Equivalently, $C$ has
finite type and $U_v^+F\cap\mathcal T_k(Z)\ne0$.

If a word in the raising divided powers sends $F$ to such a
monomial $A$, then for each connected component $J$ of $C$ some
$X_f$ with $f\in P\setminus I$ occurs in $A$ with positive
exponent and has positive $j$-weight for some $j\in J$.
In particular, when every polynomial frozen vertex is chosen,
the embedding is impossible.
\end{corollary}
\begin{proof}
For $i\in I$, the word must use color $i$. Indeed, $F$ is a
multiple of $X_i$ in the target. All other generators kill $X_i$
and preserve $X_i M$ by the product rule, whereas $\nu_i(A)=0$.
Fix a component $J$. Generators in distinct components commute,
so we may move all factors from $J$ to the end of the word, in
their original order of application. Let $E_j^{(r)}$ be the last
of these, with $r>0$. Then $A$ lies in its image and is killed by
$E_j$. The decomposition into simple $j$-strings gives
$\omega_j(\deg A)\ge r$. Every unchosen frozen variable is killed by
$E_j$. When it is invertible, its inverse is also killed by $E_j$,
so their opposite weights are both nonnegative and hence zero.
The remaining factors of $A$ have nonnegative exponents by regularity.
Their weights, counted with multiplicity, sum to the positive weight
of $A$. Thus at least one polynomial frozen factor has positive
exponent and positive $j$-weight.

The monomial $A$ is joint highest for $J$. On an affine component
the level-zero identity would force its $j$-weight to be zero.
Thus all components have finite type by Corollary~\ref{X:indefinite}, and the positive action is jointly
locally nilpotent. In particular, if $\operatorname{top}(F)$ is
an unchosen frozen monomial, the word computing it gives such an
$A$, and both conditions of Theorem~\ref{S:single-product} hold.
Conversely that theorem gives the stated monomial whenever the
embedding exists. The same theorem gives the equivalent condition
using $U_v^+F$, and the required polynomial factor of $A$ rules
out $P=I$.
\end{proof}
\section{Integrable actions from exact boundary pairs}\label{B:section-exact}

An integrable extension determines its lowering operators uniquely.
We ask when a specified source operator is the lowering operator
opposite a chosen sink.
Let $I$ be the chosen sinks. For each $i\in I$ choose a frozen
vertex $j=\sigma(i)$ that is a source in some seed. We keep all
variables at $I\cup\sigma(I)$ polynomial, so our target is
$M=k\otimes_R\mathcal U_{R,P}$ with $P\supseteq I\cup\sigma(I)$.
We call $(i,j)$ \emph{exact} if, in a seed where
$i$ is a sink, the source monomial $M_j^\star$ has full pointed
degree $e_i+e_j-m_i$ in the original algebra
\cite[(6.3)]{Fei}.

The classical divided-power test \cite[Proposition~6.4]{Fei} has
the following quantum form.
\begin{equation}\label{B:exact-raise} E_i M_j^\star=X_j\quad(i\ne j),\qquad E_i^{(2)}M_i^\star=M_i\quad(i=j). \end{equation}
Indeed, the cluster monomial $M_j^\star$ is a distinguished lift.
At height $r$ its maximal raise translates its degree by $rs_i$
with coefficient one (Theorem~\ref{B:biperfect}). The degree
$e_j-s_i$ has height $r=1+\delta_{ij}$ and gives the indicated
monomial. Conversely each right-hand side is highest, so subtracting
$rs_i$ recovers the degree $e_j-s_i$.

The classical operators of an exact pair form a Chevalley triple
\cite[Theorem~6.5]{Fei}. We compute the quantum commutator in a
form that also applies when some frozen variables are invertible.
For an exponent $\ell$ on those variables, write
$\tau_\ell X^a=X^{a+\ell}$ for based Laurent translation, as in
the proof of Proposition~\ref{G:localization}.

\begin{theorem}\label{B:paired}
For a family of exact pairs and any $P\supseteq I\cup\sigma(I)$,
the homomorphisms $H_i=(h_{\sigma(i)}^\star-h_i)/d_i$ are
adapted integral coroot gradings. With $F_i=E_{\sigma(i)}^\star$
and $K_ix=q_i^{H_i(\deg x)}x$, these operators give an integrable
$U_v(\mathfrak g'(C))$-action on $M$ whose simple divided powers
preserve $\mathcal U_{R,P}$.

For a single polynomial sink $i$ and polynomial source $j$, the
operators $E_i,\tau_\ell E_j^\star$ are the raising and lowering
operators of a type-one $U_{q_i}(\mathfrak{sl}_2)$-action on $M$ if and only if
$\tau_\ell M_j^\star$ has full pointed degree $e_j-s_i$ in a
seed optimizing $i$. In this case $d_i=d_j=d$, and the integral
coroot is
\begin{equation}\label{B:translated-coroot} H=(h_j^\star+\Lambda(\ell,-)-h_i)/d. \end{equation}
In particular, the two ordinary boundary operators form a
Chevalley pair precisely when $(i,j)$ is exact.
\end{theorem}
\begin{proof}
Suppose first that the degree is $e_j-s_i$. Deleting the invertible
vertices commutes with tropical mutation and gives a classical
exact pair, so $d_i=d_j=d$ by \cite[Theorem~6.5]{Fei}. Write
$E=E_i$, $D=\tau_\ell E_j^\star$ and
$h'=h_j^\star+\Lambda(\ell,-)$. In a source chart, $D$ translates
the usual source shift by $\ell$. Since $\ell_j=0$, its divided
powers are $\tau_{a\ell}(E_j^\star)^{(a)}$. Thus $D$ is locally
nilpotent and its divided powers preserve the target integral form.
In the opposite commutator \eqref{B:comm} its monomial is
$Z=\tau_\ell Z_j^\star$ and its grading is $h'$. Indeed, the
extra $\Lambda(\ell,a)$ in the grading cancels the change from
$\Lambda(n,a)$ to $\Lambda(n+\ell,a)$.

In a seed where $i$ is a sink, write $s=s_i$, so $Z_i=X^s$.
All $i$-exponents in $Z$ lie between $-1$ and $1$, and the only
term at exponent one is $X^{-s}$. Indeed $\tau_\ell M_j^\star$ is a
regular localized cluster monomial of degree $e_j-s$. Its maximal
raise is the monomial in \eqref{B:exact-raise}. By
\eqref{B:monomials} this raise translates all terms of greatest
$i$-exponent with their coefficients unchanged. Normalized division
by $X_j$ gives the assertion for $Z$.

Let $t=h'-h_i$ and $Y=E(Z)$. The shifts of $E,D$ have opposite
degrees, and each boundary grading takes the value $-d$ on its own
shift. Thus $h_i(Z)=h'(Z_i)=d$. On expanding the two centered
commutators, the terms with $x$ between $Z_i$ and $Z$ cancel,
leaving
\[ [E,D]x=\frac{v^{t(x)}xY-v^{-t(x)}Yx}{q_i-q_i^{-1}}. \]
In the seed torus \eqref{B:comm} gives
$E X^a=[a_i]_{q_i}X^{a+s}$ for every integral exponent $a$.
Thus $X^{-s}$ contributes $1$ to $Y$, and every term of $Y-1$
has $i$-exponent $-2$.
For a seed variable $X_a$ the contribution of $Y-1$ therefore
has strictly negative $i$-exponent. Regularity of $[E,D]X_a$
forces this contribution to vanish. Hence
$(Y-1)X_a=v^{2t(X_a)}X_a(Y-1)$ for every $a$. This identity
extends to Laurent monomials by multiplication and inversion.
Consequently $[E,D]x=(v^{t(x)}-v^{-t(x)})x/(q_i-q_i^{-1})$
on each homogeneous component.

At $v=1$ this commutator acts by $H$ in \eqref{B:translated-coroot}.
Both derivations kill the inverted frozen variables. Setting these
variables to one sends the two derivations to the classical sink
and source derivations of the reduced seed, as is seen in their
respective optimizing charts. On seed variables the commutator
therefore identifies $H$ with the pullback of the classical
exact-pair coroot. In particular $H$ is integral. Since the
shifts have opposite degrees and $H(\gamma_i)=2$, we have the
type-one rank-one relations. For a family of exact pairs take
$\ell=0$. Their classical coroots are adapted by
\cite[Theorem~6.8]{Fei}, so Theorem~\ref{X:criterion} assembles
the rank-one actions with their given lowering operators.

Conversely, suppose $E=E_i$ and $F=\tau_\ell E_j^\star$ are
Chevalley operators. Their monomial formulas give local nilpotence,
so the module is a sum of finite simple strings. Let $x=X_j$ and
$y=Fx=\tau_\ell M_j^\star$. We have $F^2x=0$, and $x$ has
$E$-height $\delta_{ij}$. In each string contributing to $y$, the
relevant component of $x$ is one step above the lowest vector and
at most $\delta_{ij}$ steps below the highest. Thus the $E$-height
$r$ of $y$ is the largest of these string lengths, either one or
two, and $r=2$ requires $i=j$. On the strings of length $r$, apply the divided-power formulas
to the components of $x$ one step above their lowest vectors.
This gives
\[ F^{(r)}E^{(r)}y=F^{(r)}E^{(r-1)}x\ne0. \]
A lowest component of $x$ and all shorter strings contribute zero
to both sides, by the simple-string formulas. Both $E^{(r)}y$ and
$E^{(r-1)}x$ are distinguished lifts of $F$-height $r$: the first
is a maximal raise, and the second is $X_j$ or $M_i$. The maximal source raise is injective on lifts of a fixed height
by Theorem~\ref{B:biperfect}, and Laurent translation preserves
this injectivity. Thus $E^{(r)}y=E^{(r-1)}x$. Its left side is
highest, so $r=1+\delta_{ij}$. Subtracting $rs_i$ from the
degree of this maximal raise gives $g(y)=e_j-s_i$, as in
\eqref{B:exact-raise}. Finally the commutator determines the toral
operator because the scalars $[m]_{q_i}$ are distinct for integral
$m$.
\end{proof}

For any $P\supseteq I\cup\sigma(I)$, substituting
$d_iH_i=h_{\sigma(i)}^\star-h_i$ into Proposition~\ref{X:product}
shows that the exact-pair action makes $M$ a module algebra for the
symmetric coproduct precisely when $h_{\sigma(i)}^\star=-h_i$ on
$\Gamma$ for every $i\in I$.

\begin{corollary}\label{B:upper}
For exact boundary pairs, both upper Kashiwara operators preserve
$\mathcal L_P$. In the localization keeping only $I\cup\sigma(I)$
polynomial, the residue of a regular lift follows the sink arrow
under upper raising and the paired source arrow under upper lowering.
If the regular lifts span $\mathcal U_{R,P}$, their residues give
the upper crystal basis on the original target with these arrows,
and their global lifts are the regular $L_g$.
\end{corollary}
\begin{proof}
Theorems~\ref{B:paired} and~\ref{X:criterion} give lattice stability.
In the indicated localization the regular lifts span by
Corollary~\ref{I:all-colors}. The upper raising arrow is the sink
arrow by Theorem~\ref{B:strings}. Since $F_i=E^\star_{\sigma(i)}$,
the unique term of greatest remaining height in $F_iL_g$ gives
both the upper lowering arrow and the source arrow by
Theorems~\ref{X:criterion} and~\ref{B:biperfect}. Inclusion
intertwines the upper operators and induces an injection on residues
by the intersection formula for $\mathcal L_P$. Under target
spanning, Theorem~\ref{X:criterion} identifies the global lifts.
\end{proof}

\begin{remark}\label{B:spanning-expectation}
We expect the spanning hypothesis in Corollary~\ref{B:upper} to hold
for every $P\supseteq I\cup\sigma(I)$. At the paired vertices the
orders of $L_g$ and $\vartheta_g$ agree by Corollary~\ref{I:all-colors}.
Thus only $f\in P\setminus(I\cup\sigma(I))$ require the equality
$\nu_f(L_g)=\nu_f(\vartheta_g)$ for every label $g$
(Proposition~\ref{I:spanning-reductions}). The first author's conjecture \cite[Conjecture~6.11]{FeiTropical} predicts that the
$F$-polynomials of the theta and quantum triangular basis elements
with label $g$ have the same Newton polytope. Since each frozen order
is $g_f$ plus the minimum on the corresponding polytope in direction
$w_f$, the conjecture would give the required equality.
\end{remark}

Suppose now that all frozen variables are polynomial and the exact
pairs cover them, so $F=I\cup I'$ with $I'=\sigma(I)$. For
$J\subseteq F$ let $B_J$ denote the classical seed pattern obtained
by deleting the frozen vertices outside $J$. The \emph{relative
Donaldson--Thomas square} $\Psi:B_I\to B_{I'}$ applies inverse
normalized reddening twice and relabels $i$ as $\sigma(i)$. Write
$\rho_J$ for the classical specialization setting $x_f=1$ for
$f\in F\setminus J$. Let $K$ be an opposite cluster involution of $B_I$
preserving the symmetrizing integers. There is a unique opposite
involution $\Phi$ of the full pattern satisfying $\rho_{I'}\Phi=\Psi K\rho_I$
\cite[Theorem~8.2]{Fei}. Let $\omega$ be the
permutation of $I$ induced by $K$ and $\phi$ the action of $\Phi$
on regular theta labels.

\begin{corollary}\label{B:relative-dt}
There is an $R$-linear involution $\widehat\Phi$ of the target
integral form given by $\widehat\Phi(L_g)=L_{\phi(g)}$ on the upper
global basis. It preserves the common lattice $\mathcal L_F$ and
commutes with the bar. In any seed $t$, it acts on based Laurent
expansions by $X_t^a\mapsto X_{t'}^{p(a)}$, where $t'=\Phi(t)$ and
$p$ is the induced relabeling of exponents.

The form $\Omega=\Lambda_t+p^*\Lambda_{t'}$ descends to the exchange
grading group $\Gamma$. For homogeneous $x,y$,
\[ \widehat\Phi(xy)=v^{\Omega(\deg x,\deg y)}\widehat\Phi(y)\widehat\Phi(x). \]
Thus $\widehat\Phi$ reverses multiplication if and only if
$p^*\Lambda_{t'}=-\Lambda_t$ for one, equivalently every, pair of
corresponding seeds.

For $i\in I$ and $r\ge0$, the lift conjugates $E_i^{(r)}$ to
$F_{\omega(i)}^{(r)}$. Let $S\subseteq I$ have finite Cartan type,
with longest Weyl group element $w_{0,S}$ and opposition $i\mapsto i^*$.
Write $\xi_S$ for Lusztig's involution on each connected $S$-component.
Then $\widehat\Phi(L_g)=L_{\xi_S(g)}$ for every regular label $g$
if and only if $\omega(i)=i^*$ on $S$ and $\phi(h)=w_{0,S}h$ for
every $S$-highest label $h$. For the latter equality it suffices
to test the primitive generators $(h,m)$ of the nonvertical extreme
rays of the epigraph $m\ge\delta(h)$ over the cone of $S$-highest
labels. Here $\delta(h)$ is the degree of the principal $F$-polynomial
at $h$ in a fixed seed, with $\deg Y_u=d_u$.
\end{corollary}
\begin{proof}
Corollaries~\ref{I:all-colors} and~\ref{B:upper} identify the
regular theta labels with the indices of the upper global basis.
Its vectors are bar invariant and span $\mathcal L_F$ over $A_0$,
so their permutation by $\phi$ defines the required linear involution.
Passage to the opposite of $t$ gives $L_g$ and $\vartheta_g$ the
same degree. Indeed, their codegrees in $t$ agree by their common
pointed degree in the predecessor seed
(Lemma~\ref{I:predecessor} and Proposition~\ref{B:opposite-family}).
Relabeling by $p$ takes this degree to that of $\vartheta_{\phi(g)}$
in $t'$. The exchange matrices now agree, so Theorem~\ref{I:localized}
identifies the based Laurent coefficients with those of $L_{\phi(g)}$.
Thus the monomial map in every seed agrees with $\widehat\Phi$ on an
$R$-basis.

Since $p$ takes each mutable exchange row $E_u$ to the negative of
the corresponding row in $t'$, compatibility with the same symmetrizer
gives $\Omega(E_u,-)=0$. For $X_t^a,X_t^b$ the scalar relating the two products is
$v^{\Lambda_t(a,b)-\Lambda_{t'}(p(b),p(a))}=v^{\Omega(a,b)}$.
It therefore depends only on the two exchange degrees, giving the
formula for homogeneous products. Conversely, testing reversal on pairs
of seed variables, including the frozen ones, gives $\Omega(e_j,e_k)=0$
for all $j,k$.

The identity $\rho_{I'}\Phi=\Psi K\rho_I$ gives
$\Phi(x_i)=x_{\sigma(\omega(i))}$. Thus in an $i$-sink seed the
shift and quantum binomial in \eqref{B:monomials} match those of
$F_{\omega(i)}^{(r)}$ in the corresponding source seed.
Corollary~\ref{B:upper} identifies the residue arrows with the theta
arrows, so the criterion for Lusztig's involution
and the finite highest-label test follow from \cite[Theorem~8.2 and Corollary~11.2]{Fei}.
\end{proof}

\section{Recovery of source partners}\label{X:partner-section}

Once the coroots are fixed, an integrable extension of the positive
action is unique by Theorem~\ref{X:criterion}. We recover the sources
of its lowering operators from the coroot reflections on the regular
theta cone. Each reflection singles out a frozen vertex, and the
pointed degree of its boundary monomial determines whether the pair
is exact in the reduced seed.

\subsection{The source determined by a coroot reflection}\label{X:recovery}

Write $F$ for the frozen vertices, $P$ for those kept polynomial and
$I\subseteq P$ for the chosen sinks. Let $\pi$ delete the coordinates
at $F\setminus P$ and call the resulting exchange data the
\emph{reduced seed}. Here we use its classical theta cone, whether
or not the reduced seed admits a compatible quantum form.

\begin{theorem}\label{X:partial-scale}
Suppose that, for each $i\in I$, the positive operator $E_i$ admits
an integrable rank-one extension with an integral exchange coroot
$\omega_i$. Then there is a unique reachable source $k=\kappa(i)\in P$
for which, in a seed where $i$ is a sink, $\pi g(M_k^\star)$ is a
positive multiple of $\pi(e_k-s_i)$. The map $\kappa:I\to P$ is
injective, and
\[ \pi g(M_k^\star)=r_i\,\pi(e_k-s_i),\qquad r_i=d_k/d_i\in\mathbb Z_{>0},
   \qquad \omega_i(\deg X_f)=\delta_{fk}-\delta_{fi}\quad(f\in F). \]
The vertex $k$ is a source in the successor of the chosen $i$-sink
seed. The pair is exact in the reduced seed precisely when $r_i=1$.
In particular, if $P=F$ and all frozen symmetrizing integers are
equal, every chosen sink has an exact partner in the original seed.

If the coroots are adapted, $r_i>1$ is possible only when
$\kappa(i)$ is unchosen. On every affine component $J$ of the
chosen Cartan matrix, $\kappa$ permutes $J$ and its pairs are
exact in the reduced seed.
\end{theorem}
\begin{proof}
At $f\in P$ the frozen theta order depends only on the mutable
degree and the coordinate at $f$ \cite[(5.1)]{Fei}. Thus
$\mathcal C_P$ is the inverse image of the reduced cone.
By Theorem~\ref{X:polyhedral} the reflection at $i$ preserves
$\mathcal C_P$ and its lineality space. In an $i$-sink seed,
that space lies in $g_i=0$ and contains every $e_f$ for
$f\notin P$. The $i$-coordinate of the reflection of $e_f$ is
$\omega_i(\deg X_f)$, so $\omega_i$ descends through $\pi$.
For a descended coroot the reflections commute with projection.
Preservation of the original cone is therefore equivalent to
preservation of the reduced cone.

Lemma~\ref{X:facet-recovery} identifies the frozen part of
$e_i^*+\omega_i$ in the reduced seed as $e_k^*$, giving the asserted
weights on $P$. The lemma locates $k$ in the successor, computes
$g(M_k^\star)$ in the $i$-sink seed by \eqref{X:scaled-partner},
and proves uniqueness and injectivity. Deleting frozen coordinates
commutes with tropical mutation, and the surviving frozen rows mutate
as before on the mutable columns. Hence $k$ is a source in the same
successor before deletion, with the asserted projected degree. In
the reduced seed the defining degree of an exact partner is $e_k-s_i$,
so exactness is equivalent to $r_i=1$.

Draw an arrow $i\to\kappa(i)$ for each chosen sink. These arrows
form cycles and paths since each vertex has at most one predecessor
and one successor. The positive integers $r_i$ telescope to one
around a cycle, so every pair on a cycle is exact after projection.
A path of positive length ends at an unchosen vertex.

Assume the coroots are adapted, and let $k=\kappa(i)$ be chosen.
If $k=i$, then $r_i=1$. Otherwise the degree identity, transported
to a $k$-source chart, shows that $r_i$ divides all mutable entries
of row $k$, since $\pi$ keeps those exponents of $M_k^\star$.
This divisibility persists under row mutation. Thus in a $k$-sink
chart $M_k=N^{r_i}$ for a mutable cluster monomial $N$.
Its $i$-height $-c_{ik}$ is $r_i$ times the $i$-height of $N$,
since the greatest $i$-exponent is additive on products in an
$i$-sink chart.
On the other hand $\omega_i(\deg M_k)=\omega_i(\gamma_k)+
\omega_i(\deg X_k)=c_{ik}+1$. Hence $r_i$ divides both $c_{ik}$ and
$c_{ik}+1$. This proves the assertion about multiplicities.

Finally, if $i$ belongs to an affine component $J$ and
$f=\kappa(i)$ lies outside $J$, the frozen-weight formula gives
$\sum_{j\in J}a_j^\vee\omega_j(\deg X_f)=a_i^\vee>0$.
This contradicts level zero (Corollary~\ref{X:indefinite}).
Thus $\kappa(J)\subseteq J$, and injectivity makes the restriction
a permutation.
\end{proof}

\begin{proof}[Completion of the proof of Theorem~\ref{X:criterion}]
Let $P_0=I\cup\kappa(I)$. The frozen-weight formula and
Corollary~\ref{X:localization} show that precisely the variables
in $P\setminus P_0$ can be inverted among the unchosen polynomial
variables while preserving the action. This embeds $M$ in the
integrable target for $P_0$. Every vertex of $P_0$ is a reachable
sink or source, so the regular lifts span this target by
Corollary~\ref{I:all-colors}. The spanning case of
Theorem~\ref{X:criterion} gives its upper global basis.

The inclusion intertwines the upper operators by the rank-one
string formulas. Hence both operators preserve
$\mathcal L_P=M\cap\mathcal L_{P_0}$. Over the discrete valuation
ring $A_\infty$ this intersection is weight graded and free as a
submodule of the free ambient lattice. It spans $M$ by clearing
the denominators of the finitely many coordinates of each vector.
\end{proof}

The target inherits its integral form and bar from this localization.
Its regular lifts are exactly the ambient global vectors lying in
$M$. Their negative divided powers satisfy the estimates of
Theorem~\ref{X:criterion} in the ambient basis, and every maximal
lowering is regular on $P$ by integral stability. These lifts form
an upper global basis of $M$ precisely when
$\nu_f(L_g)=\nu_f(\vartheta_g)$ for all labels $g$ and all
$f\in P\setminus P_0$ (Proposition~\ref{I:spanning-reductions}).
In particular, spanning is automatic if every unchosen polynomial
variable has nonzero weight for some chosen coroot.

\begin{corollary}\label{X:partners}\label{X:invertible}
Suppose every polynomial frozen vertex is a chosen sink. Prescribe
integral exchange gradings $\omega_i$ with $\omega_i(\gamma_i)=2$.
The length condition \eqref{X:length} holds for every color if and
only if the reduced seed has a permutation of exact source partners
and the $\omega_i$ are the pullbacks of its commutator gradings.
Under these equivalent conditions the permutation is unique, the
coroots are adapted, and the regular lifts form the upper global
basis of the original target. When
$P=F$ these pairs are exact in the original seed.

Without prescribed coroots one obtains the same criterion from
the numbers $N_{ia}$ in Corollary~\ref{X:reconstruct-weights}.
They must be finite, each resulting covector $w_i$ must descend
to the exchange grading group, and its reflection must preserve
$\mathcal C_P$. The frozen part of $e_i^*+w_i$ then identifies the
partner of $i$.
\end{corollary}
\begin{proof}
All regular lifts span by Corollary~\ref{I:all-colors}. Thus for
one color at a time, \eqref{X:length} is equivalent to rank-one
integrability by Theorem~\ref{X:criterion}. In
Theorem~\ref{X:partial-scale} the map $\kappa:P\to P$ is then a
permutation, so the cycle argument gives exact pairs on the reduced
seed without any prior assumption of adaptedness.

For each reduced exact pair, the exponent $e_k-s_i$ transported
to a $k$-source chart differs from $m_k^\star$ only in the
invertible coordinates. The translated pair of Theorem~\ref{B:paired}
therefore gives a rank-one action on the original target with
coroot pulled back from the reduced pair. Under the length
condition this is the prescribed coroot by single-color uniqueness
(Corollary~\ref{X:reconstruct-weights}). Conversely, the
pullbacks are adapted by \cite[Theorem~6.8]{Fei}, since deletion
preserves the simple shifts and the $i$-heights of the monomials
$M_j$. Thus the rank-one actions give the length condition and
the upper global basis by Theorem~\ref{X:criterion}. Finally,
Corollary~\ref{X:reconstruct-weights} applied one color at a time
gives the criterion in terms of $N_{ia}$ and $w_i$.
\end{proof}

\subsection{The lowering operator at the recovered source}\label{X:recovered-lowering}

To identify a lowering operator when its source has multiplicity
$r$, we take the $r$th root of the source monomial. Dividing the
source row and its symmetrizing integer by $r$ realizes this root
as a boundary monomial. Let $T$ be diagonal, the identity on mutable
vertices, with positive rational
entries on frozen vertices. Whenever $TB$ and $TD_V$ are integral,
\begin{equation}\label{X:frozen-rescaling} (TB)(TD_V)=T(BD_V)T,\qquad \mu_u(TB)=T\mu_u(B)\quad(u\in U). \end{equation}
The first identity preserves skew-symmetry of $BD_V$, and the
second follows
from $T_{uu}=1$ and $[tb]_+=t[b]_+$ for $t>0$. Thus the quantum
upper algebra, its distinguished lifts and the operators at the
unscaled vertices are unchanged. We may take the entries between
frozen vertices to be zero in the reference seed, since they affect
neither the mutable rows and columns nor the boundary monomials
\cite[\S2]{Fei}.

\begin{theorem}\label{X:root-lowering}
For any $I\subseteq P$, the positive action admits an integrable
type-one extension for some adapted integral coroots precisely when
each $i\in I$ has a reachable source $k\in P$ with
\[ \pi g(M_k^\star)=r_i\pi(e_k-s_i),\qquad r_i=d_k/d_i,\qquad r_i=1\ \text{if }k\in I. \]
The degree is taken in a seed optimizing $i$. Under these
conditions the map $i\mapsto k$ and the extension are unique. Transport $e_k-s_i$ to a $k$-source
chart and call its exponent $b$. The lowering operator there is
$F_i X^a=[a_k]_{q_i}X^{a+b-e_k}$. Prescribed coroots give this
extension precisely when the variables $x$ of one extended cluster
satisfy
\begin{equation}\label{X:source-weights} \omega_i(\deg x)=\operatorname{ht}_{E_k^\star}(x)-\operatorname{ht}_{E_i}(x)\qquad(i\in I), \end{equation}
where $\operatorname{ht}_D(x)=\max\{n:D^n x\ne0\}$.

In every such extension $S_i=L_{-s_i}$ has order $-1$ at $k$ and
order zero at all other vertices of $P$. For exchange-homogeneous
$x$, the lowering operator on the original target is
\[ F_i(x)=\frac{v^{h_i(x)+d_i\omega_i(\deg x)}xS_i
                   -v^{-h_i(x)-d_i\omega_i(\deg x)}S_i x}{q_i-q_i^{-1}}. \]
\end{theorem}
\begin{proof}
Necessity is Theorem~\ref{X:partial-scale}. Conversely, fix a
source $k$ satisfying the projected degree identity and write
$r=r_i$. For $k\ne i$ the ratio is the integral $k$-coordinate of the
projected degree, and for $k=i$ it is one.
In a $k$-source chart, tropical mutation commutes
with $\pi$ and with multiplication by $r$, giving
$r\pi b=\pi m_k^\star$. The monomial $M_k^\star$ involves only
mutable variables, so $b=m_k^\star/r+\ell_i$ with
$\ell_i\in\mathbb Z^{F\setminus P}$. Rescale row $k$ by $1/r$.
The row is divisible by $r$ on the mutable columns and the new
symmetrizing integer is $d_i$. The sink operator at $i$ stays the
same, since $k=i$ implies $r=1$.

The rescaled source monomial is $Q=X^{m_k^\star/r}$. Its translate
$\tau_{\ell_i}Q=X^b$ has full degree $e_k-s_i$ in the $i$-sink
seed. Hence Theorem~\ref{B:paired} supplies an integrable rank-one
action on the original target with the asserted source-chart
formula. Its coroot pulls back from the corresponding classical
exact pair in the reduced, rescaled seed. This construction for
each $i$ gives the separate rank-one extensions required by
Theorem~\ref{X:partial-scale}. In particular the candidates $k$
are unique and distinct.

Now rescale all the rows with $r_i>1$ in the reduced seed at once.
These are distinct unchosen rows. Divisibility on the mutable
columns persists under mutation, so we may rescale integrally in
one reference seed by \eqref{X:frozen-rescaling}. This leaves the
chosen sink operators and their Cartan matrix unchanged.
All the indicated pairs are exact in these classical data.
Their commutator gradings are adapted by \cite[Theorem~6.8]{Fei}.
For a given pair, rescaling other rows changes neither of its
boundary derivations, so these gradings pull back to the rank-one
coroots just constructed. Theorem~\ref{X:criterion} gives the full
action. Uniqueness follows from Corollary~\ref{X:reconstruct-weights}.

In a finite simple string the coroot weight is the lowering height
minus the raising height. For a weight vector the simple-string components have the same
weight, so the maxima of their two heights still have this
difference. Here the lowering height equals that for
$E_k^\star$. Both operators lower the $k$-coordinate by one in
the source chart and translate the surviving exponents injectively.
This proves \eqref{X:source-weights} on every nonzero homogeneous
vector. The degrees of the seed variables generate $\Gamma$.

Finally, frozen translation gives $S_i=X^{b-e_k}$ in the source
chart, with exponents $-e_k$ on $P$. Substituting
$h'=h_i+d_i\omega_i$ from \eqref{B:translated-coroot} in the
opposite commutator gives the formula for $F_i(x)$.
\end{proof}

Fix $i\in I$. In any seed $t$ write $a_f=B_{f,U}(t)$ for the
frozen row on the mutable columns and $T_t$ for the mutable tropical
Donaldson--Thomas transformation in projective coordinates.
The source test of Theorem~\ref{X:root-lowering} for $k\in P$ is
\[ \frac{a_k}{d_k}=\frac{T_t^{-2}(a_i)}{d_i},\qquad \nu_f(S_i)=-\delta_{fk}\quad(f\in P), \]
with $d_k=d_i$ if $k$ is chosen. An unchosen partner may have
$d_k/d_i>1$ and must still satisfy $\nu_k(S_i)=-1$. Inverted
frozen vertices impose no order condition.

In an $i$-sink seed $(-s_i)_U\le0$, so $S_i$ is a localized cluster
monomial in its successor (Lemma~\ref{I:predecessor}). The identities
$g_t(M_i)_U=-a_i$ and $g_t(S_i)=e_i-\operatorname{codeg}_t(M_i)$
hold in the sink seed by definition and persist under mutation
\cite[proof of Proposition~4.2]{Fei}.
Since $T_t^{-1}$ takes negative mutable degrees to mutable codegrees
\cite[\S4.2]{Fei}, these identities give
\[ T_t^{-2}(a_i)=T_t^{-1}\bigl(-g_t(S_i)_U\bigr)=\operatorname{codeg}_t(S_i)_U. \]
This codegree and the frozen row $a_k$ obey the same homogeneous
mutation rule, so one seed suffices for the comparison.
In the successor write $S_i=X^c$. Every recovered partner is a source
there by Theorem~\ref{X:partial-scale}. Conversely, the row equality
$a_k=r_i c_U\ge0$ makes $k$ a source. For such a vertex
$m_k^\star=(a_k,0)$, so the two equalities are exactly the
projected identity $r_i\pi(c+e_k)=\pi m_k^\star$.

At each unchosen polynomial vertex $f\ne k$, one must check that
$S_i$ is regular. This is equivalent to order zero since
$c_f=-\operatorname{codeg}(M_i)_f\le0$ by regularity of $M_i$ along $f$. When
Theorem~\ref{S:single-product} applies, the pole at $k$ prevents
the ambient upper lowering operator from restricting to the target
(Theorem~\ref{S:shifts}). The upper operators of the full action
preserve $\mathcal L_P$ by Theorem~\ref{X:criterion}.

Starting with exact pairs, multiply the row and symmetrizing integer
of an unchosen partner $k=\kappa(i)$ by an integer $n>1$.
This fixes $a_k/d_k$, $S_i$ and the chosen raising operators by
\eqref{X:frozen-rescaling}. The same action now has multiplicity
$n$ at $k$, so $i$ has no exact partner by Theorem~\ref{X:partial-scale}.
This applies to the base affine spaces with distinct
polynomial partners in \cite[Proposition~8.5]{Fei}, where spanning
follows from Corollary~\ref{I:all-colors}.

\begin{corollary}\label{X:restored-coefficients}
Suppose the chosen sinks admit an integrable action with adapted
integral coroots on one target, and let $P_0=I\cup\kappa(I)$.
For any set $P'\supseteq P_0$ of polynomial frozen vertices, an integrable
extension with adapted integral coroots exists precisely when
$\nu_f(S_i)\ge0$ for all $i\in I$ and $f\in P'\setminus P_0$.
The targets for each partner map thus form an interval under
inclusion, and distinct maps give disjoint intervals. Throughout
each interval the lowering operators agree. The regular lifts
give the upper global basis on the subinterval obtained by omitting
the additional divisors where the frozen orders of some $L_g$
and $\vartheta_g$ differ.
\end{corollary}
\begin{proof}
For each $i$, the pole of $S_i$ at $\kappa(i)\in P_0$ forces the
same partner on any $P'\supseteq P_0$ admitting an extension.
The normalized row and multiplicity already satisfy the source
test. At the additional vertices $f\in P'\setminus P_0$, the
required order zero is equivalent to regularity, as proved above.
The source-chart formula of Theorem~\ref{X:root-lowering} gives the
same lowering operators throughout the interval.
The global basis assertion follows from Theorem~\ref{X:criterion}
and Proposition~\ref{I:spanning-reductions}, since the regular lifts
span on $P_0$.
\end{proof}

A reduced exact pair ($r_i=1$) is exact before deletion precisely
when $\nu_f(S_i)=0$ also at every $f\in F\setminus P$, since these
orders are the coordinates of $\ell_i$. In the coefficient criterion of
\cite[Proposition~4.5]{Fei}, this says that every restored row
lies in the \emph{admissible coefficient cone} of the pair.

\section{The global basis of the quantized coordinate ring}\label{O:section}

Let $G$ be a connected, simply connected complex semisimple group
with Lie algebra $\mathfrak g$. The left and right regular actions on
$\mathcal O_q(G)$ induce two families of arrows on Kashiwara's
upper crystal basis. We recover its upper global vectors by
multiplying each regular distinguished lift by a scalar determined
by its two weights.

Write $w_0$ for the longest element of the Weyl group. In the
Berenstein--Zelevinsky seed for the double Bruhat cell $G^{w_0,w_0}$,
we keep the frozen minors polynomial. Let $\mathbf L^{\mathrm{reg}}$
consist of those lifts of Theorem~\ref{I:localized} with nonnegative
order at every frozen minor, and write $\mathbf G^{\mathrm{up}}$
for Kashiwara's upper global basis. In the coordinate ring convention
$q=v^2$, our coefficient rings are $R=\mathbb Z[q^{\pm1/2}]$
and $k=\mathbb Q(q^{1/2})$.

The cluster bar and the coordinate-ring bar differ by a scalar on
each \emph{biweight} space. We use the following convention for its
two weights.
Let $v_\delta$ be a vector of weight $\delta$ in a finite-dimensional
type-one module over $U_q(\mathfrak g)$, and let $\xi$ be a covector
vanishing outside its $\gamma$-weight space. We assign biweight
$(\gamma,\delta)$ to the matrix coefficient $z\mapsto\xi(zv_\delta)$
on $U_q(\mathfrak g)$. This convention gives the generalized minor
$\Delta_{\gamma,\delta}$ the biweight indicated by its indices.
In the bimodule notation of \cite[\S8.1]{KKKO},
$\mathrm{wt}_l=\delta$ and $\mathrm{wt}_r=\gamma$.
The exchange relations preserve the biweight, so each mutable exchange
row has biweight zero. All exponents of a pointed lift differ from
its degree by integral combinations of these rows. Each lift is
therefore homogeneous. We normalize the Weyl-invariant form
by $(\alpha_i,\alpha_i)=2d_i$ for each simple root $\alpha_i$,
where the $d_i$ are the standard symmetrizing integers. Write
$q_i=q^{d_i}$. A subscript $\gamma,\delta$ on a basis selects the
vectors of that biweight.

\begin{theorem}\label{O:comparison}
On each biweight space of $\mathcal O_q(G)$, the regular lifts and
Kashiwara's upper global basis satisfy
\[ q^{((\delta,\delta)-(\gamma,\gamma))/4}\,\mathbf L^{\mathrm{reg}}_{\gamma,\delta}=\mathbf G^{\mathrm{up}}_{\gamma,\delta}. \]
\end{theorem}

It suffices to treat a simple factor. Concatenating double words
for the simple factors gives a direct sum of seeds. The regular
lifts factor by Section~\ref{I:mutation}, Kashiwara's bases factor
under tensor products of the coordinate rings, and the exponent
in the theorem is additive.

For $G$ simple the double-word seeds are mutation equivalent, the
mutable part admits a reddening sequence, and each frozen vertex
admits an optimizing seed. Oya, Qin and Yakimov identify the upper
algebra with polynomial frozen minors as $\mathcal O_R(G)$
\cite[\S3.4 and Theorem~4.1]{OQY}. Hence the regular lifts form an
$R$-basis by Corollary~\ref{I:all-colors}. To apply the recognition
criterion we first compare the two bars and identify the boundary
operators with the raising operators for the regular actions.

\subsection{The two regular actions}

Write $J$ for the coordinate-ring bar, given by
$(Jc)(z)=\overline{c(\bar z)}$ with Lusztig's bar on
$U_q(\mathfrak g)$, and $J_{\mathrm{cl}}$ for the cluster bar. Let
$\Phi$ multiply the space of biweight $\mu=(\gamma,\delta)$ by
$q^{Q(\mu)}$, where $Q(\mu)=((\delta,\delta)-(\gamma,\gamma))/4$.
For a nonzero matrix coefficient, $\nu=\delta-\gamma$ belongs to
the root lattice. The root lattice has even norm and pairs integrally
with the weight lattice, so $2Q(\mu)=(\gamma,\nu)+(\nu,\nu)/2$ is
integral. Thus $\Phi$ multiplies each biweight space by a unit of $R$.

Berenstein and Zelevinsky identify $J$ on the seed torus with their
twisted bar \cite[Proposition~10.9(ii)]{BZ}. For a based monomial
of biweight $(\gamma,\delta)$, their symmetric grading is
$(\gamma,\gamma)-(\delta,\delta)=-4Q(\mu)$
\cite[(8.5), (10.13)]{BZ}. Their formula for the twisted bar
\cite[(6.6)--(6.7)]{BZ} thus gives $J\Phi=\Phi J_{\mathrm{cl}}$.

We next identify the raising operators for the two regular actions.
Let $e_i,f_i,K_i^{\pm1}$ be the Chevalley and Cartan generators
of $U_q(\mathfrak g)$. The two raising operators act on a matrix
coefficient $c$ by
\[ (e_i^L c)(z)=c(f_i z),\qquad (e_i^R c)(z)=c(z e_i). \]
The first raises $\gamma$ and the second raises $\delta$. For
the first action we turn the right module on the covector factor into
a left module using the anti-involution that exchanges $e_i,f_i$
and fixes $K_i$. With these conventions both factors carry highest-weight
crystals \cite[\S7.3]{KashiwaraGlobal}.

\begin{proposition}\label{O:boundary}
For each simple root $i$, the operators $\Phi^{-1}e_i^L\Phi$ and
$\Phi^{-1}e_i^R\Phi$ are positive boundary operators at reachable
frozen vertices of the polynomial target.
\end{proposition}
\begin{proof}
For each action we choose a double word in which the raising operator
kills all chamber minors except one. Write
$u_{\le k}$ and $w_{\le k}$ for the products of the
reflections at the negative and positive letters, respectively,
through position $k$ of a double reduced word for $(w_0,w_0)$.
The chamber minor $x_k$ at a letter of type $b$ has weights
$(u_{\le k}\varpi_b,w_0w_{\le k}\varpi_b)$, where $\varpi_b$
is the $b$th fundamental weight \cite[\S3.4]{OQY}. The initial frozen minor of
this type has weights $(\varpi_b,w_0\varpi_b)$.
The compatibility diagonal at a mutable vertex of type $b$ is
$2d_b$, and we assign the same value to its frozen vertices
\cite[(8.9)]{BZ}. Thus the parameter at such a vertex is
$v^{2d_b}=q_b$.

Choose $a$ so that $w_0\alpha_a=-\alpha_i$. For the right action
take a word beginning with $+a$, with initial frozen minor $x_j$
of type $a$ and first letter at position $p$.
At every letter position $k$, the reduced word for $w_{\le k}$
begins with $s_a$, so $w_{\le k}^{-1}\alpha_a$ is negative.
Every second weight therefore has nonnegative $i$-pairing except
that of the initial minor $x_j$, where the pairing is $-1$.
If an extremal vector has nonnegative $i$-weight, $e_i$ kills it.
On the length-one string at $j$, the divided-power normalization
gives $e_i v_{w_0\varpi_a}=v_{w_0s_a\varpi_a}$
\cite[\S2.2 and proof of Proposition~4.6(ii)]{OQY}. Thus
$e_i^R x_j=x_p$.

For the left action use a word ending with $-a$. Now $j$ is
the last position and $p$ the preceding occurrence of type $a$.
The weights of $x_j$ are $(w_0\varpi_a,\varpi_a)$, and those of
$x_p$ are $(w_0s_a\varpi_a,\varpi_a)$ since the intervening
letters fix $\varpi_a$. At an earlier letter position $k$, the
remaining negative letters form a reduced word for
$t=u_{\le k}^{-1}w_0$ ending in $s_a$. Thus
$u_{\le k}^{-1}\alpha_i=-t\alpha_a$ is positive. The first
weight of each earlier letter minor therefore pairs nonnegatively with
$\alpha_i^\vee$, as do the first weights of the initial minors. For
such an extremal weight $\gamma$, the weight $\gamma+\alpha_i$
is absent, so $\xi_\gamma f_i=0$. At $j$ the pairing is
$-1$, and $f_i v_{w_0s_a\varpi_a}=v_{w_0\varpi_a}$.
Pairing with the normalized extremal covectors gives $e_i^L x_j=x_p$
and $e_i^L x_k=0$ for $k\ne j$.

In both words the double-word exchange formula gives
$\widetilde b_{ju}=-\delta_{up}$ for every mutable $u$ in the
column convention of \cite[(8.7)]{BZ}. In our convention
$b_{uj}=-\widetilde b_{ju}=\delta_{up}$. Since
$b_{ju}d_u=-b_{uj}d_j$ and $d_j=d_p$, the frozen row is
$B_{j,U}=-e_p$. Thus $j$ is a sink with shift $s=e_p-e_j$.

To extend the agreement to all matrix coefficients, we compare the
product rules after conjugation by $\Phi$. Write
$H_i(\eta)=\langle\alpha_i^\vee,\eta\rangle$. Raising $\eta$ by
$\alpha_i$ increases $(\eta,\eta)/4$ by $d_i(H_i(\eta)+1)/2$.
Thus on $x$ of biweight $(\gamma,\delta)$ we have
\[ (\Phi^{-1}e_i^L\Phi)x=q_i^{(H_i(\gamma)+1)/2}e_i^L x,\qquad (\Phi^{-1}e_i^R\Phi)x=q_i^{-(H_i(\delta)+1)/2}e_i^R x. \]
At $x_j$ both scalars are one. By the coproduct formulas
$\Delta e_i=e_i\otimes1+K_i\otimes e_i$ and
$\Delta f_i=f_i\otimes K_i^{-1}+1\otimes f_i$, each conjugated
operator $\mathcal E$ satisfies the centered product rule
$\mathcal E(xy)=(\mathcal E x)v^{h(y)}y+v^{-h(x)}x(\mathcal E y)$
on homogeneous elements. Here $h=-(\alpha_i,\gamma)$ for the left
action and $h=-(\alpha_i,\delta)$ for the right.

For $k\ne j$, apply this rule to
$x_jx_k=v^{2\Lambda_{jk}}x_kx_j$. Using $\mathcal E x_j=x_p$ and
$\mathcal E x_k=0$ gives $h(x_k)=-\Lambda(s,e_k)$, as in the boundary
grading $2d_a e_j^*-\Lambda(s,-)$ of Lemma~\ref{B:grading}.
The shift $s$ raises the affected weight by $\alpha_i$, so
$h(s)=-2d_i=-2d_a$. The boundary grading takes this value on $s$
as well, since $s_j=-1$ and $\Lambda(s,s)=0$. Agreement on $e_p$
and $s$ gives agreement on $e_j$. The operators
agree on the seed variables, and the product rule determines their
values on inverses, hence on the entire seed torus.
\end{proof}

\subsection{Recognition by principal minors}

\begin{proof}[Proof of Theorem~\ref{O:comparison}]
We apply Proposition~\ref{I:recognition} to the rescaled lifts
and Kashiwara's upper global basis, both fixed by $J$. The lifts
form an $R$-basis of $\mathcal O_R(G)$, since $\Phi$ acts by units.
After adjoining $q^{1/2}=v$, the divided-power form defining
$\mathcal O_R(G)$ \cite[\S\S2.1--2.3]{OQY} is contained in the one
used in Kashiwara's balanced triple \cite[\S\S1.1, 8.1]{KKKO}.
Evaluating on the smaller form places every upper global vector in
$\mathbb Q\otimes_{\mathbb Z}\mathcal O_R(G)$. Its coordinates in
the rescaled lifts therefore belong to $\mathbb Q[v^{\pm1}]$.

Write $S$ for the $A_\infty$-span of the rescaled lifts. For the
two regular actions the coroot weights $H_i(\gamma)$ and $H_i(\delta)$
are integral exchange gradings. Each raising operator increases the
corresponding grading by two, and every matrix coefficient lies in
a finite-dimensional module for each action. By Proposition~\ref{O:boundary},
conjugating each raising operator by $\Phi^{-1}$ gives a boundary
operator. Theorem~\ref{X:criterion} gives stability of the common
lattice one color at a time. Here the regular lifts span the target,
so the rescaled lattice is $S$ and each upper lowering operator
takes a distinguished residue to another such residue or to zero.

At $q=0$ the Peter--Weyl isomorphism decomposes Kashiwara's lattice
as a direct sum of finite tensor product lattices indexed by
dominant weights $\lambda$. Under the two regular actions, the
lowering arrows generate each product crystal from its highest tensor
\cite[Proposition~8.1.2 and the definition following (8.1)]{KKKO}.
The lattice and residue basis in each factor are dual to the lower
ones \cite[\S\S1.3, 8.1]{KKKO}. Hence we may use the upper operators
of Theorem~\ref{X:criterion}, which give this same dual crystal basis
\cite[Proposition~3.25]{Kimura}. These operators commute with $J$ because
the Chevalley generators and the quantum binomials in their string
formulas are bar invariant. Applying $J$ and extending scalars along $q^{-1}=u^2$
therefore transports these arrows to infinity. The resulting lattice
is $\mathcal L(\mathbf G^{\mathrm{up}})$, since $J$ fixes the upper
global vectors. The highest tensor gives the \emph{principal minor} $p_\lambda=\Delta_{\lambda,\lambda}$,
normalized by $p_\lambda(1)=1$. This minor is an upper global
vector \cite[Lemma~9.1.1]{KKKO}, so it also represents the highest
residue at infinity.

It remains to identify this normalized vector among the rescaled
lifts. Choose a double word with all positive letters first. At
the last positive letter $k$ of type $a$ we have $u_{\le k}=1$
and $w_{\le k}\varpi_a=w_0\varpi_a$, since the subsequent positive
letters fix $\varpi_a$. Thus $x_k=p_{\varpi_a}$, and all fundamental
principal minors belong to one cluster. Their products satisfy
$p_\lambda p_\mu=p_{\lambda+\mu}$ \cite[Lemma~9.1.2]{KKKO}.
In particular the fundamental principal minors commute, so their
ordinary products are already based cluster monomials. Thus
$p_\lambda$ is a regular distinguished lift by
Lemma~\ref{I:predecessor} and Theorem~\ref{I:localized}, and
$\Phi p_\lambda=p_\lambda$ since its two weights agree.

Applying lowering words to $p_\lambda$ gives vectors in $S$ whose
residues span the $\lambda$-summand of Kashiwara's lattice modulo $u$.
Proposition~\ref{I:recognition} gives
$\mathcal L(\mathbf G^{\mathrm{up}})\subseteq S$. Since each
$p_\lambda$ represents a distinguished residue in $S$, the same
criterion gives $\mathbf G^{\mathrm{up}}=\Phi\mathbf L^{\mathrm{reg}}$.
\end{proof}

\begin{corollary}\label{O:qin}
The common triangular basis of $\mathcal O_q(G)$ agrees, up to
powers of $q^{1/2}$, with its global crystal basis in the sense
of Kashiwara. Thus Conjecture~1.6 of \cite{QinSimilarity} holds.
\end{corollary}
\begin{proof}
Theorem~\ref{I:qin-comparison} identifies the distinguished lifts with
the localized common triangular basis. Qin identifies the common
triangular bases of the decorated seed and the Berenstein--Zelevinsky
seed of the double Bruhat cell by an anti-isomorphism that matches
the frozen minors \cite[proof of Theorem~7.4]{QinSimilarity}.
His basis for the partial compactification consists of the elements
regular along these minors \cite[Definition~2.22]{QinSimilarity}.
These are precisely the regular lifts in $\mathcal O_q(G)$.
Theorem~\ref{O:comparison}
identifies their normalizations on every biweight space.
\end{proof}

\appendix
\section{Facet normals of the regular cone}\label{X:facets}

To recover source partners in Section~\ref{X:partner-section}, we
apply the following lemma after deleting the inverted frozen variables.
Let $\mathcal C$ be the regular theta cone of a classical
skew-symmetrizable seed with mutable vertices $U$ and polynomial
frozen vertices $F$. Assume the mutable part admits a reddening sequence.

\begin{lemma}\label{X:facet-recovery}
Let $I\subseteq F$ consist of reachable sinks. For each $i\in I$,
choose a seed optimizing $i$, write $M_i=x^{m_i}$ and $s_i=m_i-e_i$,
and let $\omega_i$ be an integral exchange grading with
$\omega_i(s_i)=2$. Suppose the reflection
$g\mapsto g-\omega_i(g)s_i$ preserves $\mathcal C$ in this seed.
Then the frozen part of $e_i^*+\omega_i$ is $e_k^*$ for a unique
$k=\kappa(i)$. This vertex is a source in the successor. With pointed
degrees taken in the chosen $i$-sink seed, we have
\begin{equation}\label{X:scaled-partner} g(M_k^\star)=r_i(e_k-s_i),\qquad r_i=d_k/d_i\in\mathbb Z_{>0}. \end{equation}
If a reachable source $h$ satisfies
$g(M_h^\star)\in\mathbb R_{>0}(e_h-s_i)$, then $h=k$.
The map $\kappa:I\to F$ is injective.
\end{lemma}

\begin{proof}
For $f\in F$ let $\mathcal P_f\subseteq\mathbb R_{\ge0}^U$ be its
\emph{boundary polytope} and $H_f(w)=\max_{a\in\mathcal P_f}wa$
its support function. We regard $a$ as a column and extend it by
zero on $F$. The inequalities defining $\mathcal C$ have normals
$n_{f,a}=e_f^*+a$, one for each vertex $a$ of $\mathcal P_f$,
including the origin \cite[Corollary~3.2 and Lemma~3.3]{Fei}.
The vertices are integral, so the $f$-coordinate one makes each
normal primitive. To see that
$n_{f,a}$ defines a facet, choose $g_U$ with unique minimum at $a$
on $\mathcal P_f$, take $g_f=-g_Ua$, and choose the other frozen
coordinates so their inequalities are strict. A mutable cluster
monomial $Q$ has order zero at every frozen divisor, so
$g(Q)_f=H_f(-g(Q)_U)$ \cite[Lemma~4.1]{Fei}.

Fix $i$, write $m=(m_i)_U=-B_{i,U}$ and $\omega=\omega_i$.
The normal $e_i^*$ takes value $-1$ on $s_i$. For any other
facet with normal $n$, choose $g$ in its relative interior. Then
$g_i>0$ and $g+g_i s_i\in\mathcal C$ by \cite[Theorem~3.4]{Fei},
so $n(s_i)\ge0$. Since $\omega(s_i)=2$, the reflection is an integral
involution taking $s_i$ to $-s_i$. Its dual sends $n$ to
$n-n(s_i)\omega$ and permutes the primitive inward normals. Thus
$e_i^*+\omega=e_k^*+c$ is the unique normal of value $1$ on $s_i$,
and all normals other than this one and $e_i^*$ vanish on $s_i$. The frozen part identifies $k$ independently of the
optimizing seed, since the degrees of the frozen variables do not
change under mutation.

Now $n_{f,a}(s_i)=ma-\delta_{fi}$. In particular $mc=1+\delta_{ik}$
and $c$ uniquely maximizes $ma$ on $\mathcal P_k$. When $k\ne i$,
we also need a nonzero vertex of $\mathcal P_i$. Choose $u$ with
$m_u>0$. The exchange relation gives $g(x'_u)_i=[b_{ui}]_+>0$
for the variable obtained by mutating at $u$, so the frozen order
formula implies $\mathcal P_i\ne\{0\}$. Every nonzero vertex there
has $ma=1$, since its normal vanishes on $s_i$. We have thus obtained
$H_f(m)=\delta_{fi}+\delta_{fk}$ for $f\in F$.

To identify the source monomial, write $b_f=B_{f,U}$ and
$D_U=\operatorname{diag}(d_u:u\in U)$. Let $p$ be the mutable degree
with mutable codegree $b_k$, and let $\eta$ be the greatest exponent
of its principal $F$-polynomial $\mathcal F$. Then $p=b_k+\eta B_\mu$
\cite[Proposition~1.6 and \S4.2]{Fei}. By the frozen order formula
and tropical duality \cite[(5.1) and Theorem~4.3]{Fei},
\[ d_kH_k(w)=\max_{\nu\in\operatorname{supp}\mathcal F}wD_U(\eta-\nu)^{\mathsf T}. \]
Since $\nu\ge0$ and $0$ occurs, $D_U\eta^{\mathsf T}/d_k$ is the
greatest vertex of $\mathcal P_k$. It maximizes $ma$ because
$m\ge0$, so it equals $c$. The exchange grading
$\omega=e_k^*-e_i^*+c$ annihilates the mutable exchange rows. Hence
\[ 0=B_{U,k}-B_{U,i}+B_\mu c=-D_U\bigl(p/d_k-b_i/d_i\bigr)^{\mathsf T}, \]
where we used $BD_V=-D_VB^{\mathsf T}$. Thus the integral mutable
degree $p$ is $-r_i m$ for $r_i=d_k/d_i$.

Since $p\le0$ and the mutable degrees of the successor variables
are the negative coordinate vectors up to permutation, there is a
mutable cluster monomial $Q$ in the successor with mutable degree
$p$ in our seed. Its mutable codegree here is $b_k$ by the choice
of $p$.
Mutable codegrees and frozen rows obey the same tropical mutation
rule \cite[proofs of Lemma~1.5 and Proposition~4.2]{Fei}. In the
successor, $B_{k,U}$ is therefore the exponent of $Q$. This makes
$k$ a source there and $Q=M_k^\star$. Its frozen coordinates are
$g(Q)_f=H_f(r_i m)=r_i(\delta_{fi}+\delta_{fk})$, proving
\eqref{X:scaled-partner}. For $k\ne i$ the $k$-coordinate gives
$r_i\in\mathbb Z_{>0}$, and for $k=i$ the ratio is one.

If a reachable source $h$ satisfies $g(M_h^\star)=t(e_h-s_i)$
with $t>0$, its frozen coordinates are $H_f(tm)=t(\delta_{fi}+
\delta_{fk})$, so $h=k$.

For injectivity, fix $k$ in the image of $\kappa$. Negating the
pointed degree of $M_k^\star$ in every seed gives a tropical point
of the opposite pattern \cite[proof of Proposition~4.2]{Fei}.
For each $i$ mapping to $k$, \eqref{X:scaled-partner} identifies
the theta function of this point in the $i$-sink seed with the
localized cluster monomial $(M_i/(x_ix_k))^{r_i}$. Its frozen orders
are $-r_i(\delta_{fi}+\delta_{fk})$, since $M_i$ has order zero at
every frozen divisor. The set of frozen poles is therefore $\{i,k\}$.
It depends only on the tropical point determined by $k$, so it
recovers $i$ even when $i=k$.
\end{proof}

\section*{Declaration of generative AI and AI-assisted technologies}
The authors used ChatGPT and Codex (OpenAI) to explore arguments, check proofs
and references, and edit the manuscript. The authors take responsibility for its content.

\bibliographystyle{amsalpha}

\end{document}